\documentclass[11pt]{article}

\usepackage[table,xcdraw]{xcolor}
\usepackage{amsmath}
\usepackage{amsfonts}
\usepackage{amsthm}
\usepackage[english]{babel} 
\usepackage{dsfont}
\usepackage{amssymb}
\usepackage{graphicx}        
\usepackage{cite}                            
\usepackage{epstopdf}
\usepackage{epsfig}
\usepackage[table]{xcolor}
\usepackage{subfloat}
\usepackage{float}
\usepackage[thinlines]{easytable}
\usepackage{array}
\usepackage{mathtools}
\usepackage{mwe} 
\usepackage{caption}
\usepackage[skip=1pt]{subcaption}
\usepackage{cite}

\allowdisplaybreaks
\numberwithin{equation}{section}

\usepackage[linesnumbered,algoruled,boxed,lined]{algorithm2e}

\def\RR{\mathbb R}
\def\EE{\mathcal F}

\def\e{\varepsilon}

\def\diag{{\rm diag}}
\def\sign{{\rm sgn}}
\def\argmin{{\rm arg}\!\min}
\def\be{\begin{equation}}
\def\ee{\end{equation}}
\def\bea{\begin{eqnarray}}
\def\eea{\end{eqnarray}}

\def\p{{y}}
\def\g{{\bar{y}}}

\def\P{{Y}}
\def\G{{\bar{Y}}}
\def\Q{{\bar{X}}}

\def\iw{m}

\newcommand{\mc}[1]{\mathcal{#1}}

\title{From particle swarm optimization to consensus based optimization: stochastic modeling and mean-field limit}

\author{Sara Grassi\footnote{Department of Mathematics \& Computer Science, University of Ferrara, Via Machiavelli 30, Ferrara, 44121, Italy (sara.grassi{@}unife.it).} \and Lorenzo Pareschi\footnote{Department of Mathematics \& Computer Science, University of Ferrara, Via Machiavelli 30, Ferrara, 44121, Italy (lorenzo.pareschi{@}unife.it).}}
\begin{document}
\maketitle

\begin{abstract}
In this paper we consider a continuous description based on stochastic differential equations of the popular particle swarm optimization (PSO) process for solving global optimization problems and derive in the large particle limit the corresponding mean-field approximation based on Vlasov-Fokker-Planck-type equations. The disadvantage of memory effects induced by the need to store the local best position is overcome by the introduction of an additional differential equation describing the evolution of the local best. A regularization process for the global best permits to formally derive the respective mean-field description. Subsequently, in the small inertia limit, we compute the related macroscopic hydrodynamic equations that clarify the link with the recently introduced consensus based optimization (CBO) methods. Several numerical examples illustrate the mean field process, the small inertia limit and the potential of this general class of global optimization methods.  

\end{abstract}


\section{Introduction}
Optimization by swarms of particles (Particle Swarm Optimization or PSO) was initially proposed to model the intelligent behavior of flocks of birds or fish schools \cite{kennedy1995particle,kennedy1997particle,shieber98}. 
As a particle-based stochastic optimization algorithm, the PSO has attracted a great deal of attention from the scientific community, producing a huge number of variants of the standard algorithm \cite{kennedy2010particle,poli2007particle,shieber98, hassan04}. Today, similarly to other metaheuristic methods \cite{Aarts:1989:SAB:61990,Back:1997:HEC:548530,Blum:2003:MCO:937503.937505,Gendreau:2010:HM:1941310}, PSO is recognized as an efficient method for solving complex optimization problems and is currently implemented in several programming languages. Among popular metaheuristic methods we recall evolutionary programming \cite{Fogel:2006:ECT:1202305},  Metropolis-Hastings sampling algorithm \cite{hastings70}, genetic algorithms \cite{Holland:1992:ANA:531075}, ant colony optimization (ACO) \cite{dorigo2005ant} and simulated annealing (SA) \cite{holley1988simulated,kirkpatrick1983optimization}. Despite its apparent simplicity, the PSO presents formidable challenges to those interested in understanding swarm intelligence through theoretical analyses. So, to date a fully comprehensive mathematical theory for particle swarm optimization is still not available.

Recently, a new class of particle based methods for global optimization based on consensus (Consensus Based Optimization or CBO) has been introduced \cite{pinnau2017consensus,carrillo2018analytical,carrillo2019consensus,fhps20-1,fhps20-2,TW,Totzeck2018ANC}. These methods are intrinsically simpler than PSO methods and have been inspired by consensus like dynamics typical of social interactions like opinion formations and wealth exchanges \cite{partos13}.
In contrast to classic metaheuristic methods, for which it is quite difficult to provide rigorous convergence to global minimizers (especially for those methods that combine instantaneous decisions with memory mechanisms), CBO methods, thanks to the instantaneous nature of the dynamics permit to exploit mean-field techniques to prove global convergence for a large class of optimization problems \cite{carrillo2018analytical,carrillo2019consensus,fhps20-2}. Despite their simplicity CBO methods seem to be powerful and robust enough to tackle many interesting high dimensional non-convex optimization problems of interest in machine learning \cite{carrillo2019consensus,fhps20-2}. 

Motivated by these results, in this manuscript we take a first step towards building a robust mathematical theory for PSO methods based on a continuous description of their dynamics. One of the main difficulties is the introduction by PSO methods, and other metaheuristic algorithms, of memory mechanisms that make their interpretation in terms of differential equations particularly challenging. To this end, the description of the PSO model through a system of stochastic differential equations is based on the introduction of an additional state variable that takes into account the memory of the single particle. In this way, the resulting time continuous PSO dynamics is defined by a system of stochastic differential equations that uses both a global best search and a local best search. 

Adopting the same regularization process for the global best as in CBO methods \cite{pinnau2017consensus, carrillo2018analytical}, it is then possible to pass to the mean field limit and derive, at a formal level, the corresponding Vlasov-Fokker-Planck equation that characterizes the behavior of the system in the limit of a large number of particles \cite{Bellomo17,Golse, Jab2, Snitzman}.
 Thanks to the new mathematical formalism based on mean field equations, it is then possible to study the behavior of PSO methods in the limit of small inertia, similarly to what done in other contexts for nonlinear Vlasov-Fokker-Planck type systems \cite{AS, DLP}. In particular, we show how in this limit the PSO dynamics is described by simplified hydrodynamic models that correspond to a generalization of CBO models including also memory effects and local best search. 
 
These results are subsequently validated by comparing a direct numerical solution of the stochastic particle systems with a finite volume discretization of the mean-field system \cite{Platen, DP}. Let us emphasize that even if, as a side results of our analysis, we will derive novel global optimization methods based on stochastic differential equations (SDEs) and mean-field partial differential equations (PDEs), it is beyond the scopes of the present manuscript to perform an extensive testing of the various methods performances and to discuss the practical algorithmic enhancements that can be adopted to increase the success rate, like for example the use of random batch methods \cite{AlPa, carrillo2019consensus, JLJ}, particle reduction techniques \cite{fhps20-2} and parameters adaptivity \cite{poli2007particle}. In contrast, our numerical test cases, will address the numerical validation of the mean field process and the small inertia limit, and the role of the various parameters involved in solving high dimensional global optimization problems for some prototype test functions.

The rest of the manuscript is organized as follows. In Section 2 we introduce the various discrete PSO models and derive the corresponding representations as SDEs using a suitable time continuous approximation of the memory process. Thanks to a regularization of the global best, in Section 3, we formally pass to the large particle limit and compute the respective Vlasov-Fokker-Planck equations describing the mean-field dynamic. Section 4 is then devoted to the study of the small inertia limit for the mean-field system that allows to recover a generalized CBO model as the corresponding hydrodynamic limit. Next, in Section 4, we report several numerical examples, validating the mean-field approximation, the small inertia limit and testing the performances of the minimizers against some prototype functions in high dimension. Finally, some conclusions and future research directions are reported in the last section.

\section{Stochastic differential models for particle swarm optimization}
In the sequel we consider the following optimization problem
\begin{equation}\label{typrob}
x^\ast \in \argmin\limits_{x\in \RR}\EE(x)\,,
\end{equation}
where $\EE(x):\mathbb R^{d} \to \mathbb R$ is a given continuous cost function, which we wish to minimize. In particular, both statistical estimation and machine learning consider the problem of minimizing an objective function in the form of a sum  \cite{bishop06,fumio00,vapnik91}
\begin{equation}\label{machlear}
\EE(x)=\frac1{n}\sum_{i=1}^n \EE_i(x).
\end{equation}
The PSO algorithm solves the above minimization problem  by starting from a population of candidate solutions, represented by particles, and moving these particles in the search space according to simple mathematical relationships on particle position and speed.  The movement of each particle is influenced by its best known local position, but it is also driven to the best known position in the search space, which is updated when the particles find better positions.  

\subsection{The original PSO method}

The method is based on introducing $N$ particles with position $x_i\in {\RR}^d$ and speed $v_i\in {\RR}^d$, $i=1,\ldots,N$. The particle positions and velocities, starting with an initial $x_i^0$ and $v_i^0$ assigned, are updated according to the following rule for $n\geq 0$
\begin{equation}
\begin{split}
x^{n+1}_i &= x_i^n + v_i^{n+1},\\
v^{n+1}_i &= v_i^n + c_1 R_1^n \left(\p_i^n-x_i^n\right)+c_2 R_2^n \left(\g^n-x_i^n\right),
\end{split}
\label{eq:pso}
\end{equation}
where $c_1, c_2 \in {\RR}$ are the {\em acceleration coefficients}, $\p_i^n$ is the {\em local best} position found by the $i$ particle up to that iteration, and $\g^n$ is the {\em global best} position found among all the particles up to that iteration. The terms $R_1^n$ and $R_2^n$ denote two $d$-dimensional diagonal matrices with random numbers uniformly distributed in $[0,1]$ on their diagonals. These numbers are generated at each iteration and for each particle.
Typically, the values of $x_i$ and $v_i$ are restricted within a specific search domain $X=[X_{min},X_{max}]^d$ and $V=[-V_{max},V_{max}]^d$.

By replacing the second equation in the first one and using the fact that in the previous step $x_i^n=x_i^{n-1}+v_i^{n}$, we get a model based on a single two-level recursive equation for the particle positions
\begin{equation}
x^{n+1}_i = 2x_i^n -x_i^{n-1} + c_1 R_1^n \left(\p_i^n-x_i^n\right)+c_2 R_2^n \left(\g^n-x_i^n\right).
\end{equation}
There are several ways to define the functions $\p_i^n$ and $\g^n$. In the original PSO method, these are defined by the following relationships
\begin{eqnarray}
\nonumber
\p_i^0&=&x_i^0,\\ 
\label{eq:psoEvolution}
\p_i^{n+1} &=& \left\{
\begin{array}{lcl}
\p_i^n  & \hbox{if}   & \EE(x_i^{n+1}) \geq \EE(x_i^n),  \\ 
x_i^{n+1}  & \hbox{if}   & \EE(x_i^{n+1}) < \EE(x_i^n);  
\end{array}
\right.\\ 
\nonumber
\g^0&=&\hbox{argmin}\{\EE(x_1^0),\EE(x_2^0),\ldots,\EE(x_N^0)\},\\ 
\nonumber
\g^{n+1} &=& \hbox{argmin}\{\EE(x_1^{n+1}),\EE(x_2^{n+1}),\ldots,\EE(x_N^{n+1}),\EE(\g^n)\}. 
\end{eqnarray}

\subsection{The stochastic differential PSO system}
In order to derive a time continuous version of the discrete PSO method \eqref{eq:pso}, we rewrite it in the form
\begin{equation}
\begin{split}
x^{n+1}_i &= x_i^n + {v_i^{n+1}},\\
v^{n+1}_i &= v_i^n + \frac{c_1}{2}\left(\p_i^n-x_i^n\right)+\frac{c_2}{2}\left(\g^n-x_i^n\right)+\frac{c_1}{2}\widetilde R_1\left(\p_i^n-x_i^n\right)+\frac{c_2}{2}\widetilde R_2\left(\g^n-x_i^n\right),
\end{split}
\label{eq:pso2}
\end{equation} 
where $\widetilde R_k=(2R_k-1)$, $k=1,2$. We can interpret \eqref{eq:pso2} as a semi-implicit time discretization method for SDEs of time stepping $\Delta t=1$ where the implicit Euler scheme has been used for the first equation and the Euler-Maruyama method is used for the second one. Note that, the particular distribution of the random noise will not change the corresponding stochastic differential system provided the noise has the same mean value and variance. In the case of the PSO model \eqref{eq:pso2}, since the random terms are uniformly distributed in $[-1,1]$, the mean value is $0$ and the corresponding variance $1/3$. 

We can then write its time continuous formulation as a second order system of SDEs in It\^o form
\begin{equation}
\begin{split}
dX^i_t &= V^i_t dt,\\
dV^i_t &= \lambda_1\left(\P_t^i-X^i_t\right)dt+\lambda_2\left(\G_t-X^i_t\right)dt+\sigma_1 D(\P_t^i-X^i_t)dB^{1,i}_t+\sigma_2 D(\G_t-X^i_t)dB^{2,i}_t,
\end{split}
\label{eq:psoc}
\end{equation} 
with 
\begin{equation}
\lambda_k=\frac{c_k}{2},\quad \sigma_k = \frac{c_k}{2\sqrt{3}},\quad k=1,2 
\label{eq:param}
\end{equation}
and
\begin{equation}
D(X_t)=\diag\left\{(X_t)_1,(X_t)_2,\dots,(X_t)_d\right\},
\end{equation}
a $d$-dimensional diagonal matrix. 
In \eqref{eq:psoc} the vectors 
$B^k_t=\left((B^k_t)_1,(B_t^k)_2,\dots,(B_t^k)_d\right)^T$, $k=1,2$ denote
$d$ independent 1-dimensional Brownian motions and depend on the $i$-th particle. 

One critical aspect is the definition of the best positions $\P_t^i$ and $\G_t$ which in the PSO method make use of the past history of the particles. In \cite{hassan04} the authors observed that the local best can be rewritten as
\[
\p_i^{n+1} = \p_i^n + \frac12\left(x_i^{n+1}-\p_i^n\right)S(x_i^{n+1},\p_i^n),
\]
where
\be
S(x,y)=\left(1+\sign\left(\EE(y)-\EE(x)\right)\right).
\ee
Therefore, for a positive constant $\nu$ we can approximate the above equation with the following differential system
\begin{equation}
d\P_t^i = \nu \left(X^i_t-\P^i_t\right)S(X^i_t,\P^i_t)dt,
\label{eq:lbest}
\end{equation}
with $Y^i_0=X^i_0$
and consequently define
\begin{equation}
\G_t = \hbox{argmin}\left\{\EE(\P^1_t),\EE(\P^2_t),\ldots,\EE(\P^N_t)\right\}.
\end{equation}
Note that, equation \eqref{eq:lbest} does not describe 
the evolution of the local best, but rather a time continuous approximation of its evolution. 
\subsection{Stochastic differential PSO model with inertia}
To optimize the search algorithm, the value $c_k= 2$, $k=1,2$ was adopted in early PSO research. This value, which corresponds to $\lambda_k=1$ and $\sigma_k=1/\sqrt{3}$, $k=1,2$ in the SDEs form, however, leads to unstable dynamics with particle speed increase without control. The use of hard bounds on velocity in $[-V_{\max},V_{\max}]^d$ is one way to control the velocities. However, the value
of $V_{\max}$ is problem-specific and difficult to determine. For this reason, the PSO model has been  considered with a modified term which reads as \cite{shieber98}
\begin{equation}
\begin{split}
x^{n+1}_i &= x_i^n + v_i^{n+1},\\
v^{n+1}_i &= \iw  v_i^n + c_1 R_1^n \left(\p_i^n-x_i^n\right)+c_2 R_2^n \left(\g^n-x_i^n\right),
\end{split}
\label{eq:psoi}
\end{equation}
where $\iw \in (0,1]$ is the inertia weight. The above system can be rewritten as
\begin{equation}
\begin{split}
x^{n+1}_i &= x_i^n + v_i^{n+1},\\
\iw  v^{n+1}_i &= \iw  v_i^n - (1-\iw ) v^{n+1} + c_1 R_1^n \left(\p_i^n-x_i^n\right)+c_2 R_2^n \left(\g^n-x_i^n\right).
\end{split}
\label{eq:psoi2}
\end{equation}
In this case, we can interpret the second equation as a semi-implicit Euler-Maruyama method, that is implicit in $v_i$ and explicit in $x_i$, hence the corresponding SDEs system reads
\begin{eqnarray}
\nonumber
dX^i_t &=& V^i_t dt,\\[-.5cm]
\label{eq:psoci}
\\
\nonumber
\iw  dV^i_t &=& -\gamma V^i_t dt +\lambda_1\left(\P_t^i-X^i_t\right)dt +\lambda_2\left(\G_t-X^i_t\right) dt+\sigma_1 D(\P_t^i-X^i_t)dB^{1,i}_t+\sigma_2 D(\G_t-X^i_t)dB^{2,i}_t,
\end{eqnarray} 
where $\gamma=(1-\iw ) \geq 0$. Thus, the constant $\gamma$ acts effectively as a friction coefficient, and can be related to the fluidity of the medium in which particles move. System \eqref{eq:psoci} is reminescent of other second order stochastic particle system with inertia \cite{AS, DLP}. However, note that here, the inertia weight $\iw $ and the friction coefficient $\gamma$ are not independent.

In practice, in the PSO method \eqref{eq:psoi} the parameter $\gamma$ is often initially set to some low value, which corresponds to a system where particles move in a low viscosity medium and perform extensive exploration, and gradually increased to a higher value closer to one, where the system is more dissipative and would more easily concentrate into local minima.
Most PSO approaches, nowadays, are based on \eqref{eq:psoi} (or some variant) which is usually referred to as canonical PSO method to distinguish it from the original PSO method \eqref{eq:pso} (see \cite{poli2007particle}). Similarly we will refer to \eqref{eq:psoc}-\eqref{eq:lbest} as the original stochastic differential PSO (SD-PSO) system and to \eqref{eq:psoci}-\eqref{eq:lbest} as the canonical SD-PSO system.  
We emphasize that these stochastic systems, when discretized according to the methods described above (namely implicit in $V^i_t$ and explicit in $X^i_t$) and with the choice $\Delta t=1$ correspond to the original discrete PSO methods.


\section{Mean-field description of particle swarm optimization}
In this section we introduce a modified version of the canonical stochastic differential PSO system for which we can formally compute its mean field limit. We first consider the case in absence of memory effects and then we extend the results to the general case.

\subsection{Regularized PSO dynamic without memory effects}
To simplify the mathematical description, let us consider a PSO approach where the dynamic is instantaneous without memory of the local best positions and the global best has been regularized as in \cite{pinnau2017consensus}. The corresponding second order system of SDEs takes the form
\begin{eqnarray}
\nonumber
dX^i_t &=& V^i_t dt,\\[-.5cm]
\label{eq:psociri}
\\
\nonumber
\iw  dV^i_t &=& -\gamma V^i_t dt+\lambda\left(\Q^{\alpha}_t-X^i_t\right)dt+\sigma D(\Q^{\alpha}_t-X^i_t)dB^i_t,
\end{eqnarray} 
where $\Q^{\alpha}_t$ is the weighted average
\begin{equation}\label{ValphaE}
\Q^{\alpha}_t =\frac1{N^\alpha}{\sum_{i=1}^{N} X_t^i\omega_\alpha(X_t^i)},\quad N^{\alpha} = {\sum_{i=1}^{N}\omega_\alpha(X_t^i)},
\quad  \omega_\alpha(X_t):=e^{-\alpha\EE(X_t)}\,.
\end{equation}
The choice of the weight function $\omega_\alpha$ in \eqref{ValphaE} comes from  the  well-known Laplace principle, a classical result in large deviation theory, which states that for any probability measure $\rho\in\mc{P}(\RR^d)$ compactly supported, it holds
\begin{equation}\label{Laplace}
\lim\limits_{\alpha\to\infty}\left(-\frac{1}{\alpha}\log\left(\int_{\RR^d}e^{-\alpha\EE(x)}d\rho(x)\right)\right)=\inf\limits_{x\, \in\, \rm{supp }(\rho)} \EE(x)\,.
\end{equation}
Therefore, for large values of $\alpha \gg 1$ the regularized global best $\Q^{\alpha}_t \approx \Q_t$, where
\[
\Q_t = \argmin\left\{\EE(X_t^1),\EE(X_t^2),\ldots,\EE(X_t^N)\right\}.
\] 
We emphasize that the stochastic particle system \eqref{eq:psociri} has locally Lipschitz coefficients, thus it admits strong solutions and pathwise uniqueness holds up to any finite time $T>0$, see \cite{CS, Dur}. The above system of SDEs in the sequel is considered in a general setting, without necessarily satisfying the PSO constraint \eqref{eq:param}.  

Thanks to the smoothness of the right-hand side in \eqref{eq:psociri}, we can formally derive
the mean-field description of the microscopic system (see \cite{Golse, Jab2, Snitzman}). 
%
Introducing the $N$-particle probability density
\[
f^{(N)}(x_1,\ldots,x_N,v_1,\ldots,v_N,t),
\]
we consider the dynamics of the first marginal 
\[
f_1^{(N)}(x_1,v_1,t) = \int f^{(N)}(x_1,\ldots,x_N,v_1,\ldots,v_N,t)\,d\Omega_1,
\]
where $d\Omega_1=dx_2\ldots\,dx_N\,dv_2\ldots\,dv_N$ is the volume element, and make the so-called propagation of chaos assumption on the marginals. More precisely, one assumes that for  $N\gg 1$ sufficiently large the $N$-particle probability density $f^{(N)} \approx f^{\otimes N}$, i.e. the random pairs $(X_1,V_1)$, $\ldots$, $(X_N,V_N)$ are approximatively independent and each with the same distribution $f(x,v,t)$. As a consequence
\begin{equation}\label{ValphaEc}
\Q^{\alpha}_t \approx \Q^{\alpha} (\rho)=\frac{\int_{\RR^d} x \omega_\alpha(x) \rho(x,t)\,dx}{\int_{\RR^d} \omega_\alpha(x) \rho(x,t)\,dx},\qquad \rho(x,t)=\int_{\RR^d} f(x,v,t)\,dv,
\end{equation}
and the evolution of the distribution $f(x,v,t)$ obeys the nonlinear Vlasov-Fokker-Planck equation
\begin{equation}\label{PDEi}
\begin{split}
\partial_t f + v \cdot \nabla_x f = 
\nabla_v\cdot\left(\frac{\gamma}{\iw } v f + \frac{\lambda}{\iw } (x - \Q^{\alpha}(\rho) )f+\frac{\sigma^2}{2\iw ^2}D(x - \Q^{\alpha}(\rho) )^2\nabla_v f\right)
\end{split}
\end{equation}
where we used the identity
\[
\sum_{k=1}^d \frac{\partial^2}{\partial v^2_k}\left((x-\Q^{\alpha}(\rho))^2_k f\right) = \nabla_{v}\cdot \left(D(x - \Q^{\alpha}(\rho))^2\nabla_{v} f\right)
\]
with $D(x - \Q^{\alpha}(\rho) )^2$ the diagonal matrix given by the square of $D(x - \Q^{\alpha}(\rho) )$. 
Equation \eqref{PDEi} represents the mean-field PSO (MF-PSO) model without local best and 
should be accompanied by initial (and boundary)
data, and normalization
\[
\int_{\RR^d\times \RR^d} f(x,v,t)\,dx\,dv = 1.
\]
We refer to \cite{CFT, BCC, Jab, Golse, Snitzman} and the references therein, for more details and rigorous results about mean-field models of Vlasov-Fokker-Planck type. Note, however, that the presence of $\Q^{\alpha}(\rho)$ makes the Vlasov-Fokker-Planck equation nonlinear and nonlocal. This is nonstandard in the literature and raises several analytical and numerical questions (see \cite{carrillo2018analytical, fhps20-2}).

\subsection{Regularized PSO dynamic with memory and local best}
We consider the second order system of SDEs corresponding to the canonical PSO method where the global best and local best have been regularized as follows 
\begin{eqnarray}
\nonumber
dX^i_t &=& V^i_t dt,\\
d\P_t^i &=& \nu \left(X^i_t-\P^i_t\right)S^\beta(X^i_t,\P^i_t)dt,
\label{eq:psocir}
\\
\nonumber
\iw dV^i_t &=& -\gamma V^i_tdt +\lambda_1\left(\P_t^i-X^i_t\right)dt +\lambda_2\left(\G^\alpha_t-X^i_t\right)dt\\
\nonumber
&& +\sigma_1 D(\P_t^i-X^i_t)dB^{1,i}_t+\sigma_2 D(\G^\alpha_t-X^i_t)dB^{2,i}_t,
\end{eqnarray} 
where, similarly to the previous case, we  introduced the following regularization of the global best position 
\begin{equation}\label{ValphaE2}
\G^{\alpha}_t =\frac{\sum_{i=1}^{N} Y_t^i\omega_\alpha(Y_t^i)}{\sum_{i=1}^{N}\omega_\alpha(Y_t^i)}\,, \qquad  \omega_\alpha(Y_t):=e^{-\alpha\EE(Y_t)}\,.
\end{equation}

Furthermore, in the right hand side of \eqref{eq:psocir} we have replaced the $\sign(x)$ function with a sigmoid, for example the hyperbolic tangent $\tanh(\beta x)$ for $\beta\gg 1$, and consider $S^\beta(y,x)=1+\tanh\left(\beta(\EE(y)-\EE(x))\right)$. Thanks to these regularizations, also the stochastic particle system \eqref{eq:psocir} has locally Lipschitz coefficients and therefore it admits strong solutions and pathwise uniqueness holds for any finite time $T>0$. Even in this case, the system of SDEs \eqref{eq:psocir} is generalized without restricting the search parameters to the PSO constraint \eqref{eq:param}.

In order to derive a mean field description of system \eqref{eq:psocir}, we must introduce an additional dependence from the memory variables in the $N$-particle probability density
\[
f^{(N)}(x_1,\ldots,x_N,y_1,\ldots,y_N,v_1,\ldots,v_N,t),
\]
and consider the dynamics of the first marginal 
\[
f_1^{(N)}(x_1,y_1,v_1,t) = \int_{\RR^d} f^{(N)}(x_1,\ldots,x_N,y_1,\ldots,y_N,v_1,\ldots,v_N,t)\,d\Omega_1,
\]
where now $d\Omega_1=dx_2\ldots\,dx_N\,dy_2\ldots\,dy_N\,dv_2\ldots\,dv_N$ is the volume element. 
Again assuming propagation of chaos, namely that for sufficiently large $N\gg 1$ the $N$- particle probability density factorizes $f^{(N)} \approx f^{\otimes N}$, i.e the random triples $(X_t^i,Y_t^i,V_t^i)$ are independent and with the same distribution $f(x,y,v,t)$, we have
\begin{equation}
\G_t^{\alpha}\approx \G^{\alpha}(\bar\rho) =\frac{\int_{\RR^d} y\,\omega_\alpha(y)\bar\rho(y,t)\,dy}{\int_{\RR^d}\omega_\alpha(y)\bar\rho(y,t)\,dy},\qquad \bar\rho(y,t) = \int_{\RR^d\times \RR^d} f(x,y,v,t)\,dx\,dv.
\label{eq:gbc}
\end{equation}
Additionally, the distribution $f(x,v,t)$ satisfies the nonlinear Vlasov-Fokker-Planck equation
\begin{eqnarray}
\nonumber
&&\partial_t f + v \cdot \nabla_x f + \nabla_y \cdot \left(\nu(x-y)S^\beta(x,y)f\right)= 
\\
\label{PDE}
&&\qquad \nabla_v\cdot\left(\frac{\gamma}{\iw} v f + \frac{\lambda_1}{\iw} (x - y)f
  + \frac{\lambda_2}{\iw} (x - \G^\alpha(\bar\rho) )f\right.\\
  \nonumber
&&\qquad \qquad \left.+\left(\frac{\sigma_2^2}{2{\iw^2}}D(x - \G^\alpha(\bar\rho))^2+\frac{\sigma_1^2}{2{\iw}^2}D(x - y)^2\right)\nabla_v f\right).
\end{eqnarray}
For consistency, initially we assume $f(x,y,v,0)=f_0(x,y,v)$ with $f_0(x,y,v)$ compactly supported and $f_0(x,y,v)\neq 0$ only for $x=y$.
As already mentioned, the rigorous proof of the mean-field limit is an open problem for these interacting particle systems due to the nonlinear terms and the difficulty of managing the multiplicative noise in \eqref{eq:psociri} and \eqref{eq:psocir}. 

\section{Small inertia limit of particle swarm optimization}
In this section we consider the asymptotic behavior of the previous Vlasov-Fokker-Planck equations modelling the PSO dynamic in the small inertia limit. We will derive the corresponding macroscopic equations which permit to recover the recently introduced consensus based optimization (CBO) methods \cite{carrillo2019consensus}. We refer to \cite{DLP} for a theoretical background concerning the related problem of the overdamped limit of nonlinear Vlasov-Fokker-Planck systems.


\subsection{The case without memory effects}
Let us first consider the simplified setting in absence of local best.
To illustrate the limiting procedure, let us observe that for small values of $m \ll 1$ from the second equation in \eqref{eq:psociri} we formally get 
\[
V^i_tdt =\lambda\left(\Q^{\alpha}_t-X^i_t\right)dt+\sigma D(\Q^{\alpha}_t-X^i_t)dB^i_t,
\] 
where we used the fact that $\gamma = 1 - m \approx 1$. Substituting the above identity into the first equation in \eqref{eq:psociri} gives the first order CBO dynamic \cite{carrillo2019consensus}
\begin{equation}
dX^i_t = \lambda\left(\Q^{\alpha}_t-X^i_t\right)dt+\sigma D(\Q^{\alpha}_t-X^i_t)dB^i_t.
\label{eq:CBO}
\end{equation} 
Therefore, the CBO models based on a multiplicative noise can be understood as reduced order approximations of canonical SD-PSO dynamics. Note, however, that in \eqref{eq:CBO} the values of $\lambda$ and $\sigma$ are independent and does not necessarily satisfy the PSO constraints \eqref{eq:param}. 

In the sequel we will develop these arguments in the case of the nonlinear Vlasov-Fokker-Planck equation \eqref{PDEi} describing the mean-field limit dynamic associated to \eqref{eq:psociri}. For notation clarity we denote the small inertia value $m=\e > 0$ in \eqref{PDEi},  
and re-write the scaled Vlasov-Fokker-Planck system in the form
\begin{equation}\label{PDEis}
\begin{split}
\partial_t f + v \cdot \nabla_x f + 
\frac1{\e}\nabla_v\cdot\left(-\e  v f + \lambda (\Q^{\alpha}(\rho)-x)f\right)= L_{\e}(f)
\end{split}
\end{equation}
where we used the fact that $\gamma=1-\e$ and define
\[
\begin{split}
L_{\e}(f)&=\frac1{\e}\nabla_v\cdot\left(v f + \frac{\sigma^2}{2\e} D(x - \Q^{\alpha}(\rho) )^2\nabla_v f\right)\\
&= \frac1{\e} \sum_{j=1}^d \frac{\sigma^2}{2} (x_j - \Q^{\alpha}_j(\rho))^2 \frac{\partial}{\partial v_j}\left(\frac{2 f v_j}{\sigma^2(x_j - \Q^{\alpha}_j(\rho))^2} + \frac1{\e}\frac{\partial f}{\partial v_j}\right).
\end{split}
\]
Let us now introduce the local Maxwellian with unitary mass and zero momentum
\[
\begin{split}
{\mathcal M}_\e(x,v,t)&= \prod_{j=1}^d M_{\e}(x_j,v_j,t), \\
M_{\e}(x_j,v_j,t) &= \frac{\e^{1/2}}{\pi^{1/2}\sigma |x_j-\Q^\alpha_j(\rho)|} 
\exp\left\{-\frac{\e v_j^2}{\sigma^2(x_j-\Q^\alpha_j(\rho))^2}\right\},
\end{split}
\]
then we have
\[
L_{\e}(f) = \frac1{\e^2}\sum_{j=1}^d \frac{\sigma^2}{2} (x_j - \Q^{\alpha}_j(\rho))^2 \frac{\partial}{\partial v_j}\left(f\frac{\partial}{\partial v_j}\log\left(\frac{f}{M_{\e}(x_j,v_j,t)}\right)\right).
\]
Therefore $L_{\e}(f)$ is of order $1/\e^2$ and we can write for small values of $\e \ll 1$
\begin{equation}
f(x,v,t)=\rho(x,t){\mathcal M}_\e(x,v,t).
\label{eq:maxw1}
\end{equation}
Let us now integrate equation \eqref{PDEis} with respect to $v$, we get
\[
\begin{split}
\frac{\partial \rho}{\partial t} + \nabla_x \cdot (\rho u) &= 0\\
\frac{\partial \rho u}{\partial t} + \int_{\RR^d} v\left(v\cdot \nabla_x f\right)\,dv &= \frac{1-\e}{\e} \rho u + \frac1{\e} \lambda (\Q^{\alpha}(\rho)-x) \rho
\end{split}
\]
where
\[
\rho u = \int_{\RR^d} f(x,v,t) v\,dv.
\]
Now assuming \eqref{eq:maxw1} we can compute for $\e \ll 1$ the $i$-th component as
\[
\begin{split}
\int_{\RR^d} v_i\left(v\cdot \nabla_x \left(\rho(x,t){\mathcal M}_\e(x,v,t)\right) \right)\,dv &=  \sum_{j=1}^d \frac{\partial}{\partial x_j} \left(\rho(x,t) \int_{\RR^d} v_i (v_j {\mathcal M}_\e(x,v,t))\,dv\right)\\
&=  \frac{\partial}{\partial x_i} \left(\rho(x,t) \int_{\RR} v^2_i {M}_\e(x_i,v_i,t)\,dv_i\right)\\
&= \frac{\sigma^2}{2\e}\frac{\partial}{\partial x_i} \left(\rho(x,t) (x_i-\Q^\alpha_i(\rho))^2\right)
\end{split}
\]
which provides the second order macroscopic model
\begin{equation}
\begin{split}
\frac{\partial \rho}{\partial t} + \nabla_x \cdot (\rho u) &= 0\\
\frac{\partial (\rho u)_i}{\partial t} + \frac{\sigma^2}{2\e} \frac{\partial}{\partial x_i} \left(\rho(x,t) (x_i-\Q^\alpha_i(\rho))^2\right) &= -\frac{1-\e}{\e} (\rho u)_i + \frac1{\e} \lambda (\Q_i^{\alpha}(\rho)-x_i) \rho.
\label{eq:macro}
\end{split}
\end{equation}
Formally, as $\e\to 0$, from the second equation in \eqref{eq:macro} we get 
\[
(\rho u)_i = \lambda (\Q_i^{\alpha}(\rho)-x_i) \rho -\frac{\sigma^2}{2} \frac{\partial}{\partial x_i} \left(\rho(x,t) (x_i-\Q^\alpha_i(\rho))^2\right),
\]
which substituted in the first equation yields the mean-field CBO system \cite{carrillo2019consensus}
\begin{equation}
\frac{\partial \rho}{\partial t} + \nabla_x \cdot \lambda (\Q^{\alpha}(\rho)-x) \rho = \frac{\sigma^2}{2}\sum_{j=1}^d \frac{\partial^2}{\partial x^2_j} \left(\rho(x,t) (x_j-\Q^\alpha_j(\rho))^2\right).
\label{eq:CBOp}
\end{equation}
Therefore, in the small inertia limit we expect the macroscopic density in the PSO system \eqref{PDEi} to be well approximated by the solution of the CBO equation \eqref{eq:CBOp}. We remark that this is not the case for the original CBO method proposed in \cite{pinnau2017consensus} where the noise is not in component-wise form.

\subsection{The general case with memory}
Next, we consider the same small inertia scaling in the general case with dependence from the local best. Again, we can first illustrate the result by considering the behaviour for $m \ll 1$ of the SD-PSO system \eqref{eq:psocir}. We formally get from the third equation
\[
V^i_tdt =\lambda_1\left(\P_t^i-X^i_t\right)dt +\lambda_2\left(\G^\alpha_t-X^i_t\right)dt+\sigma_1 D(\P_t^i-X^i_t)dB^{1,i}_t+\sigma_2 D(\G^\alpha_t-X^i_t)dB^{2,i}_t,
\] 
which inserted into the first equation in \eqref{eq:psocir} corresponds to a novel first order CBO dynamic with local best
\begin{eqnarray}
\nonumber
dX^i_t &=& \lambda_1\left(\P_t^i-X^i_t\right)dt+\lambda_2\left(\G^\alpha_t-X^i_t\right)dt+\sigma_1 D(\P_t^i-X^i_t)dB^{1,i}_t+\sigma_2 D(\G^\alpha_t-X^i_t)dB^{2,i}_t,\\[-.1cm]
\label{eq:psocirlb}
\\[-.2cm]
\nonumber
d\P_t^i &=& \nu \left(X^i_t-\P^i_t\right)S^\beta(X^i_t,\P^i_t)dt.
\end{eqnarray} 
In contrast with the model recently introduced in \cite{TW} the above first order CBO method avoids backward time integration through the use of an additional differential equation. We remark that at the SDEs level, by analogous arguments as the one presented in this paper, in principle even the CBO model \cite{TW} can be derived as the small inertia limit of the corresponding PSO model where memory effects are modeled as in  \cite{TW}.

Concerning the corresponding MF-PSO limit we can essentially perform analogous computations as in the previous section. Thus, after setting $m=\e>0$ we consider the scaled system
\begin{equation}
\begin{split}
\partial_t f + v \cdot \nabla_x f &+ \nabla_y \cdot \left(\nu(x-y)S^\beta(x,y)f\right)\\
&+\frac1{\e}\nabla_v\cdot\left(-\e v f + {\lambda_1}(y-x)f
  + {\lambda_2}(\G^\alpha(\bar\rho) - x)f\right) = L_\e(f), 
  \end{split}
  \label{eq:PDEs}
\end{equation}
where now 
\[
\begin{split}
L_\e(f)&= \frac1{\e}\nabla_v\cdot\left(v f + \frac{\sigma_2^2}{2{\e}}D(x - \G^\alpha(\bar\rho))^2\nabla_v f+\frac{\sigma_1^2}{2{\e}}D(x - y)^2\nabla_v f\right)\\
&=\frac1{2\e}\sum_{j=1}^d \Sigma(x_j,y_j,t)^2\frac{\partial}{\partial v_j}\left(
\frac{2fv_j}{\Sigma(x_j,y_j,t)^2}+\frac1{\e}\frac{\partial f}{\partial v_j}\right)
\end{split}
\]
and we use the notation
\[
\Sigma(x_j,y_j,t)^2 = {\sigma_2^2}(x_j - \G^\alpha_j(\bar\rho))^2+{\sigma_1^2}(x_j - y_j)^2.
\]
Then, introducing the local Maxwellian
\[
\begin{split}
{\mathcal M}_\e(x,y,v,t)&= \prod_{j=1}^d M_{\e}(x_j,y_j,v_j,t), \\
M_{\e}(x_j,y_j,v_j,t) &= \frac{\e^{1/2}}{\pi^{1/2} |\Sigma(x_j,y_j,t)|} 
\exp\left\{-\frac{\e v_j^2}{\Sigma(x_j,y_j,t)^2}\right\},
\end{split}
\]
with unitary mass and zero momentum we have
\[
L_{\e}(f) = \frac1{2\e^2}\sum_{j=1}^d \Sigma(x_j,y_j,t)^2 \frac{\partial}{\partial v_j}\left(f\frac{\partial}{\partial v_j}\log\left(\frac{f}{M_{\e}(x_j,y_j,v_j,t)}\right)\right).
\]
We can thus write for $\e \ll 1$
\begin{equation}
f(x,y,v,t)=\rho(x,y,t){\mathcal M}_\e(x,y,v,t),
\label{eq:maxw2}
\end{equation}
and after integrating \eqref{eq:PDEs} with respect to $v$ and using the approach \eqref{eq:maxw2}, we get the second order macroscopic model
\begin{equation}
\begin{split}
\frac{\partial \rho}{\partial t} + \nabla_x \cdot (\rho u) + \nabla_y \cdot \left(\nu(x-y)S^\beta(x,y)\rho\right)&= 0\\
\frac{\partial (\rho u)_i}{\partial t} + \frac{\sigma^2}{2\e}\frac{\partial}{\partial x_i} \left(\rho(x,t) \Sigma(x_i,y_i,t)^2\right) &= -\frac{1-\e}{\e} (\rho u)_i + \frac1{\e} \left(\lambda_1 (y_i-x_i)+\lambda_2 (\G_i^\alpha(\bar\rho) - x_i) \right)\rho.
\label{eq:macro2}
\end{split}
\end{equation}
Formally, as $\e\to 0$, the above system reduces to a novel mean-field CBO system with local best
\begin{equation}
\begin{split}
\frac{\partial \rho}{\partial t} + \nabla_x \cdot \left(\lambda_1 (y-x)+\lambda_2 (\G^\alpha(\bar\rho) - x) \right)\rho + \nabla_y \cdot \left(\nu(x-y)S^\beta(x,y)\rho\right)\\
 = \frac{1}{2}\sum_{j=1}^d \frac{\partial^2}{\partial x^2_j} \left(\rho(x,t)\left( \sigma_1^2 (x_j-y_j)^2+\sigma_2^2 (x_j-\G_j^\alpha(\bar\rho))^2\right)\right).
 \end{split}
\label{eq:CBOlb}
\end{equation}

\section{Numerical examples}
In this section we present several numerical tests in order to verify the validity of the previous theoretical analysis, namely the mean field limit and the small inertial limit, and to analyze the performance of the methods based on SD-PSO against various prototype global optimization functions.

\begin{figure}[tb]
\begin{minipage}{\linewidth}
\centering
\subcaptionbox{Ackley}{\includegraphics[scale=0.5]{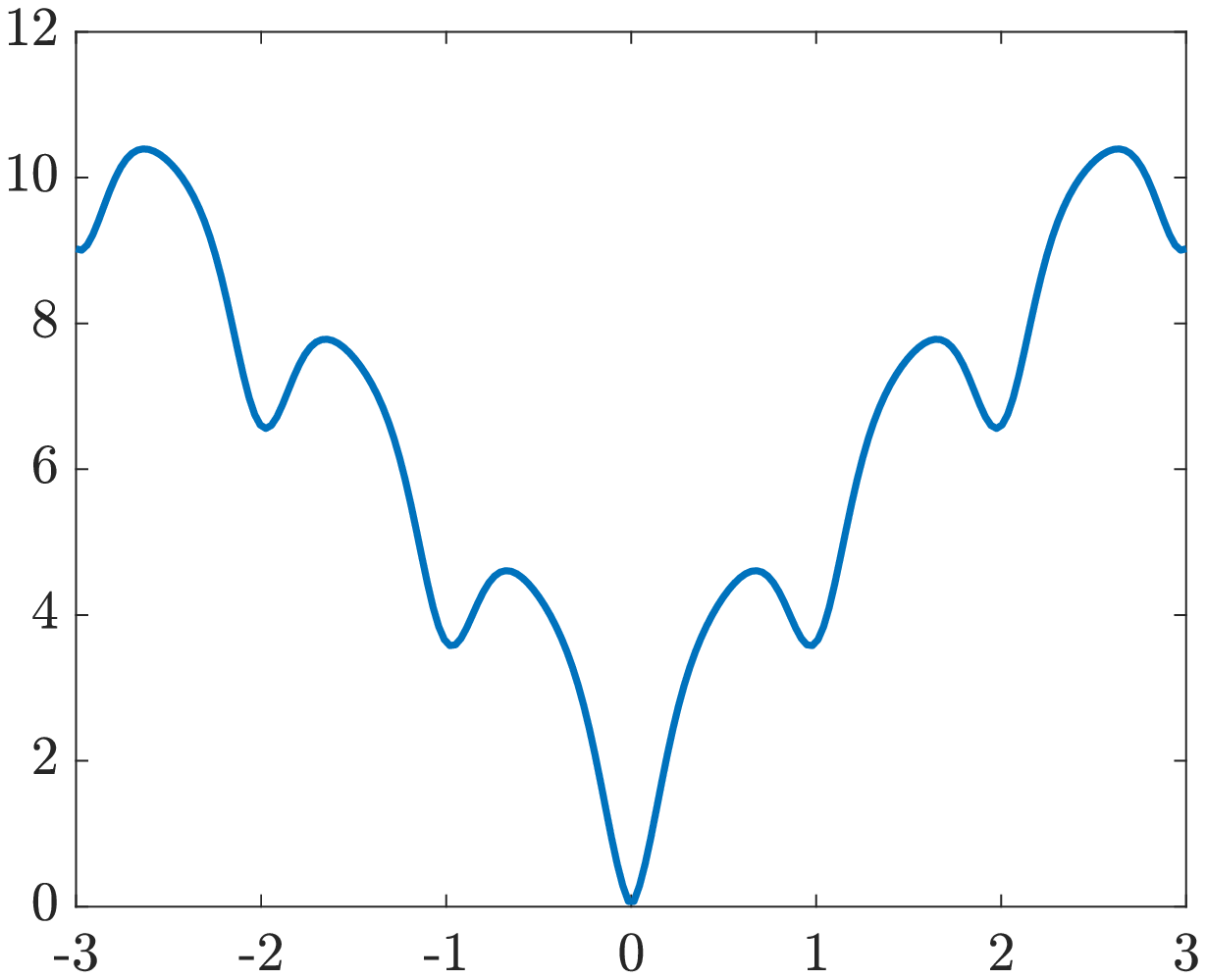}}\qquad
\subcaptionbox{Rastrigin}{\includegraphics[scale= 0.5]{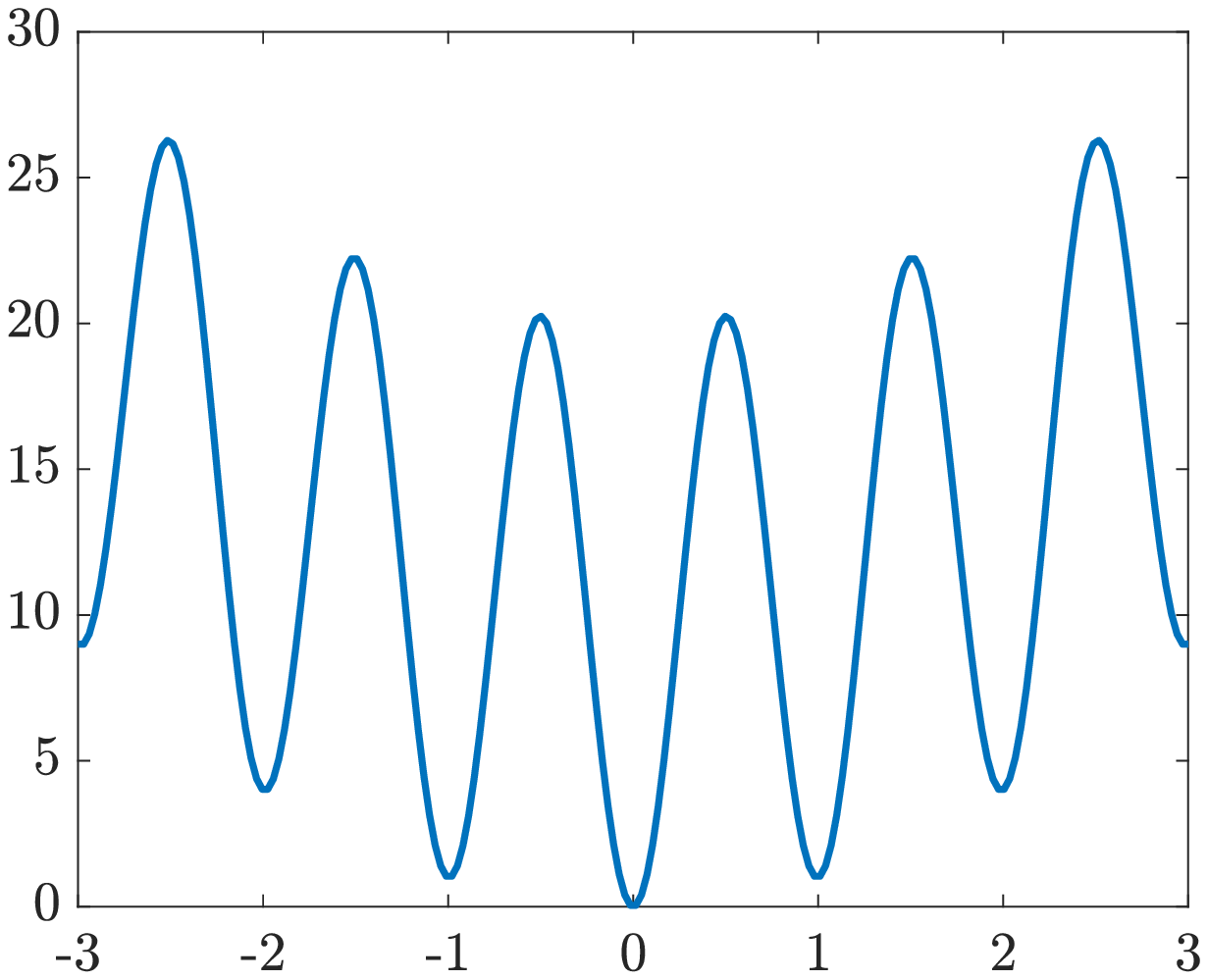}} 
\caption{One-dimensional Ackley and Rastrigin functions in the interval $\left[ -3, 3 \right]$ with global minimum in the origin.}
\label{Fig1}
\end{minipage}
\end{figure}

\subsection{Validation of the mean field limit}
In the following we present three numerical test cases to validate the mean field limit process in one dimension by considering as prototype functions for global optimization the Ackley function \eqref{eq:Ackley} and the Rastrigin function \eqref{eq:Rastrigin}. Both functions have multiple local minima that can easily trap the particle dynamics (see Figure \ref{Fig1}) and have been used recently to test consensus based particle optimizers \cite{pinnau2017consensus,carrillo2019consensus,fhps20-1}.

First we have considered the case without memory effect formulated by the SD-PSO system \eqref{eq:psociri} that uses only the action of the global best compared to the solution obtained using the mean field limit \eqref{PDEi}. The particle system \eqref{eq:psociri} is solved by 
\begin{eqnarray}
\nonumber
X^{n+1}_i &=& X^{n}_i + \Delta t \ V^{n+1}_i ,
\\[-.25cm]
\label{eq:psoiDiscr}
\\[-.25cm]
\nonumber
\iw V^{n+1}_i &=& \iw V^{n}_i - \gamma \Delta t \ V_i^{n+1} + \lambda \Delta t  \left(\Q_\alpha^n-X^n_i\right)  + \sigma \sqrt{\Delta t} \ D(\Q_\alpha^n-X^n_i) \ \theta^n_{i},
\end{eqnarray}
where $\theta_{i} \sim \mathcal{N}(0,1)$ and the last equation can be rewritten as
\vspace{5pt}
\begin{eqnarray}
\nonumber
V^{n+1}_i &=& \left( \frac{m}{m + \gamma \ \Delta t}\right) V^{n}_i + \frac{\lambda \ \Delta t }{m + \gamma \ \Delta t}\left(\Q_{\alpha}^n-X^n_i\right)+ \frac{\sigma \ \sqrt{\Delta t} }{m + \gamma \ \Delta t} D(\Q_{\alpha}^n-X^n_i) \ \theta^n_{i}.
\end{eqnarray}
The corresponding MF-PSO limit \eqref{PDEi} has been discretized using a dimensional splitting where the transport part is solved through a second order backward semi-Lagrangian method and the remaining Fokker-Planck term is discretized using an implicit central scheme. This permits to avoid restrictive CFL conditions and to obtain second order accuracy in space and velocity. 

In the second and third test cases we introduced the memory variable, initially with the action of the local best only, and then with both local and global dynamics. In this case, the SD-PSO system \eqref{eq:psocir} is solved by 
\begin{eqnarray}
\nonumber
\P^{n+1}_i &=& \P^{n}_i + \nu \ \Delta t \left(X^{n+1}_i-\P^{n}_i\right)S^\beta(X^{n+1}_i,\P^{n}_i), \\
\nonumber
X^{n+1}_i &=& X^{n}_i + \Delta t \ V^{n+1}_i ,\\[-.3cm]
\label{eq:psoDiscr}
\\[-.3cm]
\nonumber
V^{n+1}_i &=& \left( \frac{m}{m + \gamma \ \Delta t}\right) V^{n}_i + \frac{\lambda_1 \ \Delta t }{m + \gamma \ \Delta t}\left(\P^n_i-X^n_i\right)+ \frac{\lambda_2 \ \Delta t }{m + \gamma \ \Delta t}\left(\G_\alpha^n-X^n_i\right)\\
\nonumber
&& \frac{\sigma_1 \ \sqrt{\Delta t} }{m + \gamma \ \Delta t} D(\P^n_i-X^n_i) \ \theta^n_{1,i} + \frac{\sigma_2 \ \sqrt{\Delta t} }{m + \gamma \ \Delta t} D(\G_\alpha^n-X^n_i) \ \theta^n_{2,i},
\end{eqnarray}
where $\theta_{1,i}$, $\theta_{2,i} \sim \mathcal{N}(0,1)$. Note that, the above discretization is equivalent to the discrete PSO system \eqref{eq:psoi} under assumptions \eqref{eq:param} for $ \Delta t = 1 $, $\nu=0.5$, and taking the limit $\alpha$, $\beta \rightarrow \infty$ so that $\P_i^n$, $\G_{\alpha}^n$ match the local and global best definitions in \eqref{eq:psoEvolution}.
The limiting MF-PSO equation \eqref{PDE} is solved by a further dimensional splitting where the additional memory term is discretized using a Lax-Wendroff method that permits to achieve overall second order accuracy. We mention here that we tested also various other approaches for the discretization of the differential memory term. However, in our numerical results we have found essential for the accuracy of the mean-field solution in presence of local best, to discretize the differential term modeling particles' memory using a second order low dissipative scheme. Finally, concerning the time approximation, we implemented both conventional splitting as well as second order Strang splitting without noticing relevant differences in the results. We also tested several set of parameters and initial data (uniform, Gaussian) without observing significant changes with respect to the selection of results reported in the sequel. 

In all test cases we used $N=5 \times 10^5$ particles, a mesh size for the mean field solver of $90 \times 120$ points for $(x,v)\in [-3,3]\times [-4,4]$, and whenever present, the mesh and domain size in $y$ have been taken identical to those in $x$. The choice of the particle number was based on having a good compromise between the convergence to the mean-field limit and the possibility to still visually distinguish the two solutions in the figures. In the deterministic discretization the boundary conditions have been implemented assuming $f(x,v,t)=0$ or $f(x,y,v,t)=0$ outside the computational domain.  

\subsubsection*{Case \#1: MF-PSO without memory effects}
\begin{figure}[H] 
\begin{minipage}{\linewidth}
\centering
\subcaptionbox{Particle solution, $t = 0.5$}{\includegraphics[scale=0.35]{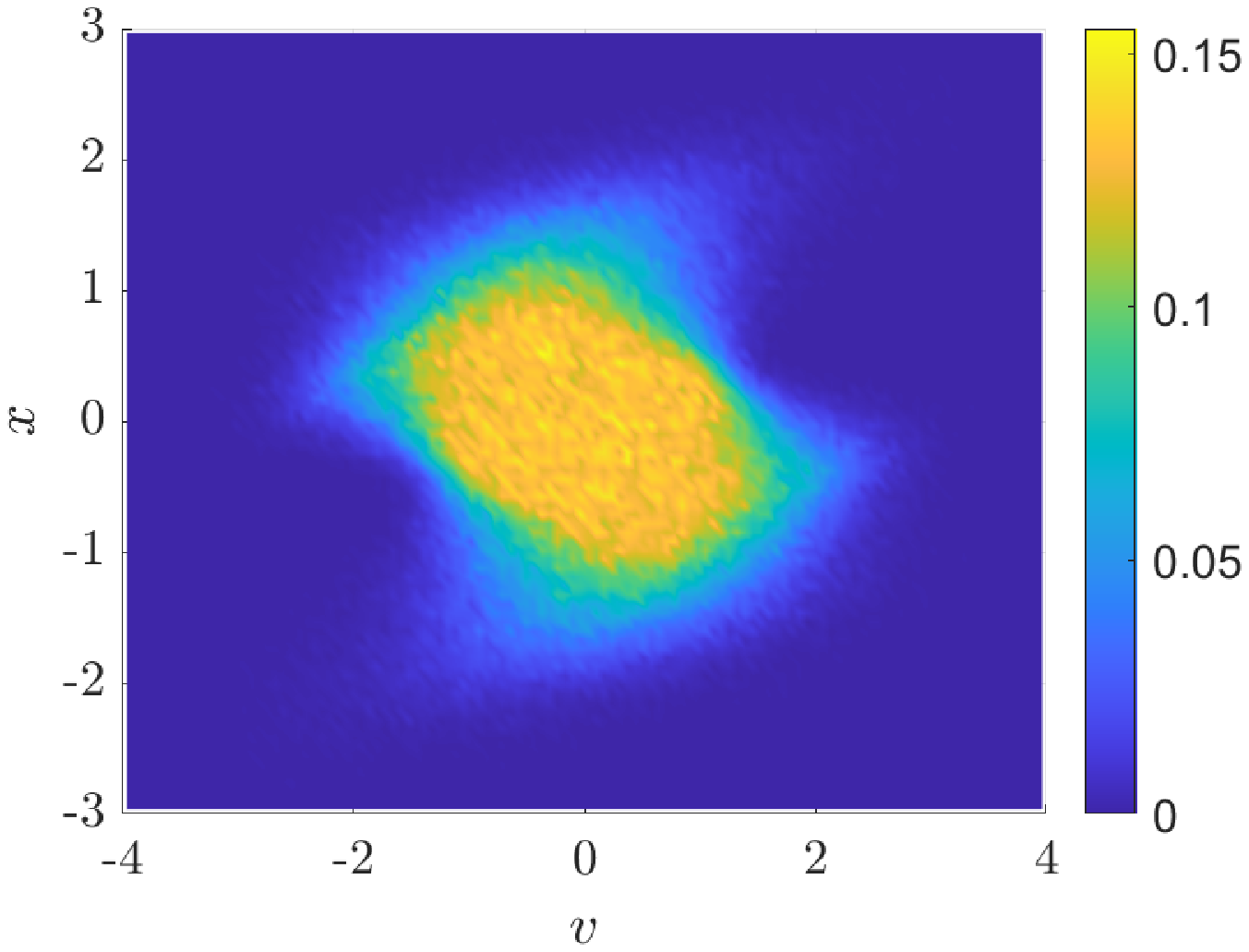}}
\subcaptionbox{Particle solution, $t = 1$}{\includegraphics[scale=0.35]{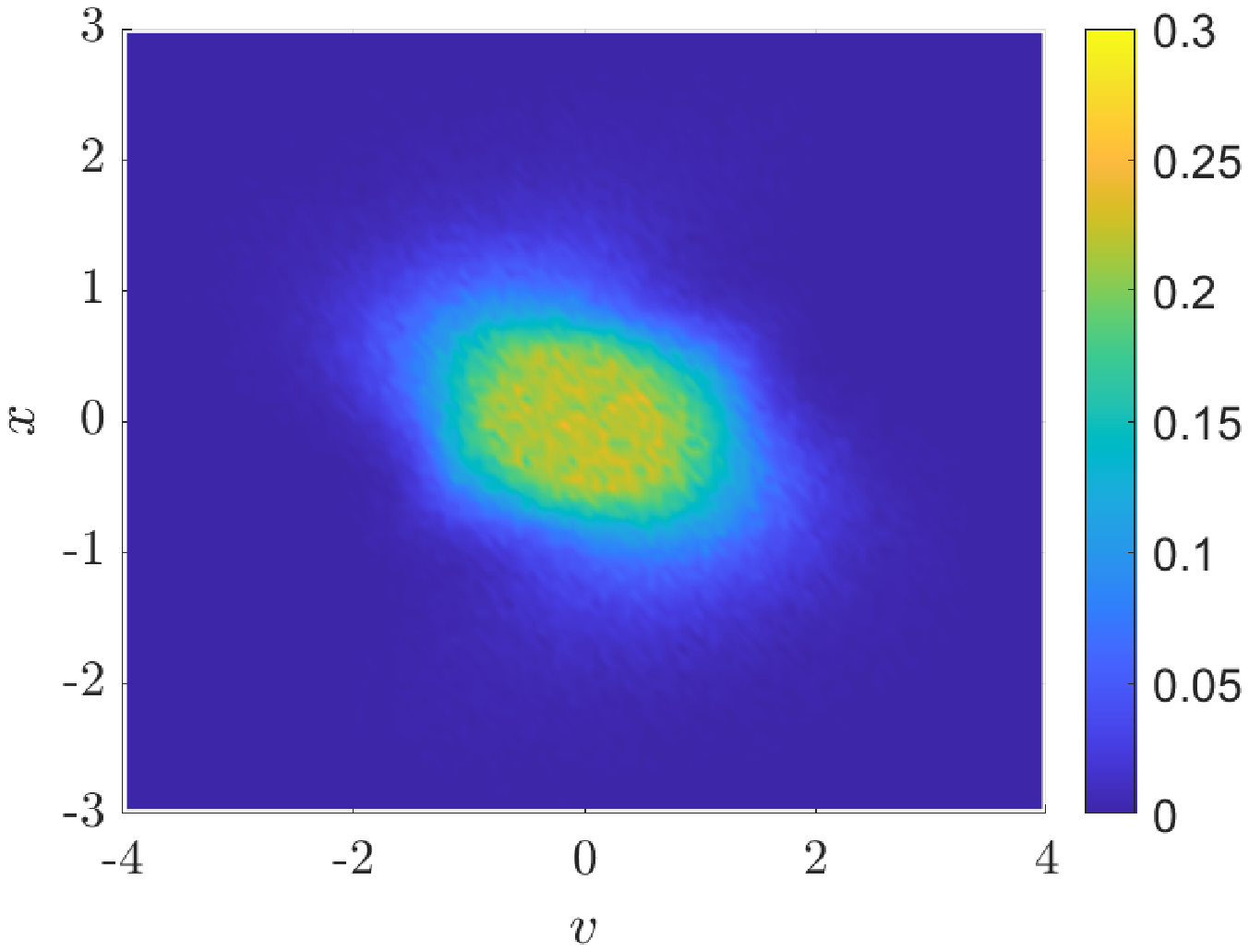}}
\subcaptionbox{Particle solution, $t = 3$}{\includegraphics[scale=0.35]{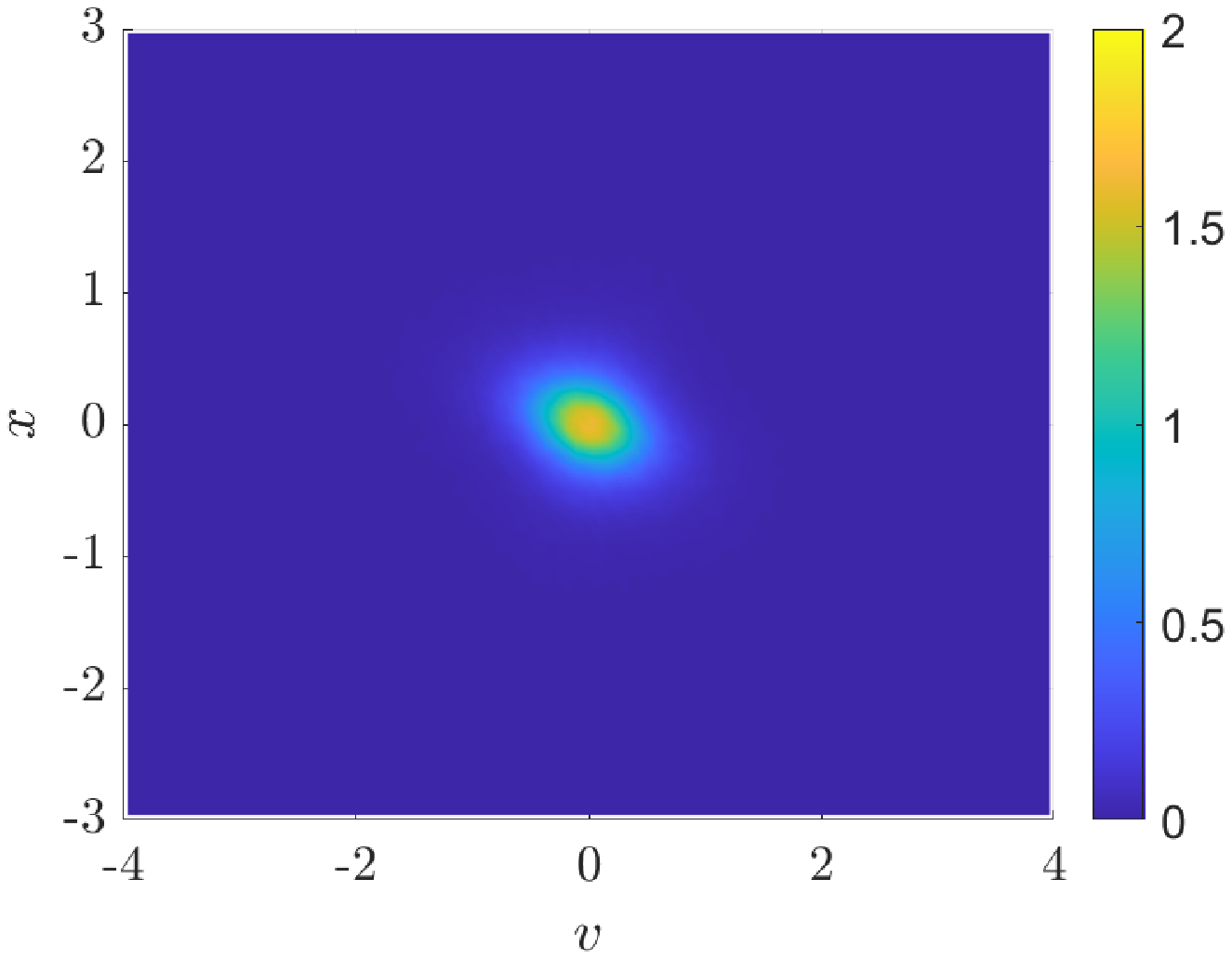}}\\
\subcaptionbox{Mean-field solution, $t = 0.5$}{\includegraphics[scale=0.35]{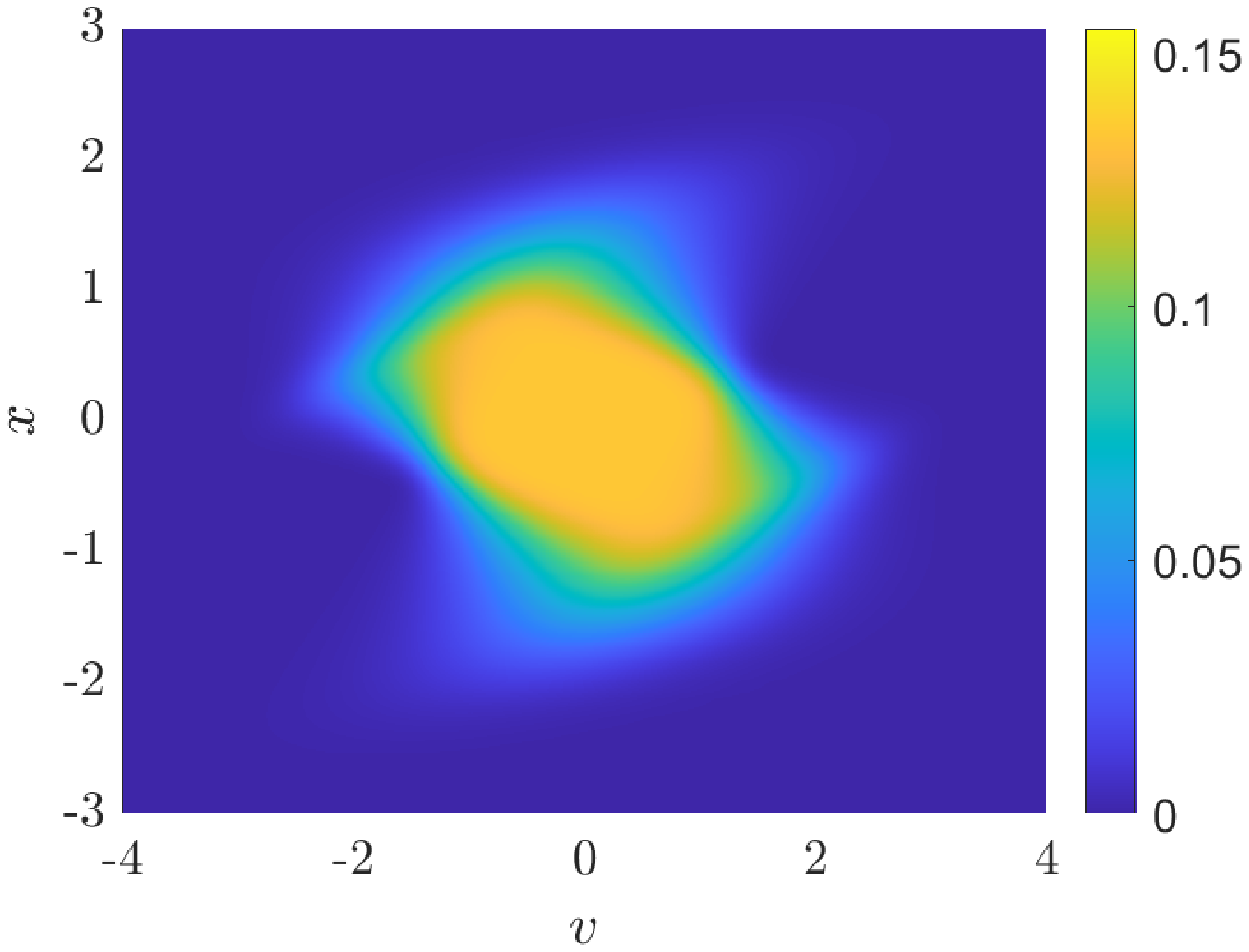}} 
\subcaptionbox{Mean-field solution, $t = 1$}{\includegraphics[scale=0.35]{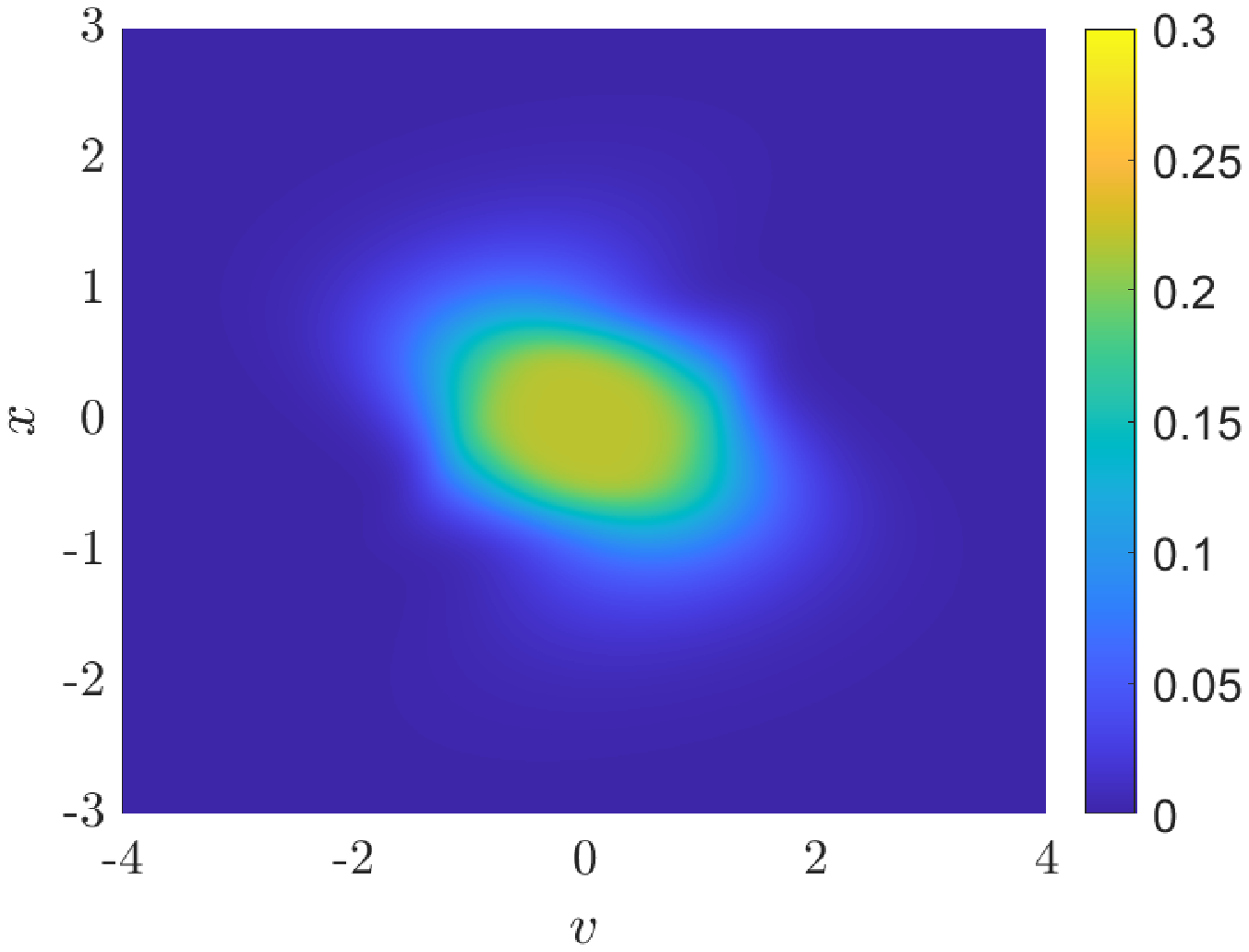}} 
\subcaptionbox{Mean-field solution, $t = 3$}{\includegraphics[scale=0.35]{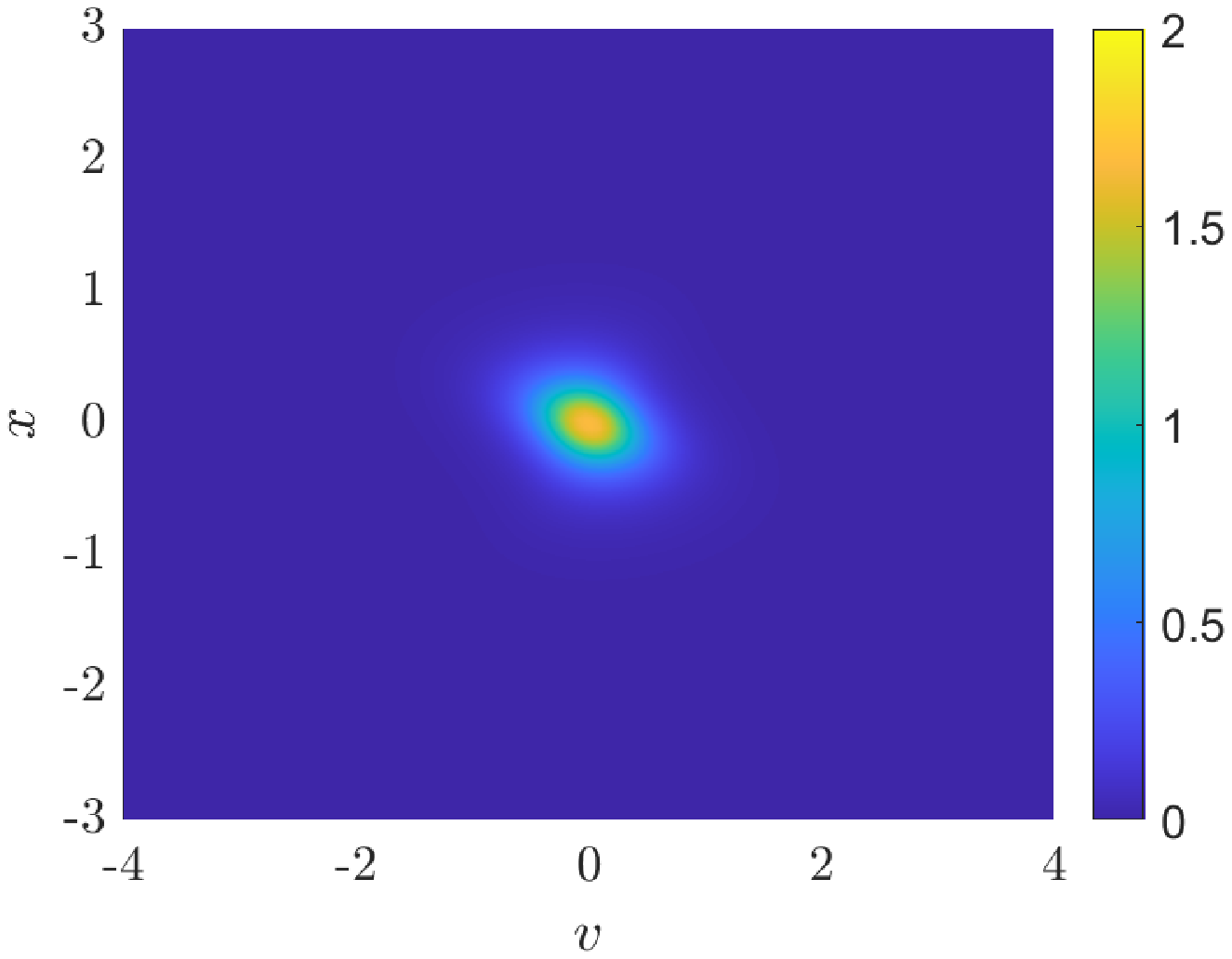}} 
\caption{Case \#1 (no memory). Optimization of the one-dimensional Ackley function with minimum in $x=0$. First row: solution of the SD-PSO system \eqref{eq:psociri}. Second row: solution of the MF-PSO limit \eqref{PDEi}.} 
\label{Fig2}
\end{minipage}
\vspace{5pt}
\\
\begin{minipage}{\linewidth}
\centering
\subcaptionbox{$\rho(x,t)$, $t = 0.5$}{\includegraphics[scale=0.35]{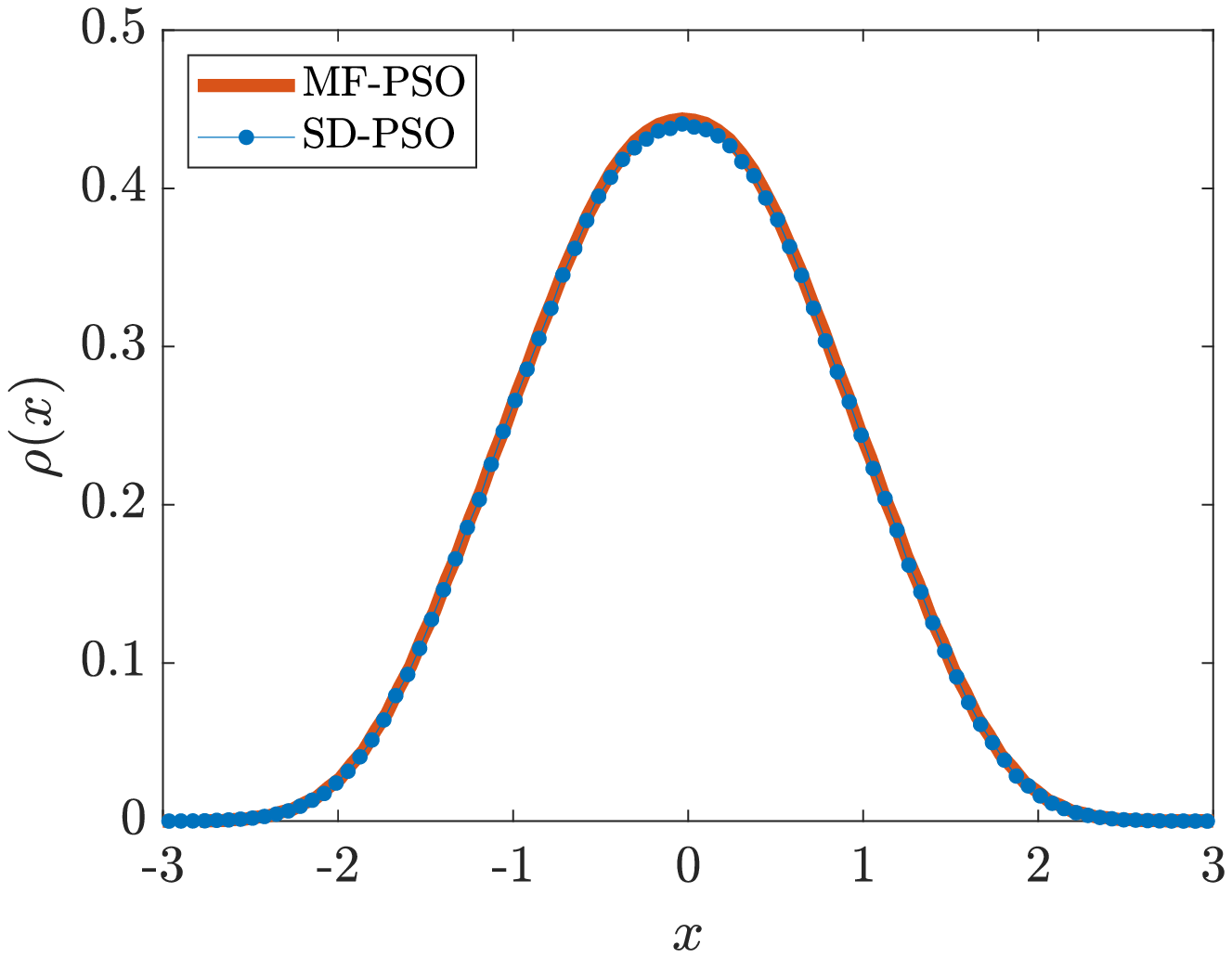}}\
\subcaptionbox{$\rho(x,t)$, $t = 1$}{\includegraphics[scale=0.35]{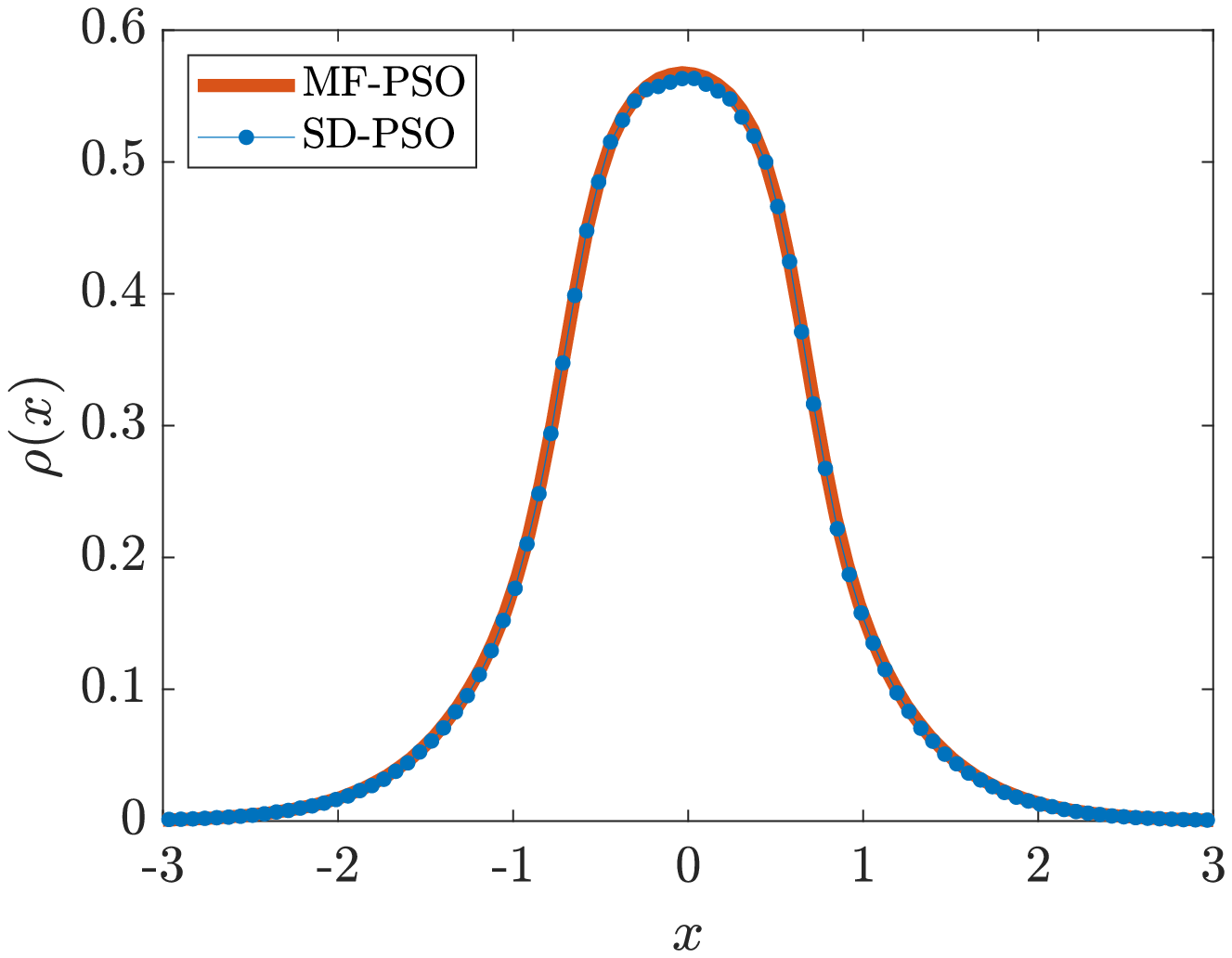}} \
\subcaptionbox{$\rho(x,t)$, $t = 3$}{\includegraphics[scale=0.35]{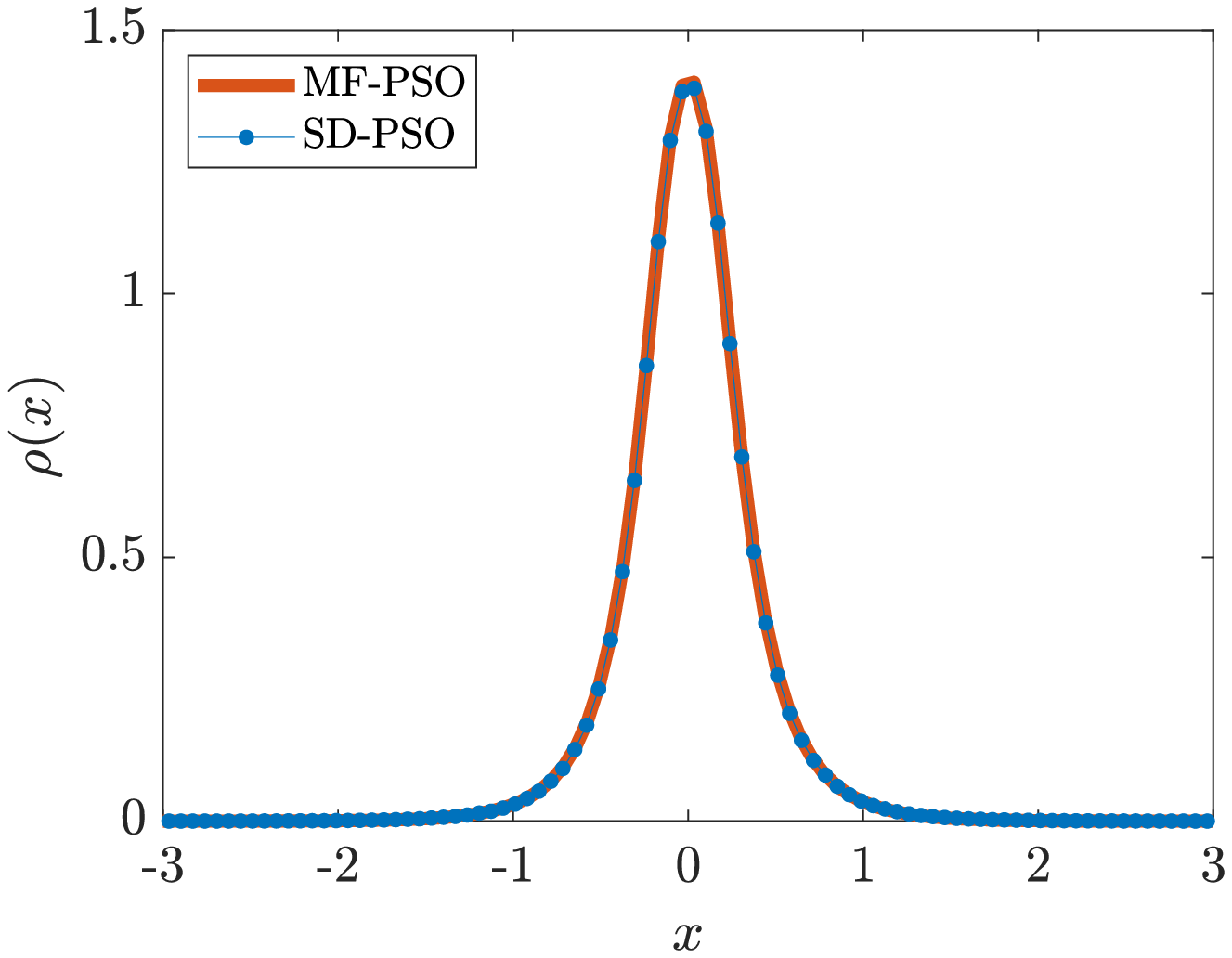}}\\
\caption{Case \#1 (no memory). Evolution of the density $\rho(x,t)$ of the SD-PSO system \eqref{eq:psociri} and the MF-PSO limit \eqref{PDEi} for the one-dimensional Ackley function with minimum in $x=0$.}
\label{Fig4}
\end{minipage}
\end{figure}

We consider the optimization process of the Ackley function with global minimum in the origin $x=0$. Here we report the results obtained with
\begin{equation}
\gamma = 0.5,\quad  \lambda = 1, \quad \sigma= {1}/{\sqrt{3}},\quad \alpha = 30.
\label{eq:paramt} 
\end{equation}
Note that, the values of $\lambda$ and $\sigma$ are compatible with the usual choice $ c_k = 2 $ in \eqref{eq:param}. In Figure \ref{Fig2} we report the contour plots of the evolution, at times $ t = 0.5 $, $ t = 1 $ and $ t = 3 $, of the particle distribution computed through \eqref{eq:psoiDiscr} and by the direct discretization of the mean-field equation \eqref{PDEi}. The initial distribution is taken uniform in all simulations.

To emphasize the good agreement between the results obtained from the resolution of the large particle limit of the SD-PSO model and the results of its corresponding MF-PSO, in Figure \ref{Fig4} we report the evolution in time of the marginal density $\rho(x,t) = \int_{\RR^d} f(x,v,t)\,dv$. The convergence towards a Dirac delta centered in the origin is very similar in both dynamics.

\begin{figure}[H]
\begin{minipage}{\linewidth}
\centering
\subcaptionbox{Particle solution, $t = 0.5$}{\includegraphics[scale=0.35]{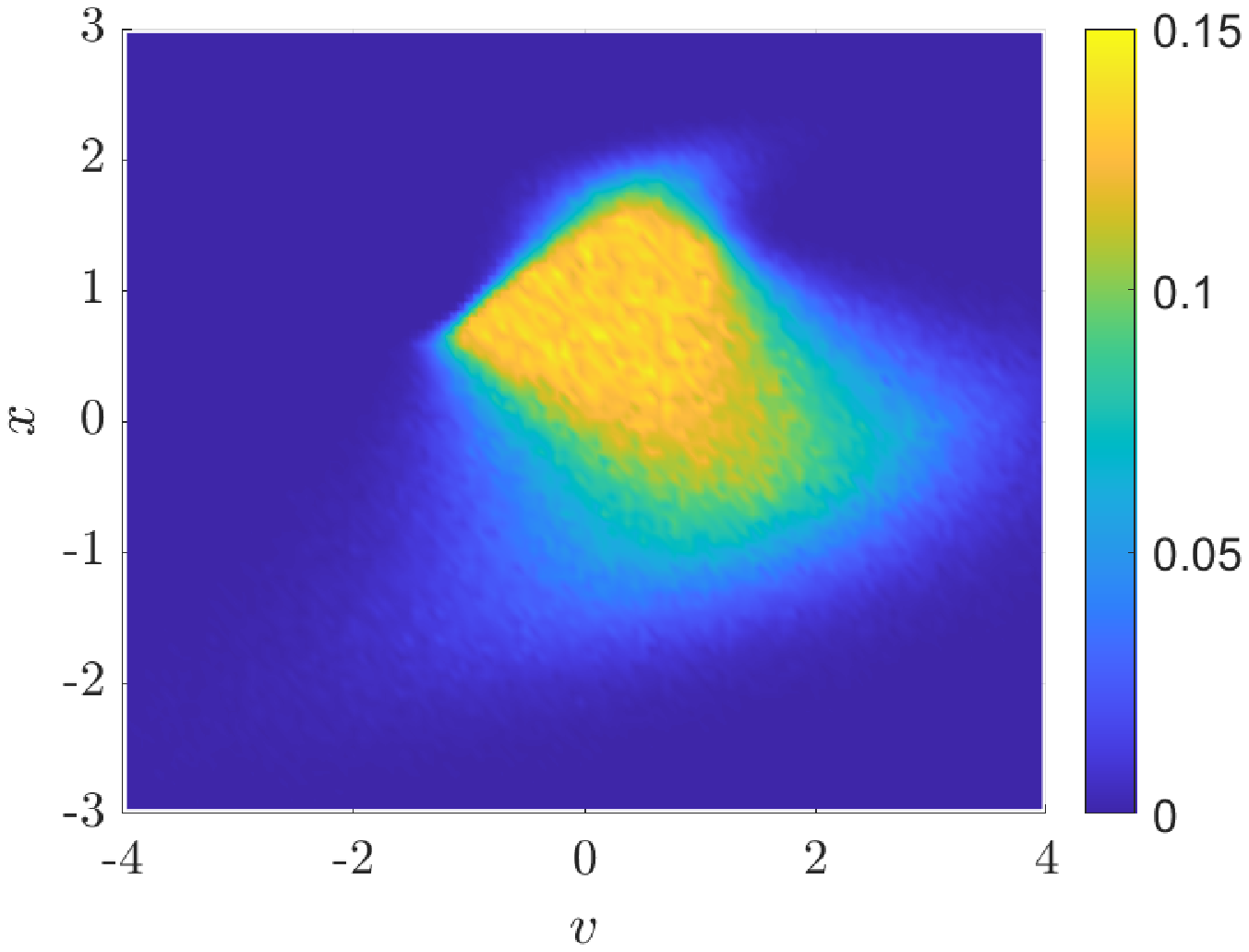}}
\subcaptionbox{Particle solution, $t = 1$}{\includegraphics[scale=0.35]{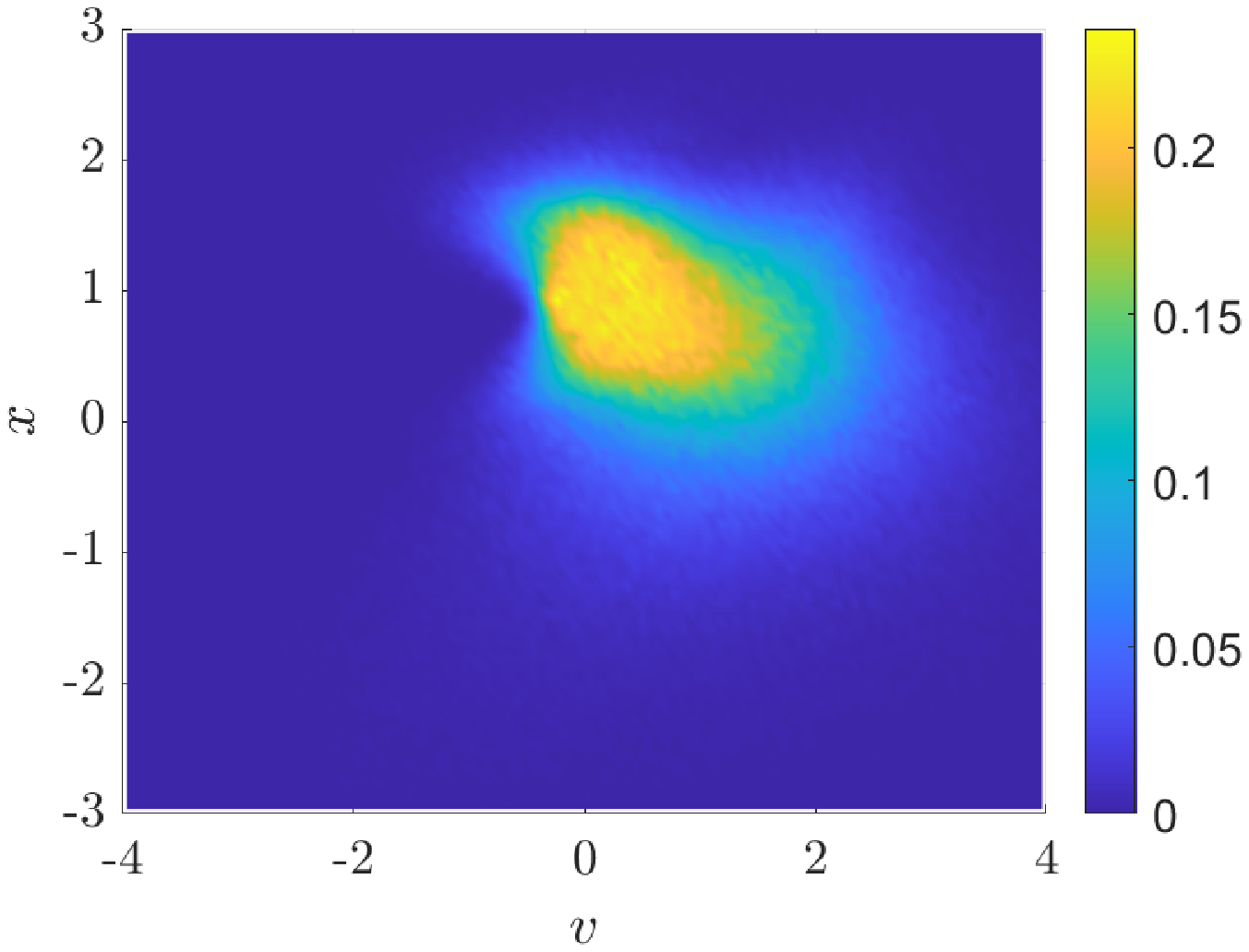}}
\subcaptionbox{Particle solution, $t = 3$}{\includegraphics[scale=0.35]{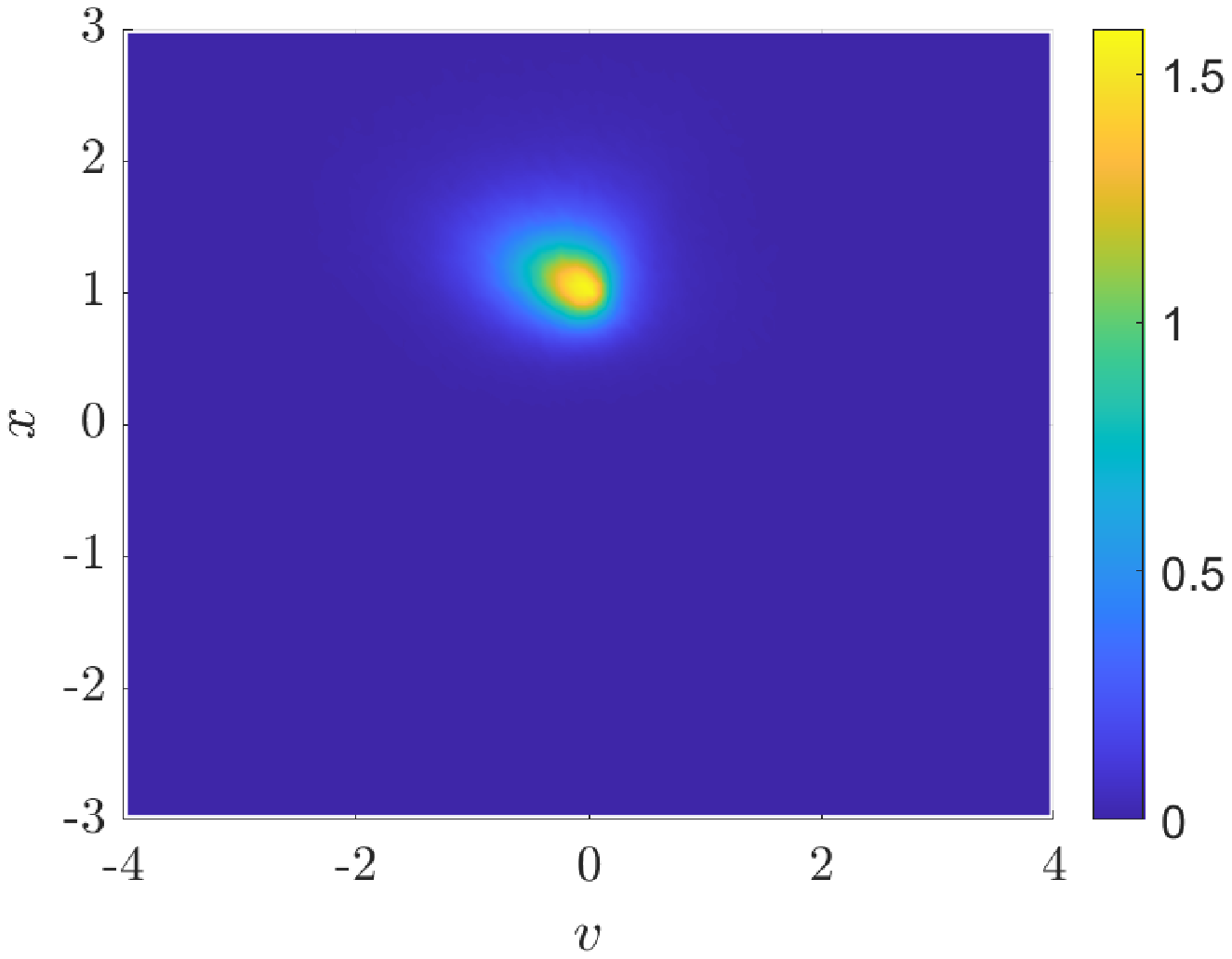}}\\
\subcaptionbox{Mean-field solution, $t = 0.5$}{\includegraphics[scale=0.35]{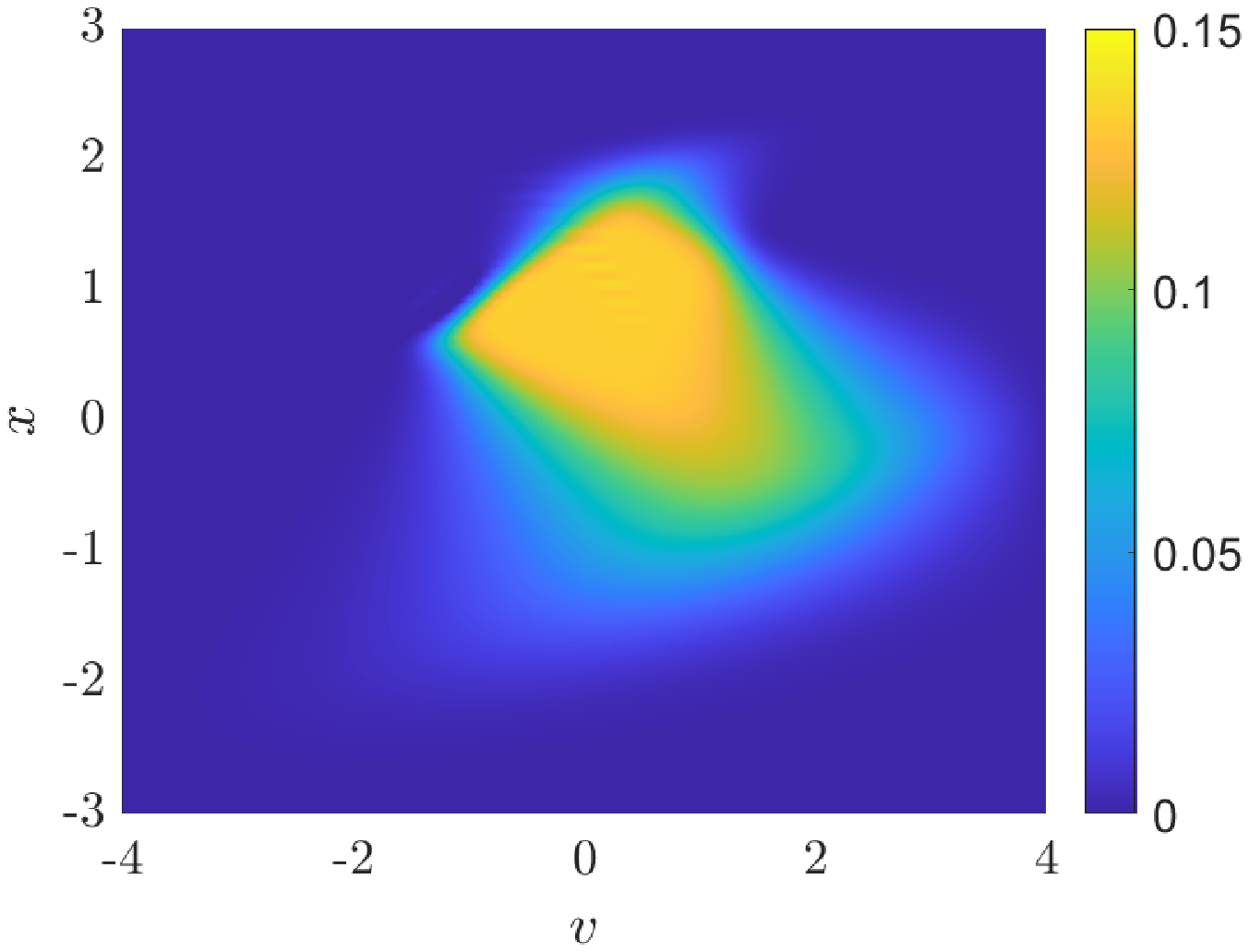}} 
\subcaptionbox{Mean-field solution, $t = 1$}{\includegraphics[scale=0.35]{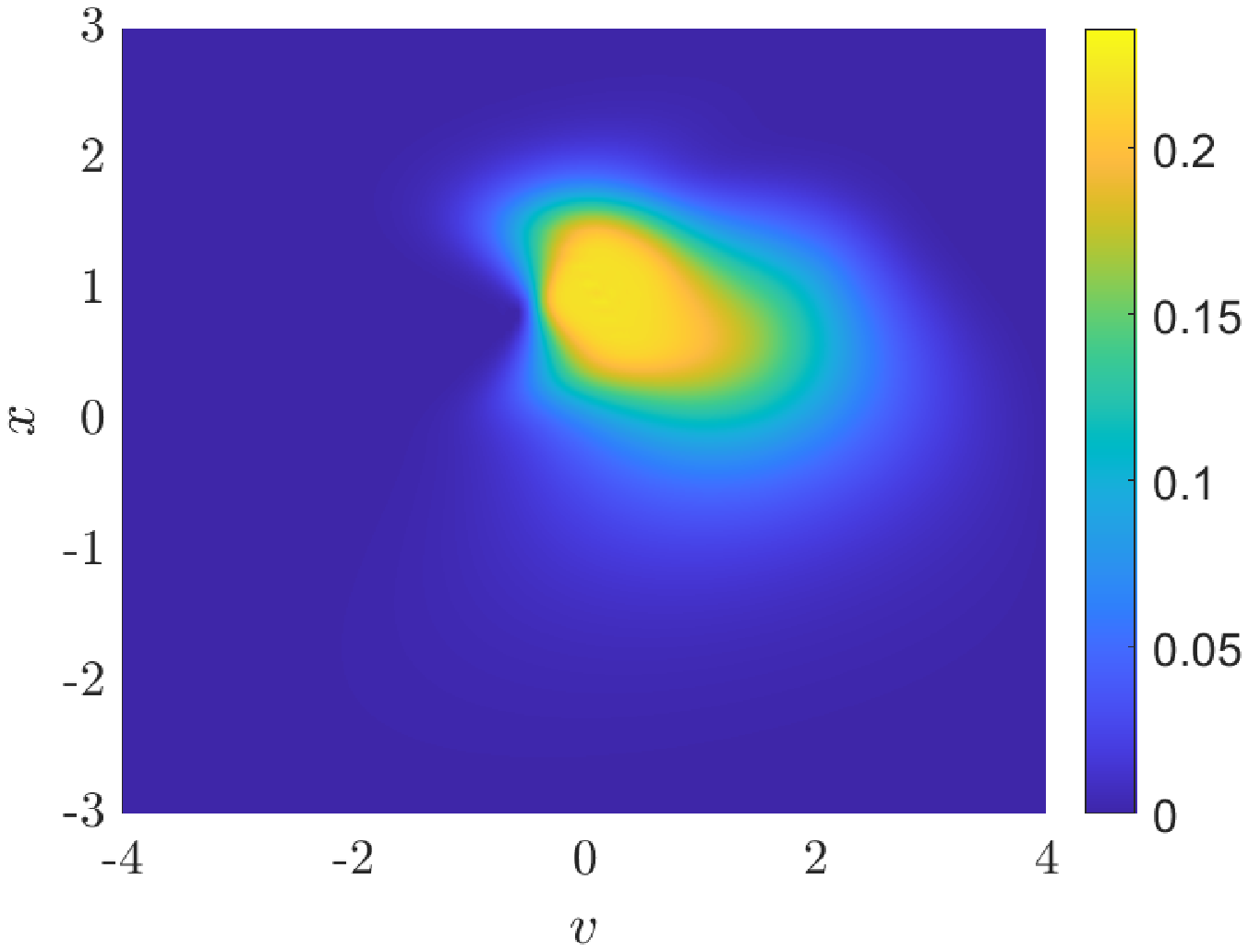}}
\subcaptionbox{Mean-field solution, $t = 3$}{\includegraphics[scale=0.35]{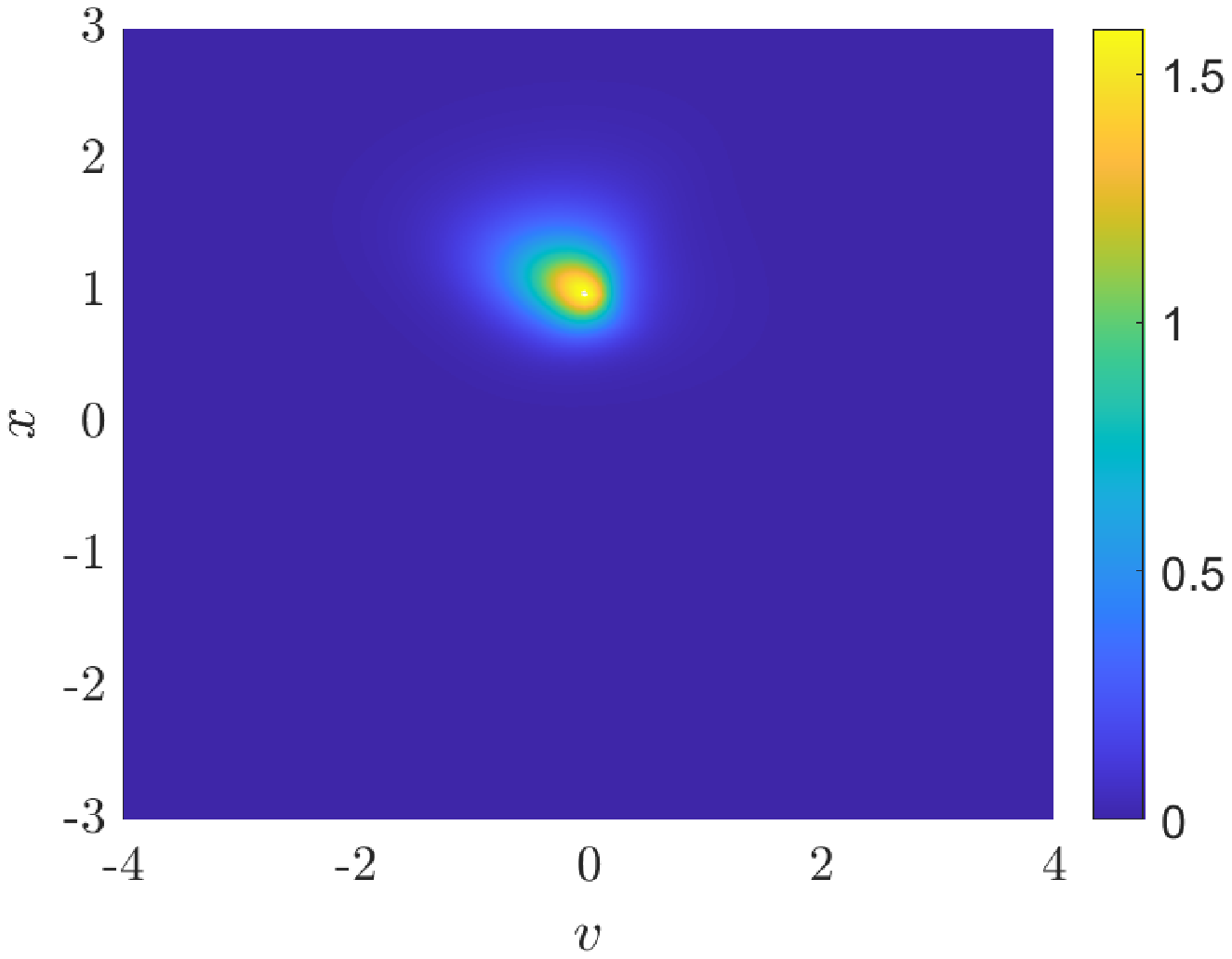}} 
\caption{Case \#1 (no memory). Optimization of the one-dimensional Ackley function with minimum in $x=1$. First row: solution of the SD-PSO system \eqref{eq:psociri}. Second row: solution of the MF-PSO limit \eqref{PDEi}.} 
\label{Fig5}
\end{minipage}
\vspace{5pt}
\\
\begin{minipage}{\linewidth}
\centering
\subcaptionbox{$\rho(x,t)$, $t = 0.5$}{\includegraphics[scale=0.35]{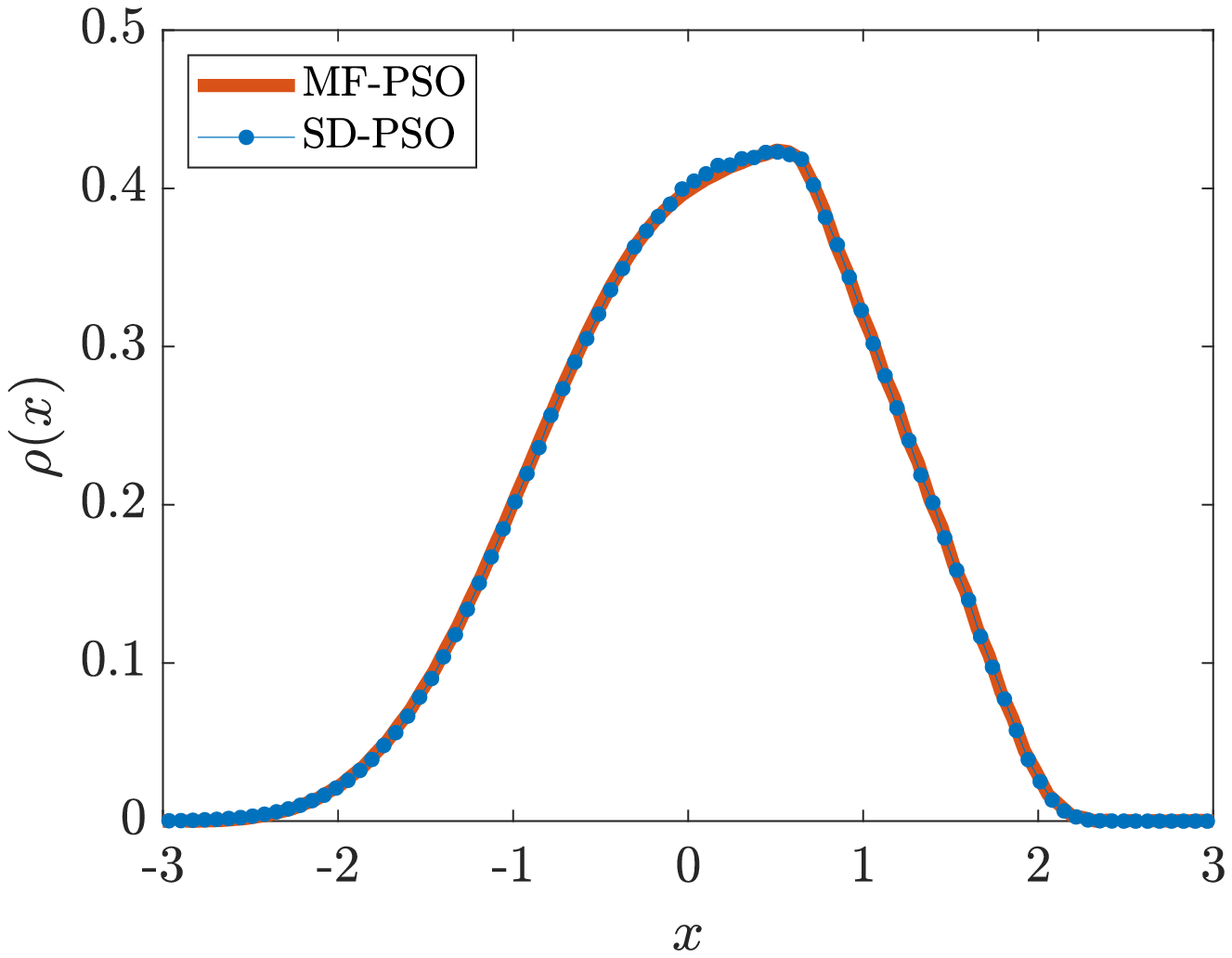}}\
\subcaptionbox{$\rho(x,t)$, $t = 1$}{\includegraphics[scale=0.35]{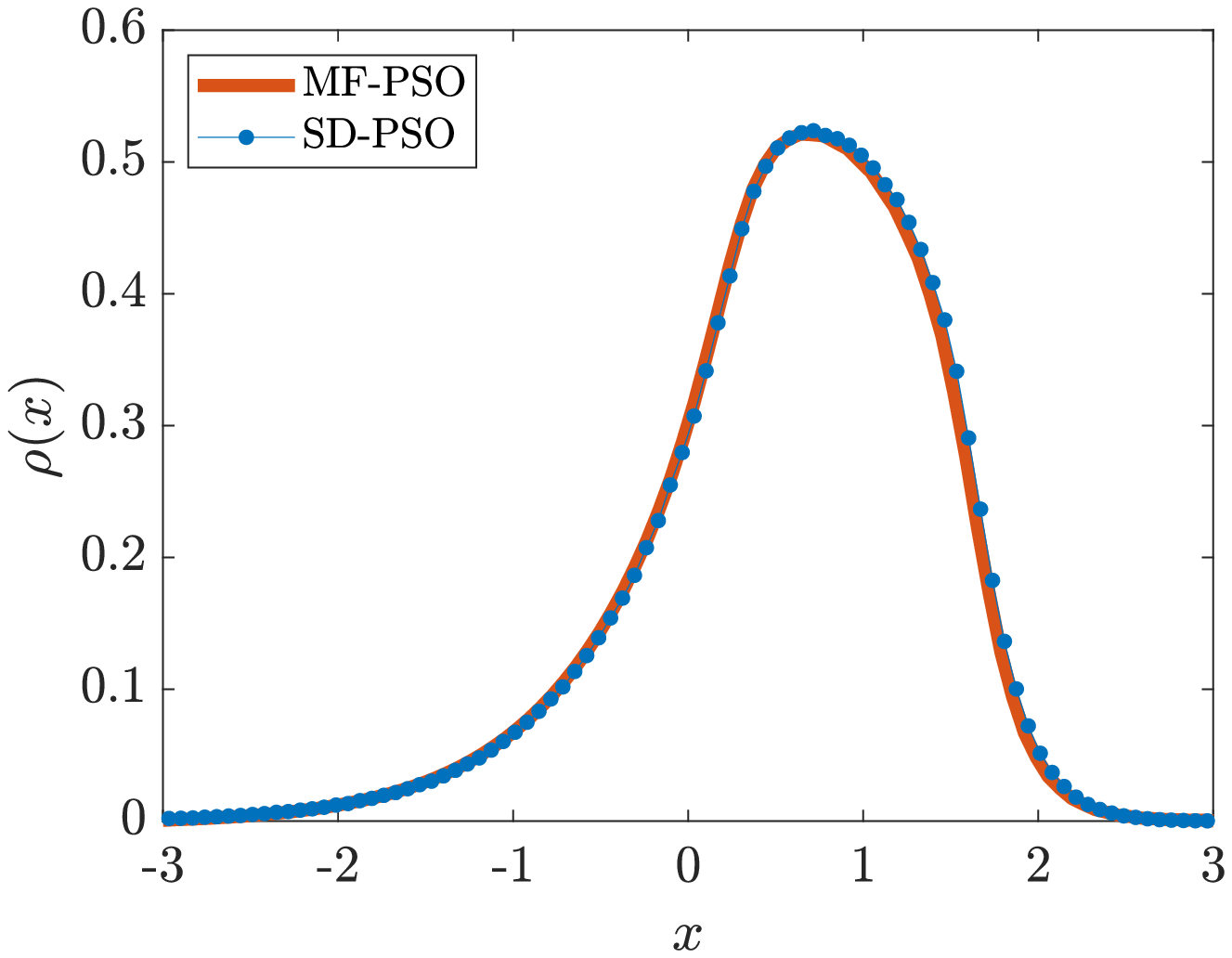}} \
\subcaptionbox{$\rho(x,t)$, $t = 3$}{\includegraphics[scale=0.35]{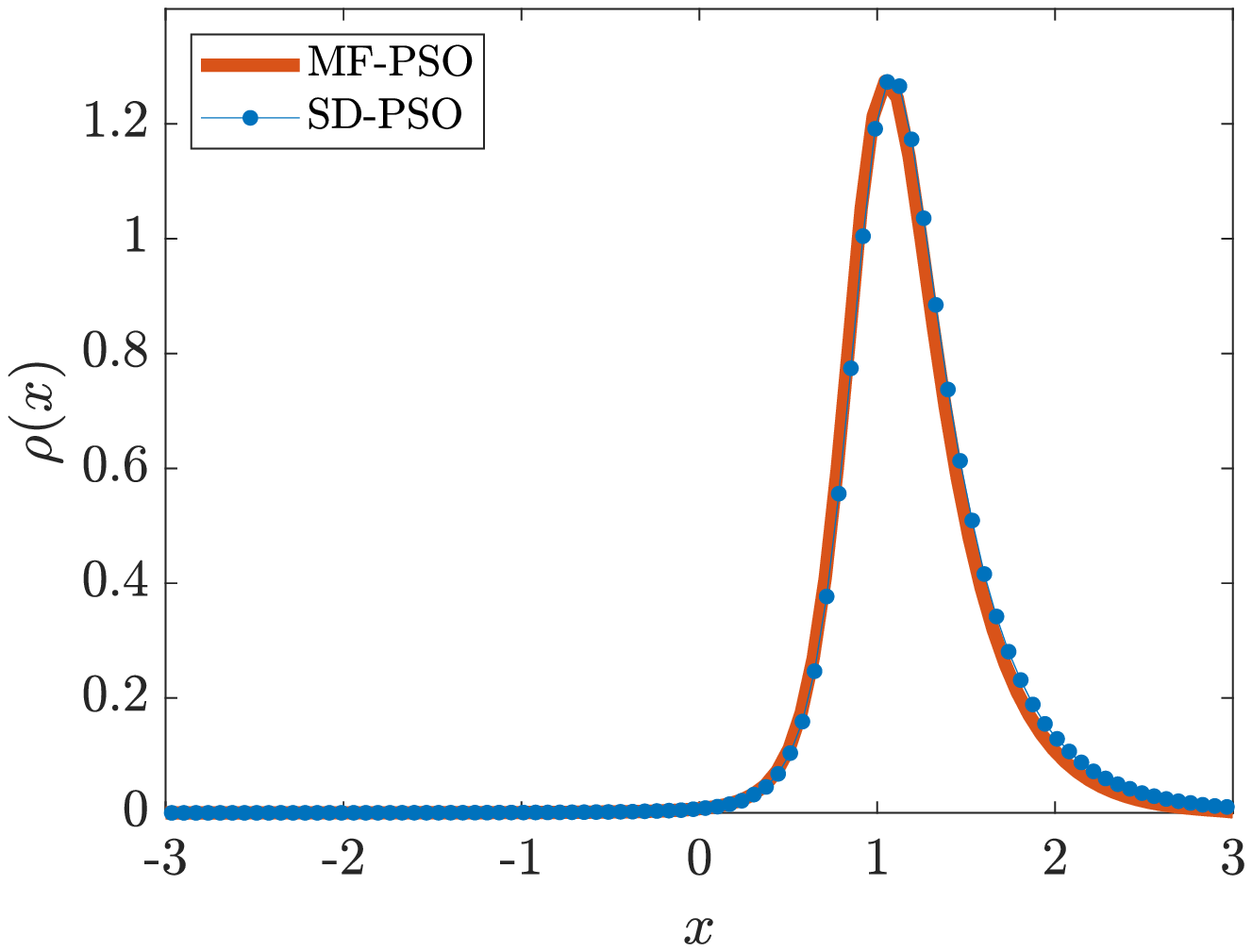}}\\
\caption{Case \#1 (no memory). Evolution of the density $\rho(x,t)$ of the SD-PSO system \eqref{eq:psociri} and the MF-PSO limit \eqref{PDEi} for the one-dimensional Ackley function with minimum in $x=1$.}
\label{Fig6}
\end{minipage}
\end{figure}

In Figures \ref{Fig5} and \ref{Fig6} we report the same results but applied to the Ackley function with minimum in $x = 1$. Even in this case the plot of the density in Figure \ref{Fig6} show the excellent agreement between the stochastic differential system and the mean-field limit. It is interesting to note that the particle distribution evolves asymmetrically in this case and initially exceed the value $ x = 1 $ before moving backward to reach the global minimum. 

\subsubsection*{Case \#2: MF-PSO with memory and only local best dynamics}
In the second test case we introduce the dependence from the memory variable and report a comparison between the solution of the discretized stochastic particle model \eqref{eq:psoDiscr} and the solver of the mean field limit \eqref{PDE} where we assume $\lambda_2$, $\sigma_2 = 0$, namely only the local best is present. The same parameters \eqref{eq:paramt} have been used together with $\beta = 30$ and $\nu = 0.5$ for the local best.
Initially the local best values are assumed to be equal to the particle positions. 

\begin{figure}[H] 
\begin{minipage}{\linewidth}
\centering
\subcaptionbox{Particle solution, $t = 0.5$}{\includegraphics[scale=0.35]{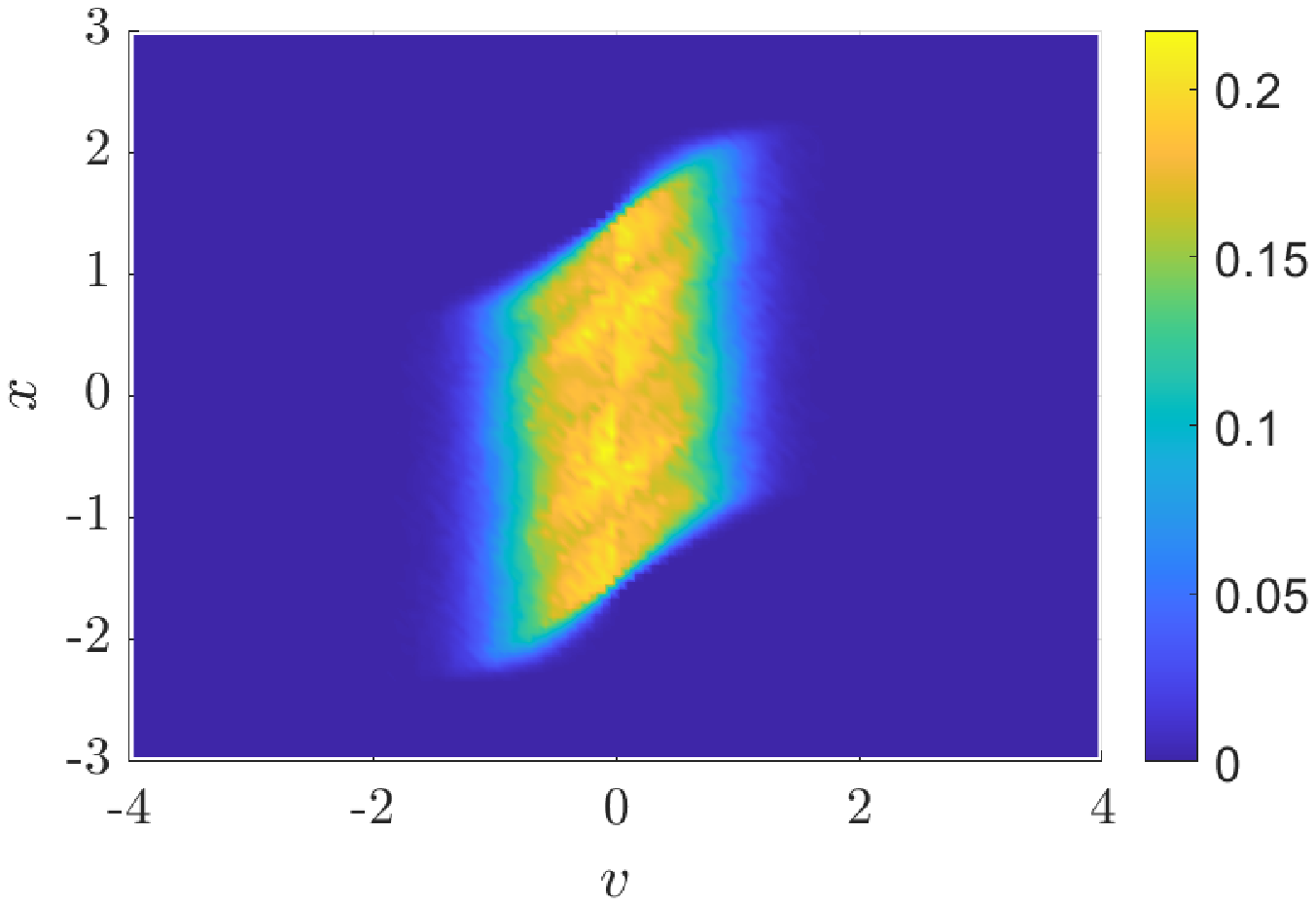}}
\subcaptionbox{Particle solution, $t = 3$}{\includegraphics[scale=0.35]{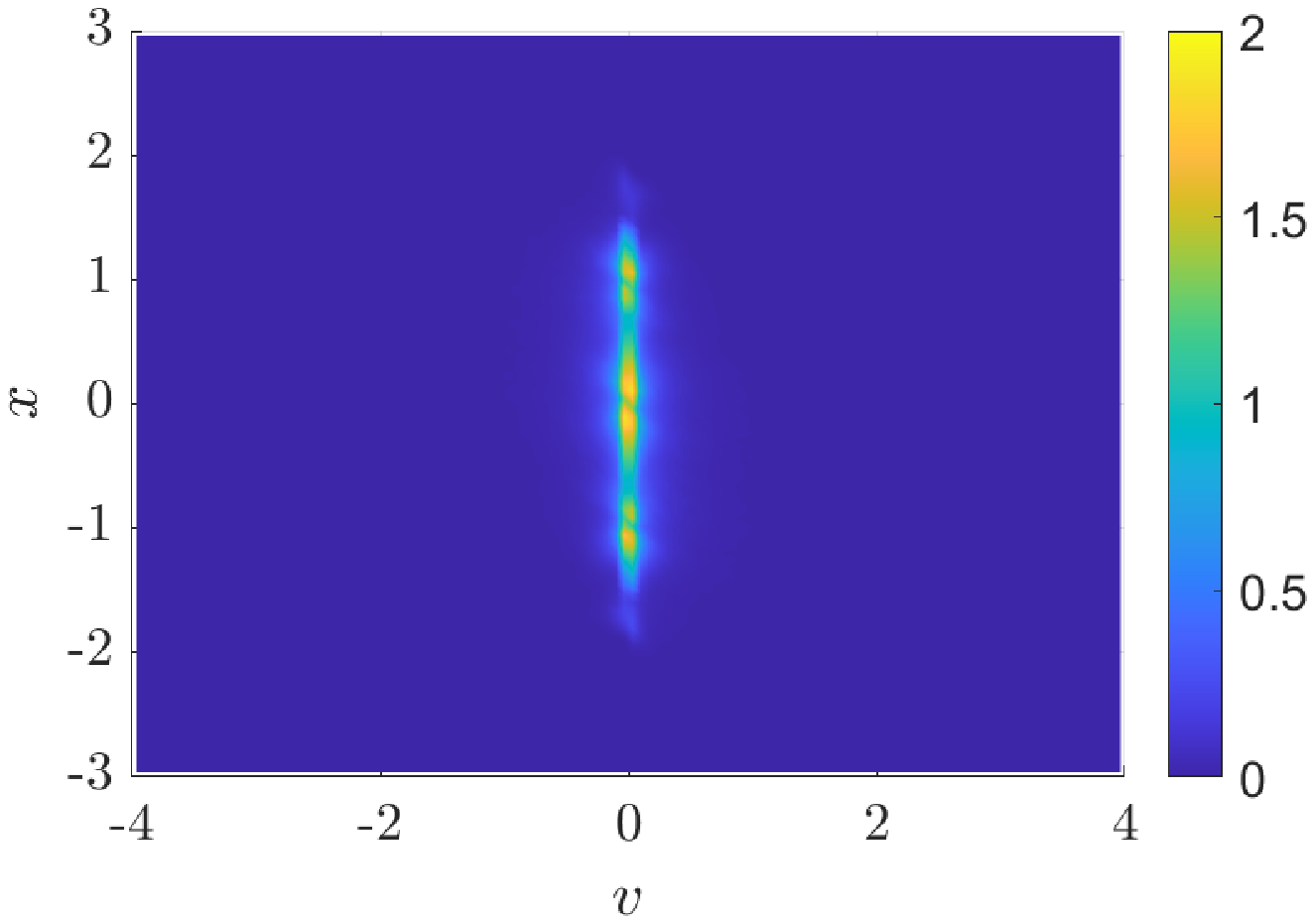}}
\subcaptionbox{Particle solution, $t = 6$}{\includegraphics[scale=0.35]{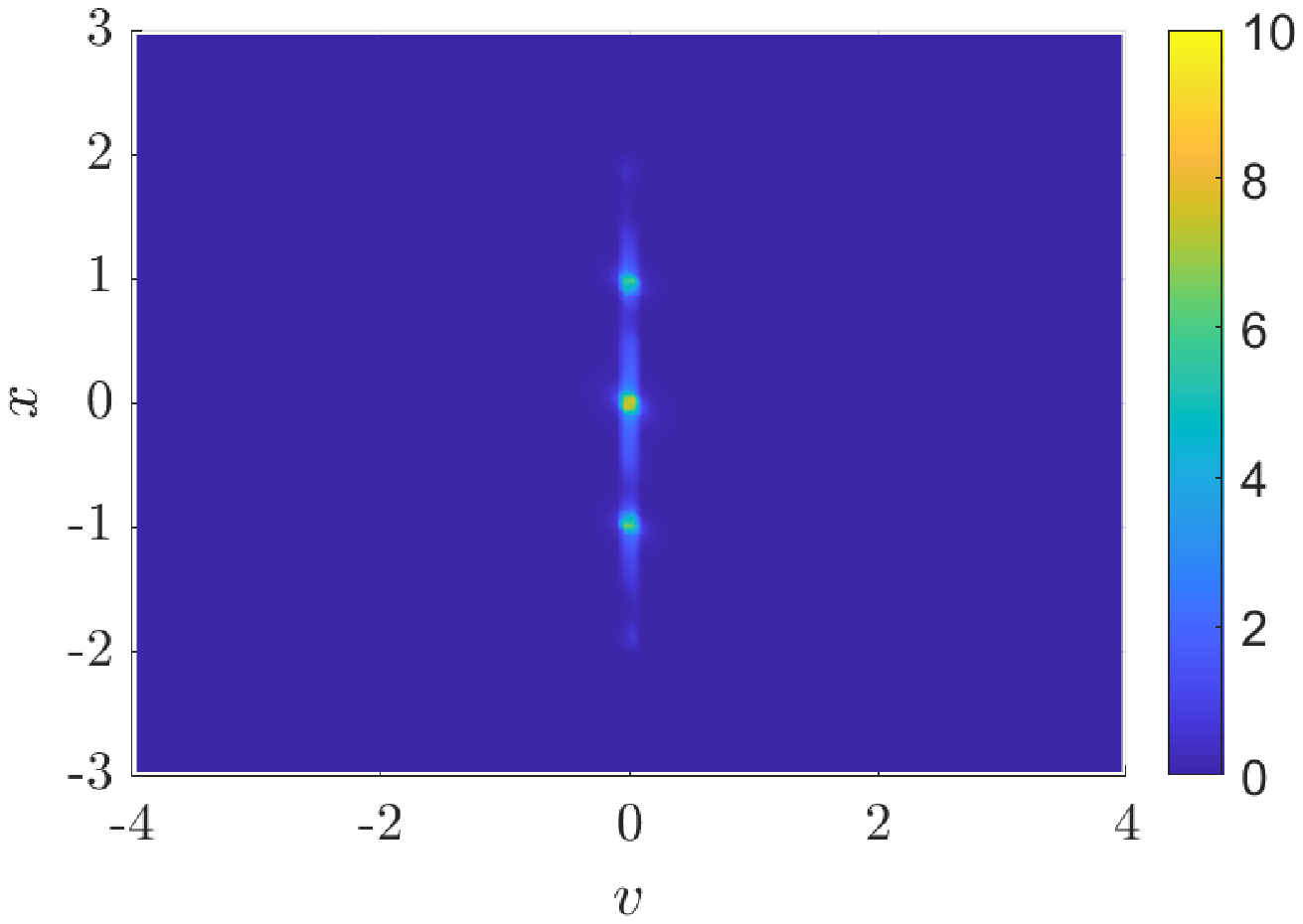}}\\
\subcaptionbox{Marginal $\rho(x,v,t)$, $t = 0.5$}{\includegraphics[scale=0.35]{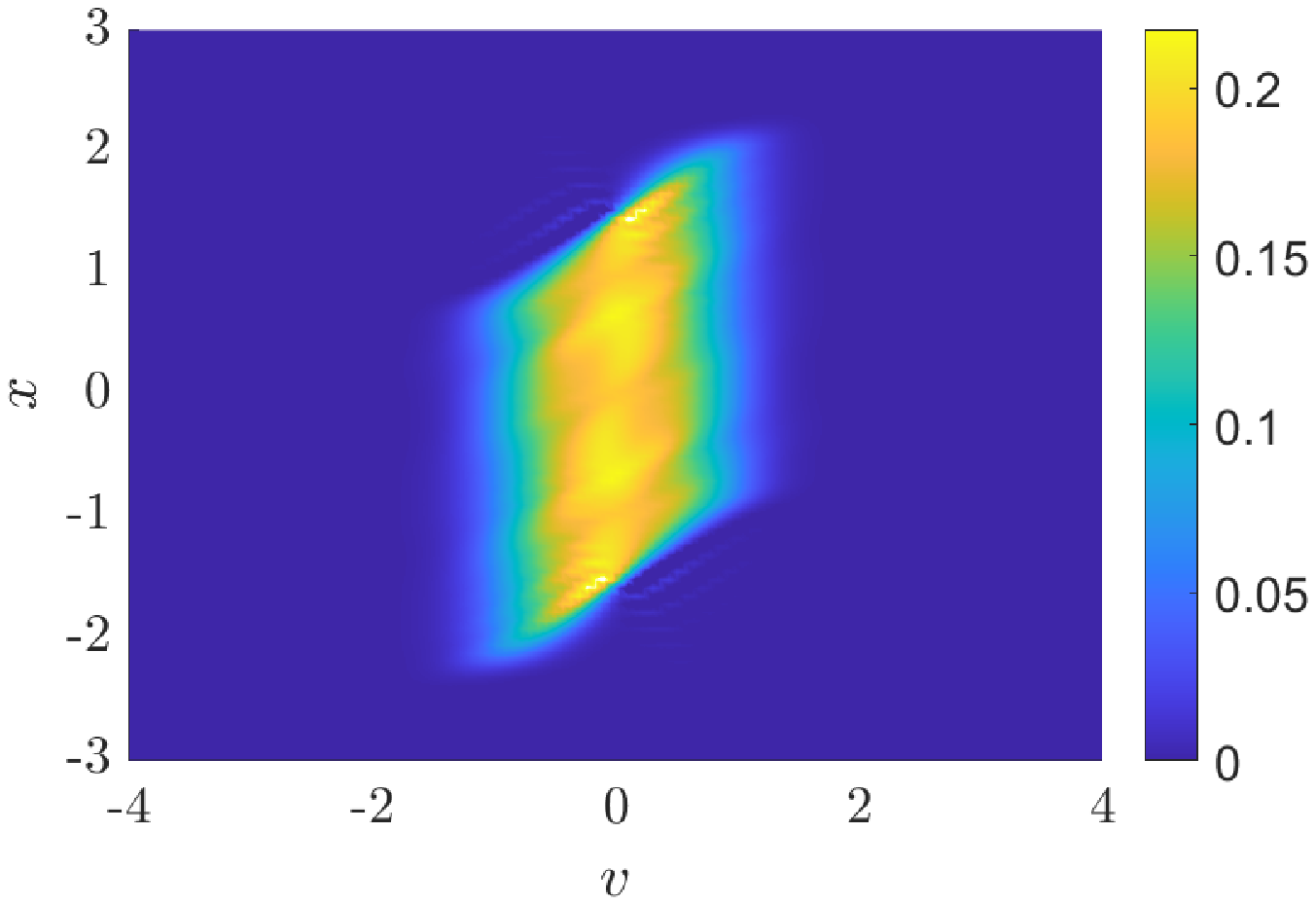}} 
\subcaptionbox{Marginal $\rho(x,v,t)$, $t = 3$}{\includegraphics[scale=0.35]{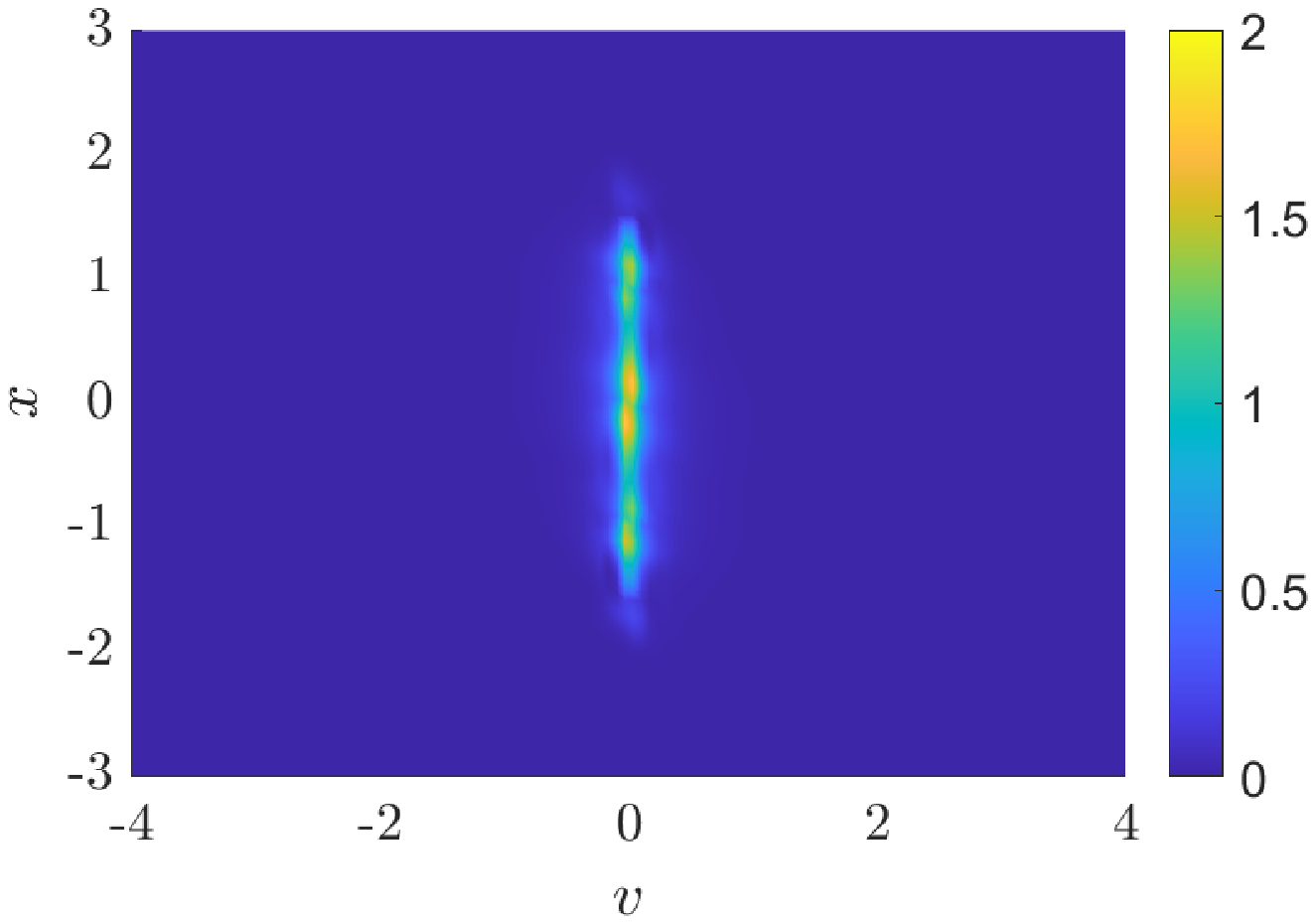}}
\subcaptionbox{Marginal $\rho(x,v,t)$, $t = 6$}{\includegraphics[scale=0.35]{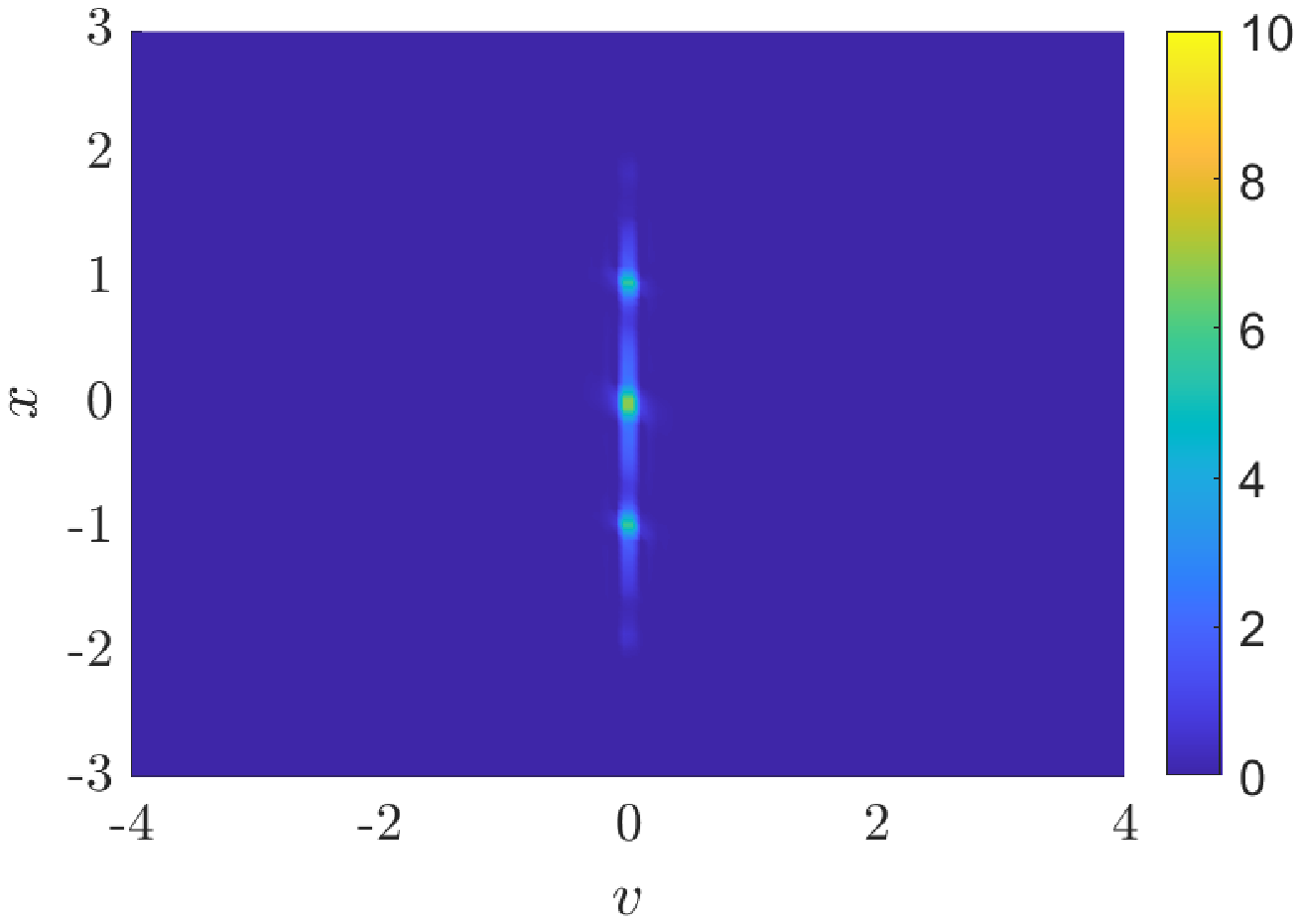}} 
\caption{Case \#2 (local best only). Optimization of the one-dimensional Ackley function with minimum in $x=0$. First row: solution of the SD-PSO system \eqref{eq:psocir}. Second row: solution of the MF-PSO limit \eqref{PDE}.} 
\label{Fig7}
\end{minipage}
\vspace{5pt}
\\
\begin{minipage}{\linewidth}
\centering
\subcaptionbox{$\rho(x,t)$, $t = 0.5$}{\includegraphics[scale=0.35]{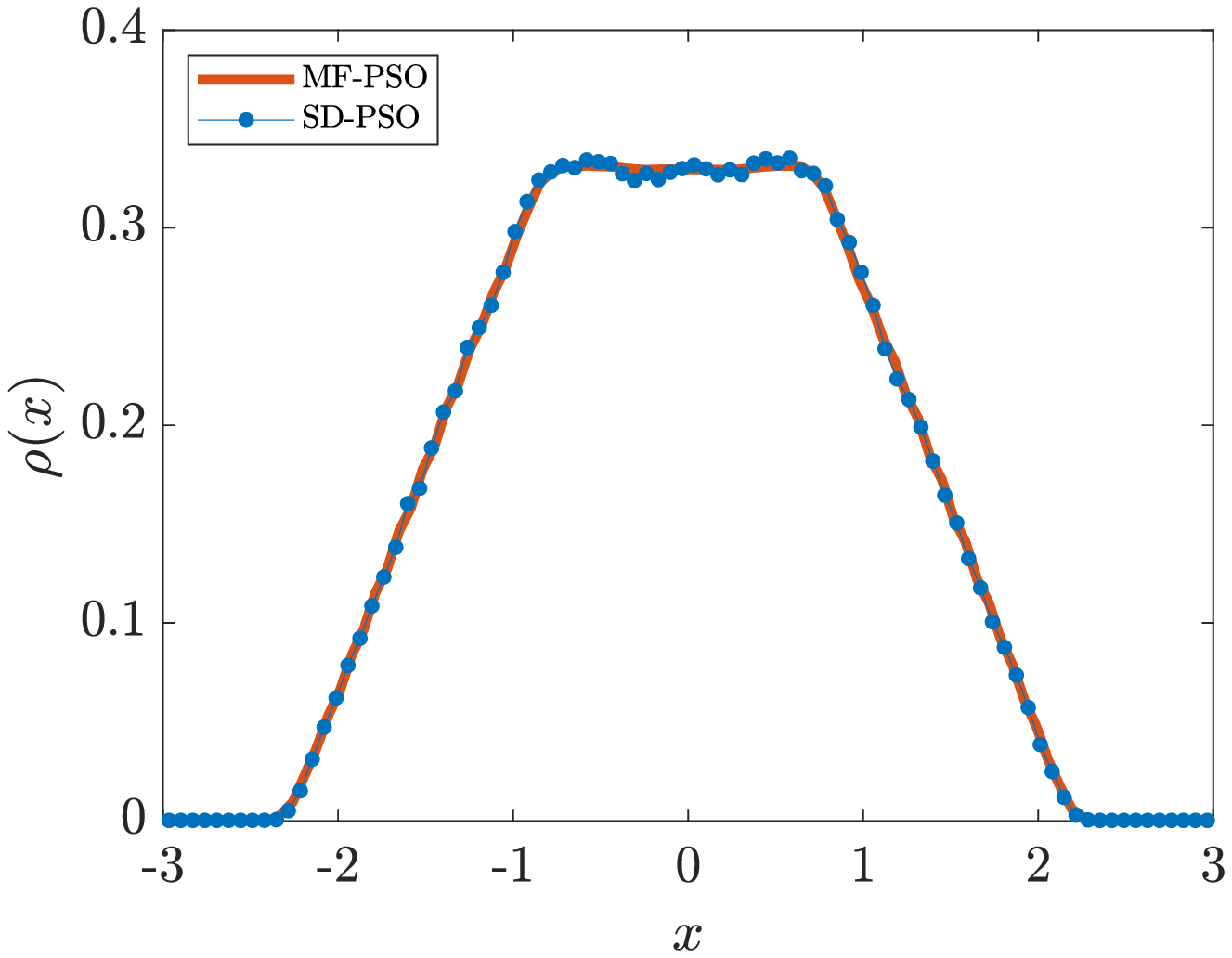}}\
\subcaptionbox{$\rho(x,t)$, $t = 3$}{\includegraphics[scale=0.35]{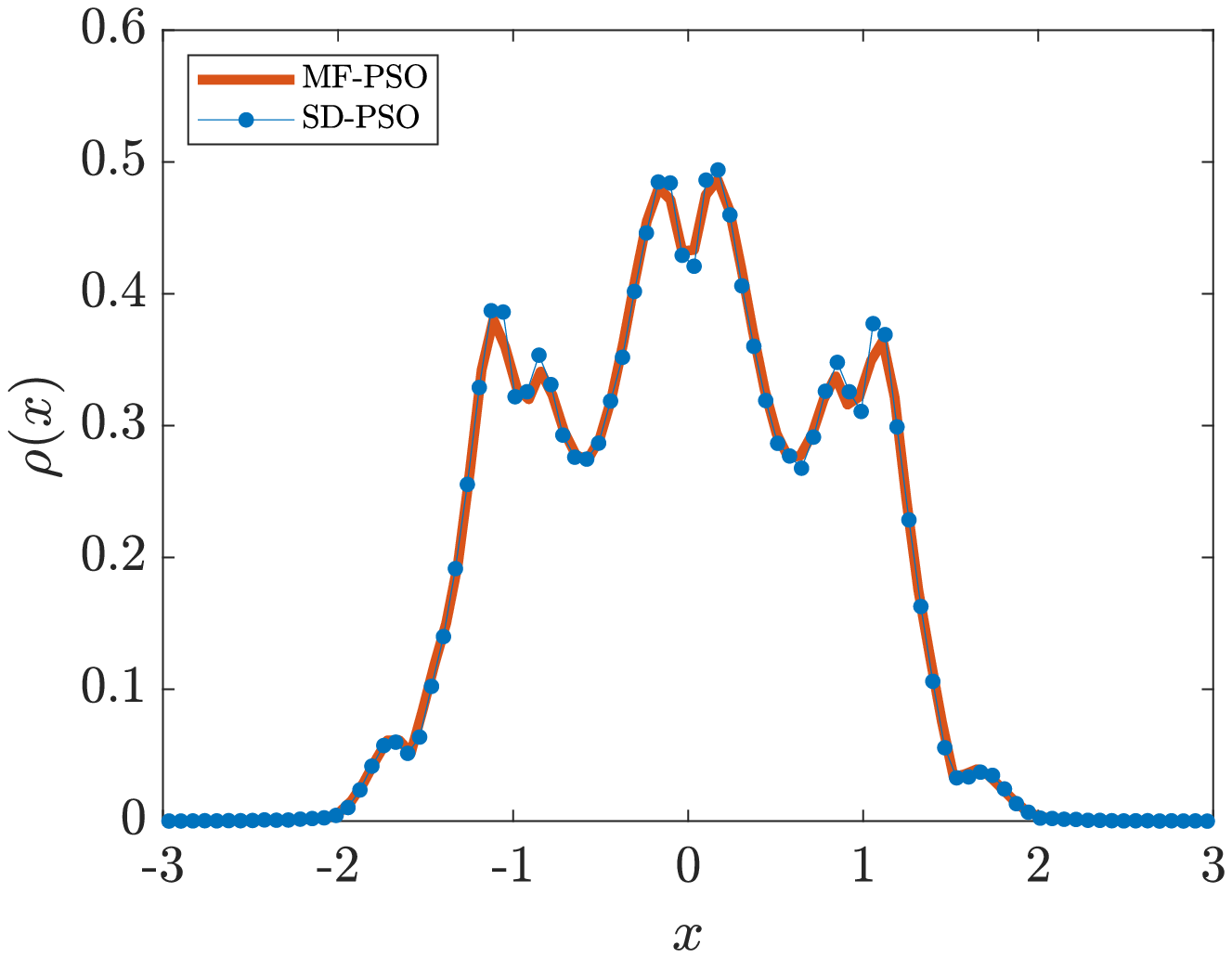}} \
\subcaptionbox{$\rho(x,t)$, $t = 6$}{\includegraphics[scale=0.35]{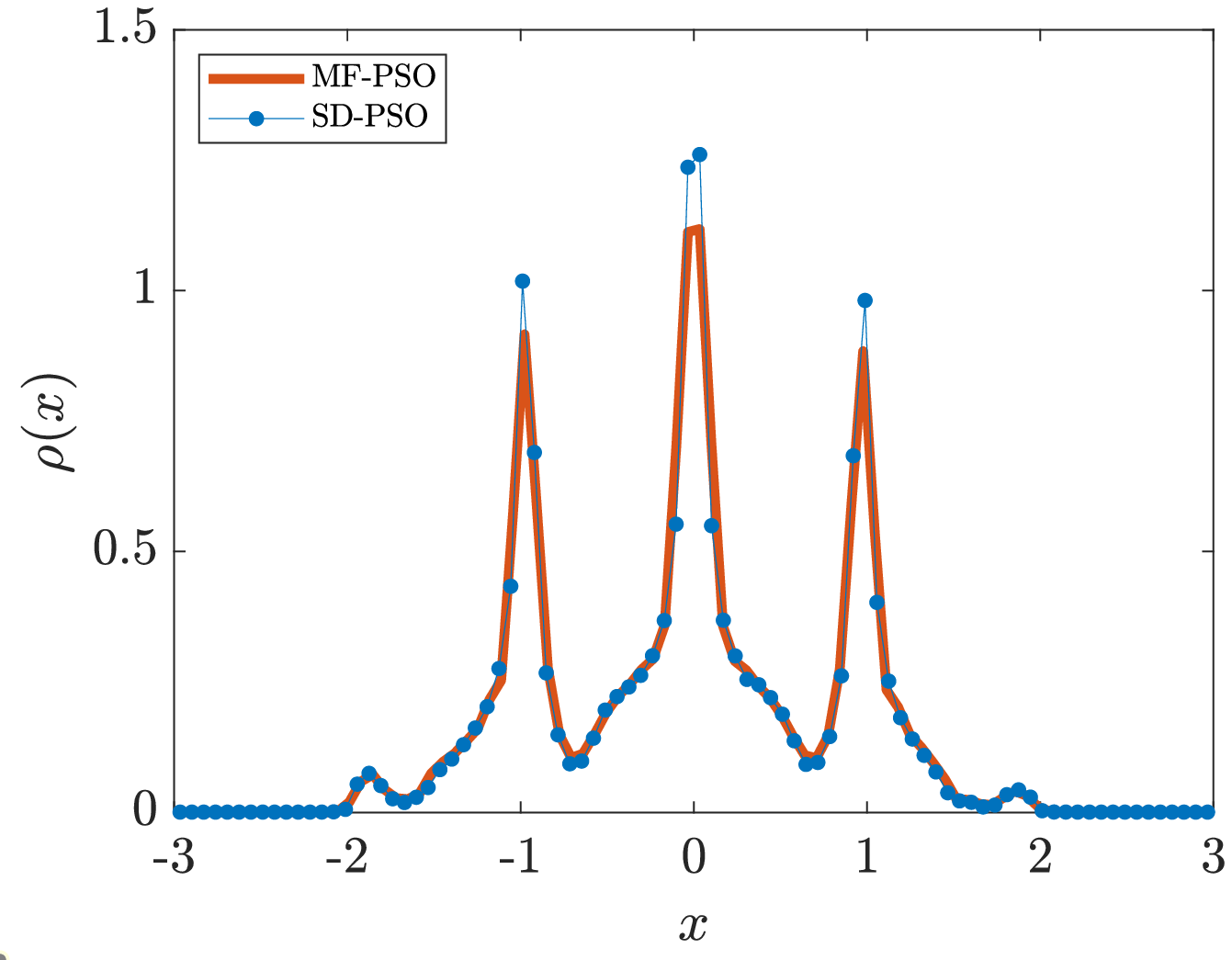}} \\
\caption{Case \#2 (local best only). Evolution of the density $\rho(x,t)$ of the SD-PSO system \eqref{eq:psocir} and the MF-PSO limit \eqref{PDE} for the one-dimensional Ackley function with minimum in $x=0$. } 
\label{Fig8}
\end{minipage}
\end{figure}

\begin{figure}[H] 
\begin{minipage}{\linewidth}
\centering
\subcaptionbox{Particle solution, $t = 0.5$}{\includegraphics[scale=0.35]{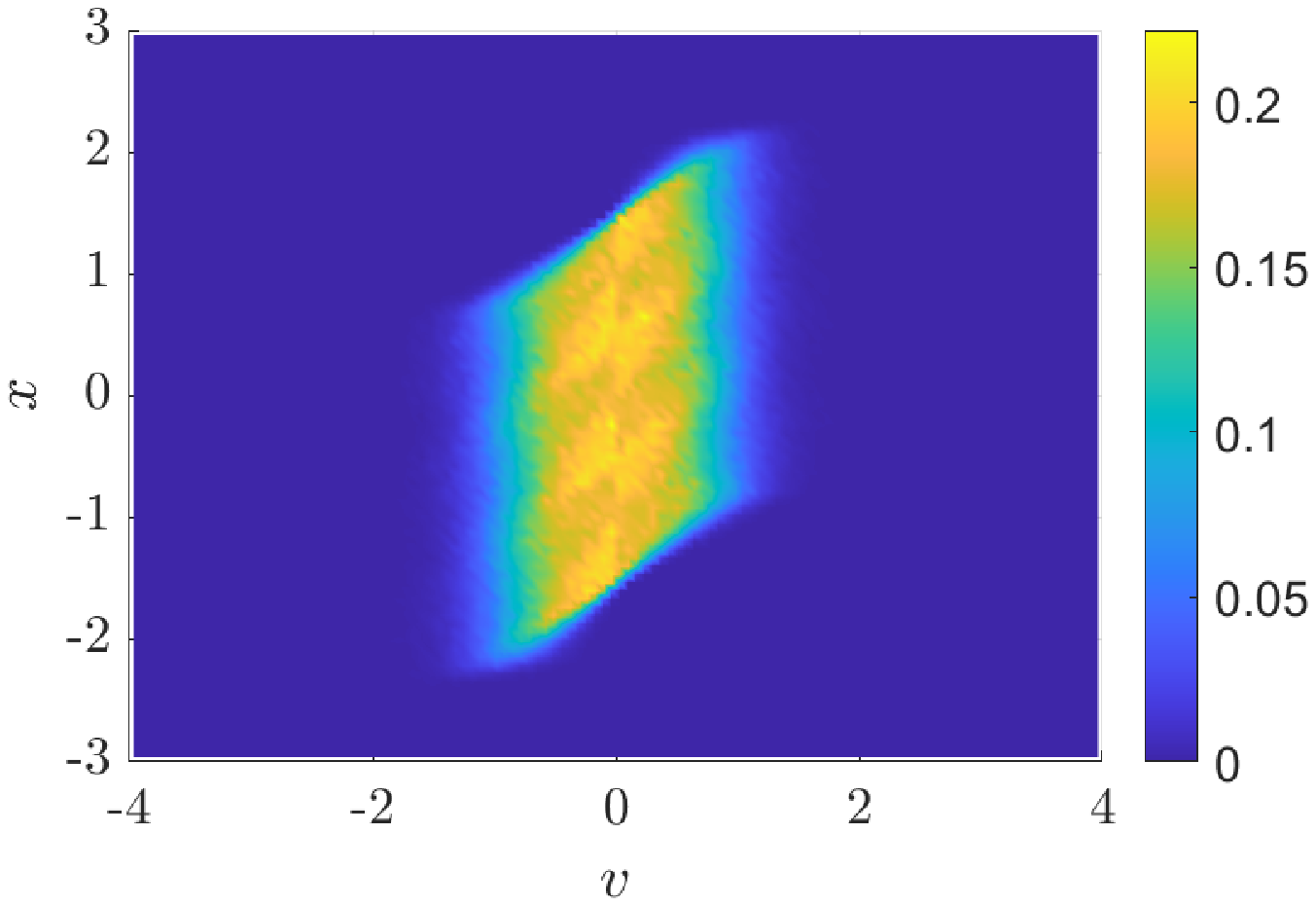}}
\subcaptionbox{Particle solution, $t = 3$}{\includegraphics[scale=0.35]{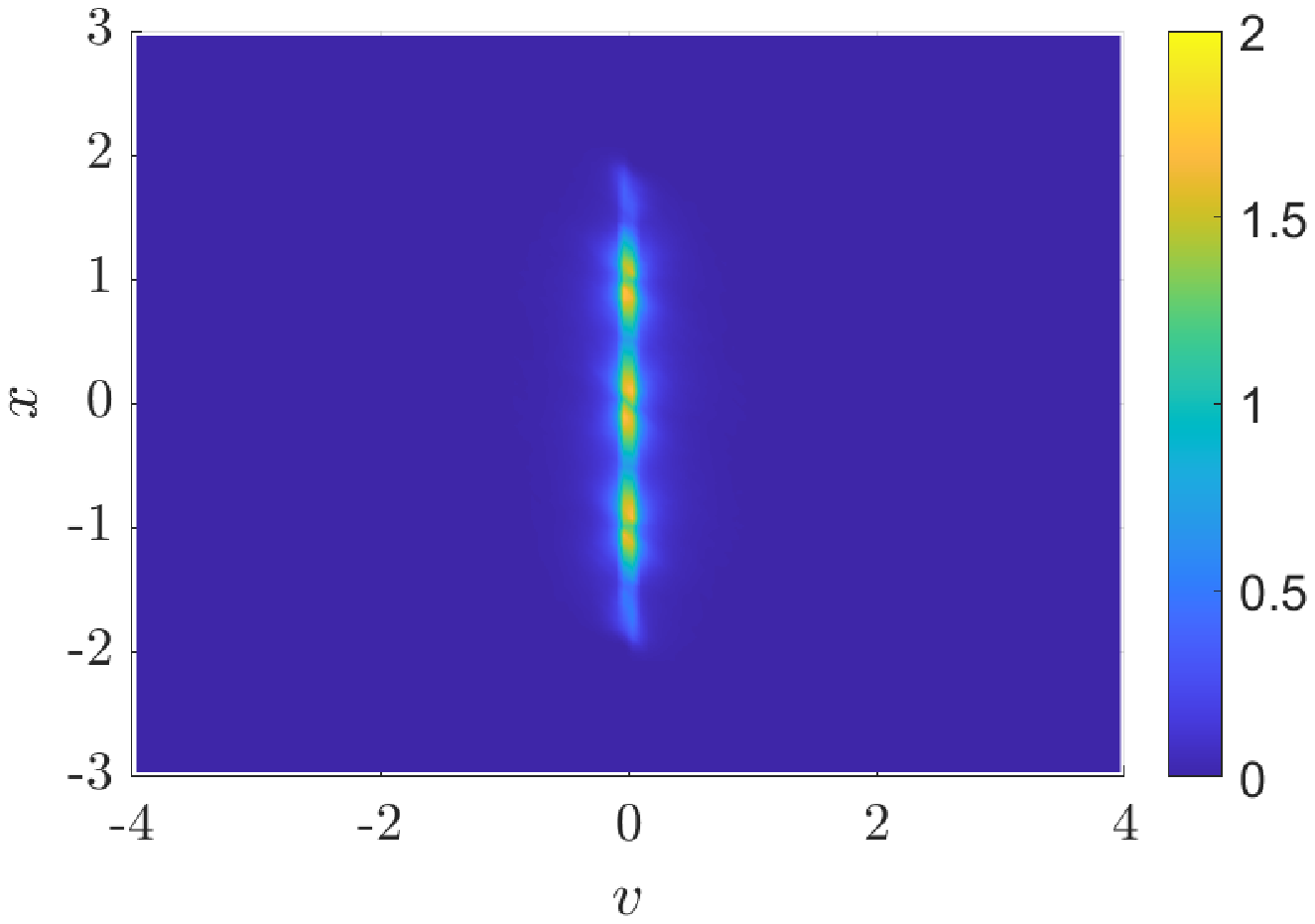}}
\subcaptionbox{Particle solution, $t = 6$}{\includegraphics[scale=0.35]{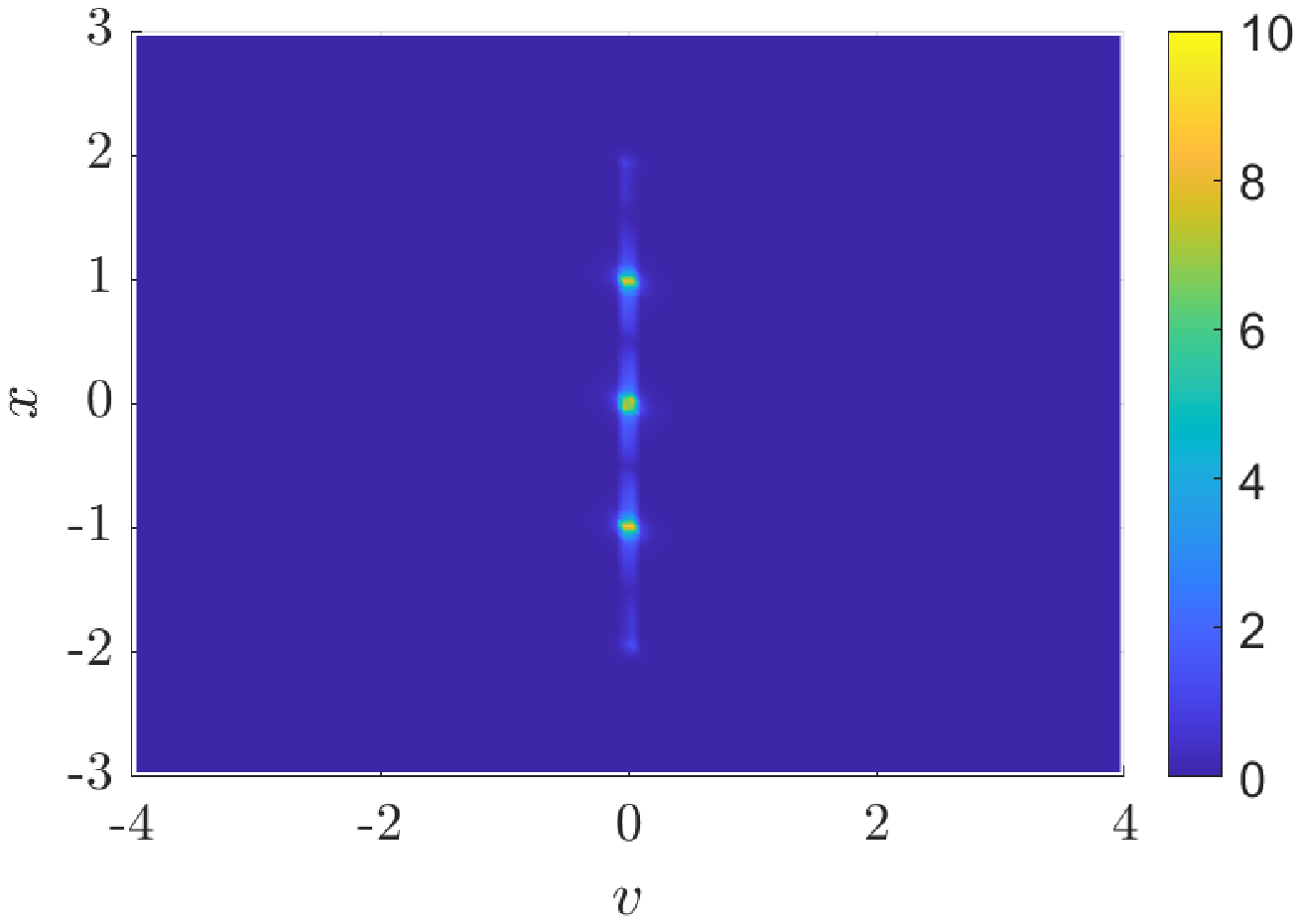}}\\
\subcaptionbox{Marginal $\rho(x,v,t)$, $t = 0.5$}{\includegraphics[scale=0.35]{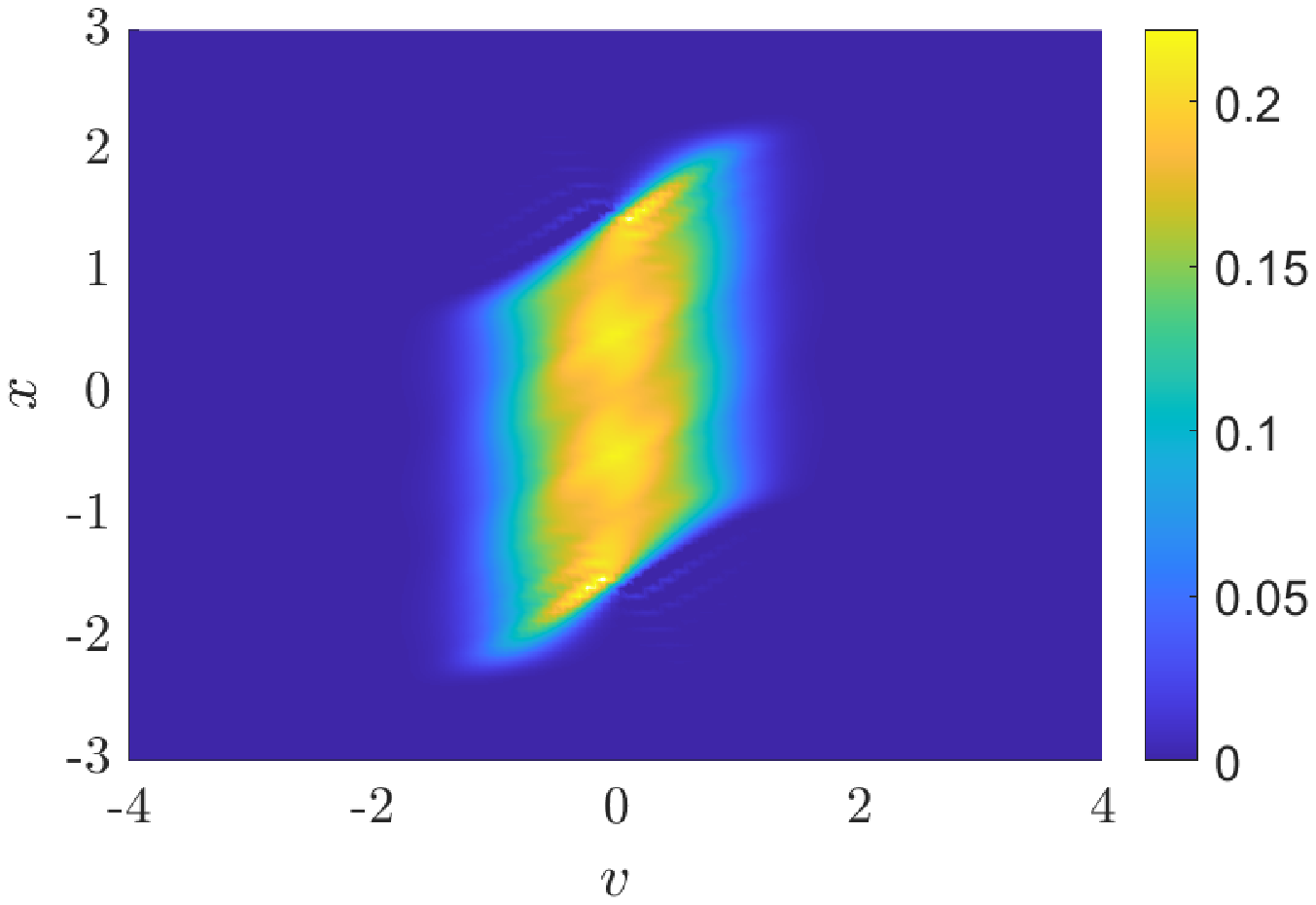}}
\subcaptionbox{Marginal $\rho(x,v,t)$, $t = 3$}{\includegraphics[scale=0.35]{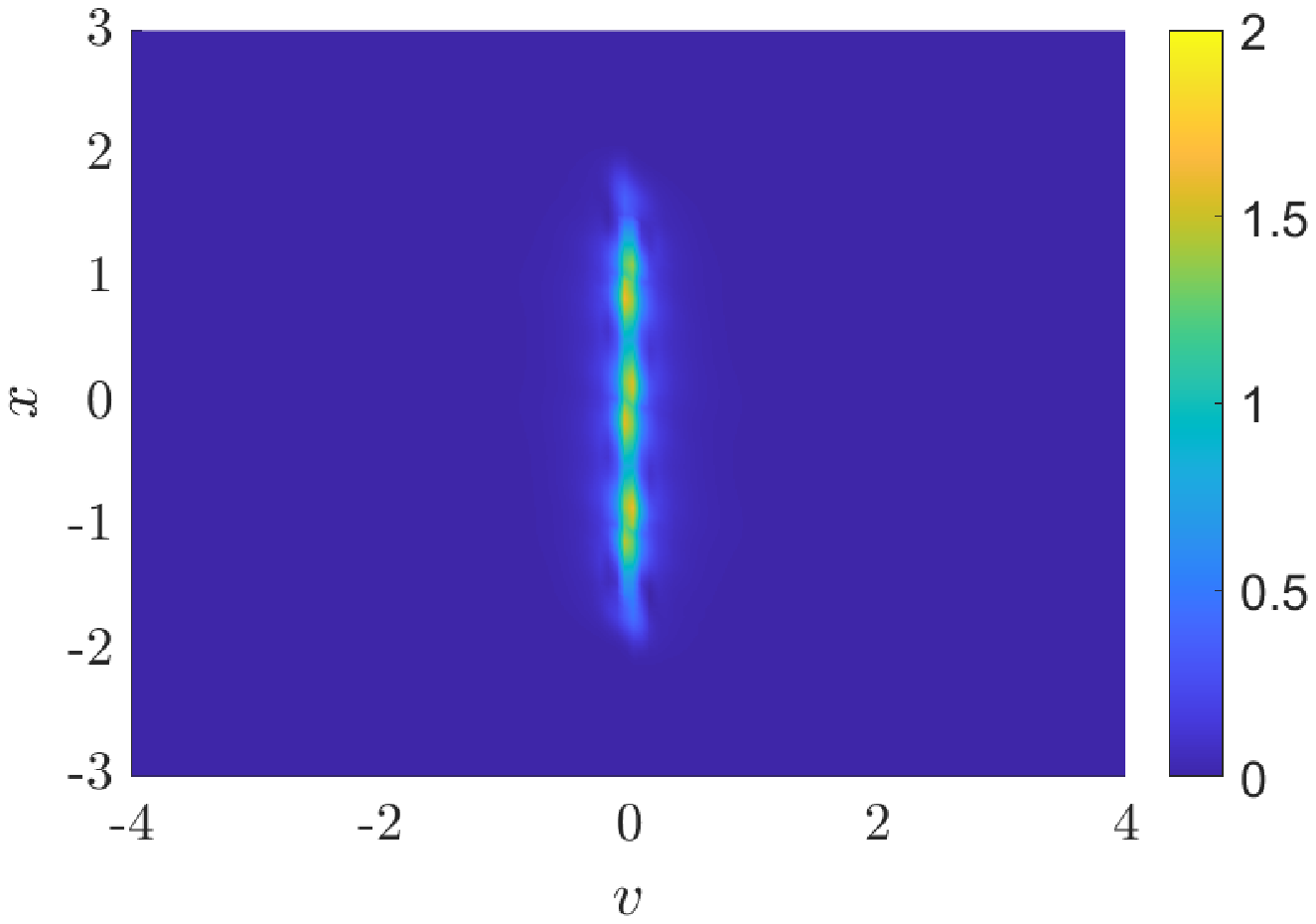}} 
\subcaptionbox{Marginal $\rho(x,v,t)$, $t = 6$}{\includegraphics[scale=0.35]{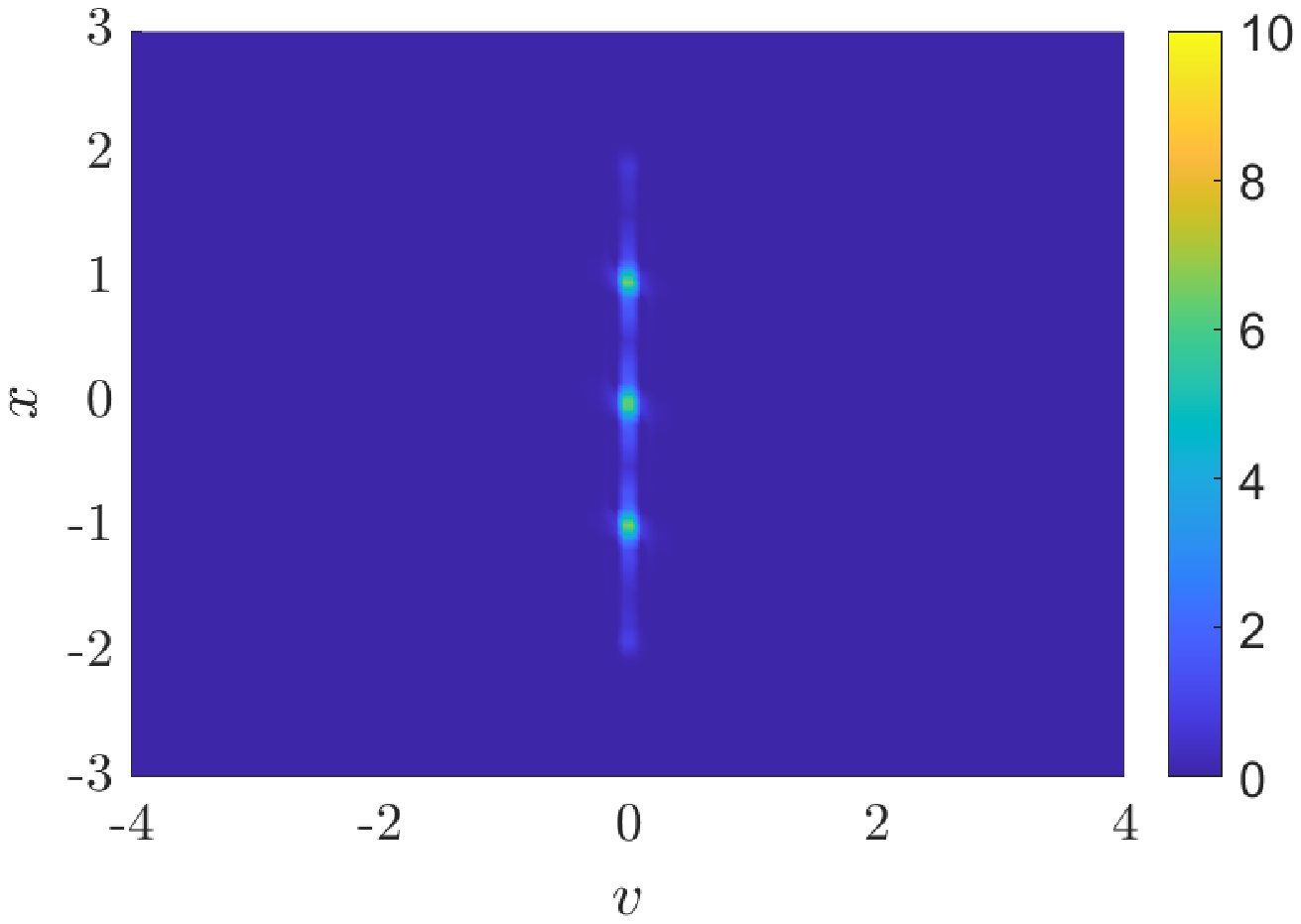}} 
\caption{Case \#2 (local best only). Optimization of the one-dimensional Rastrigin function with minimum in $x=0$. First row: solution of the SD-PSO system \eqref{eq:psocir}. Second row: solution of the MF-PSO limit \eqref{PDE}.} 
\label{Fig9}
\end{minipage}
\vspace{5pt}
\\
\begin{minipage}{\linewidth}
\centering
\subcaptionbox{$\rho(x,t)$, $t = 0.5$}{\includegraphics[scale=0.35]{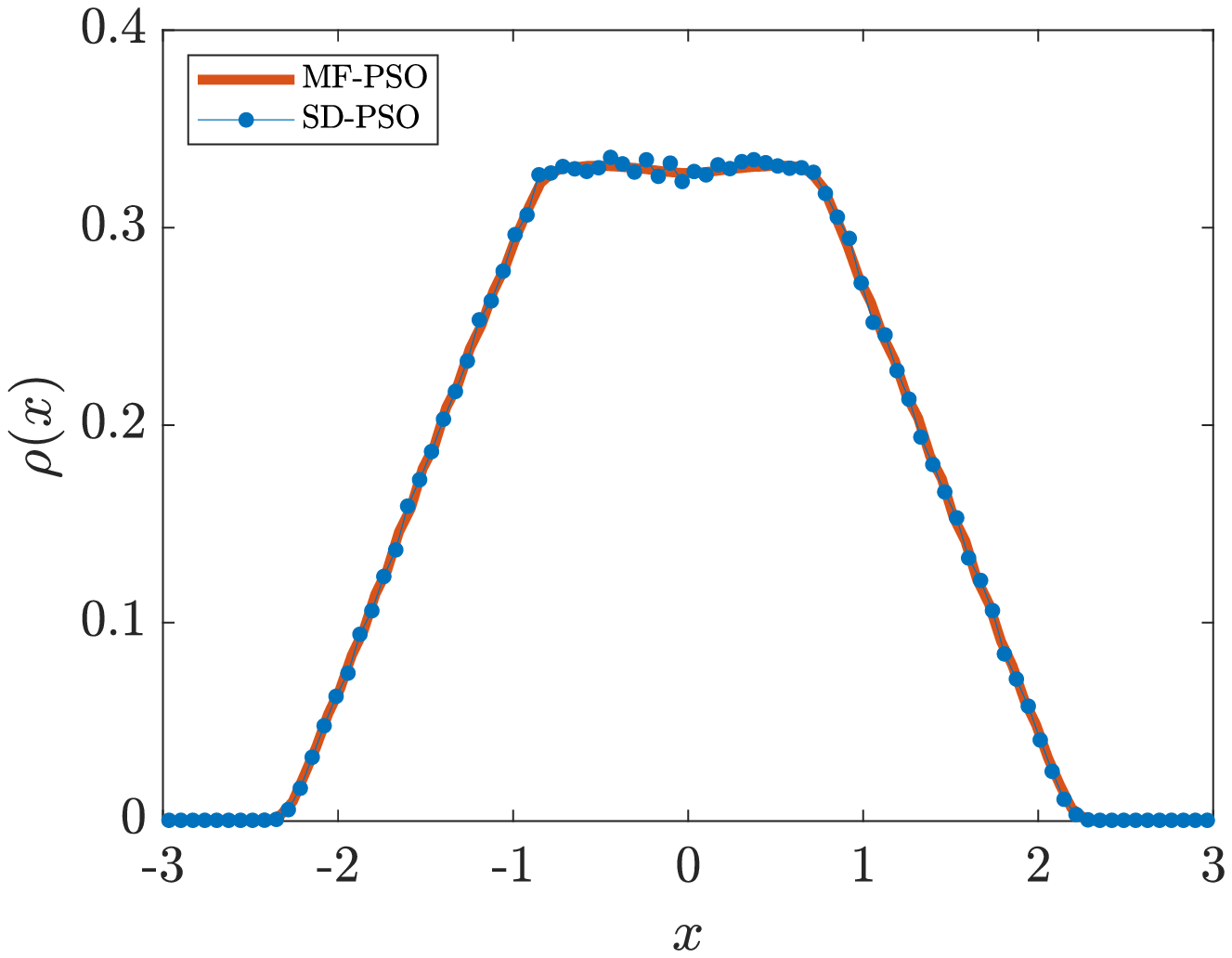}}\ 
\subcaptionbox{$\rho(x,t)$, $t = 3$}{\includegraphics[scale= 0.35]{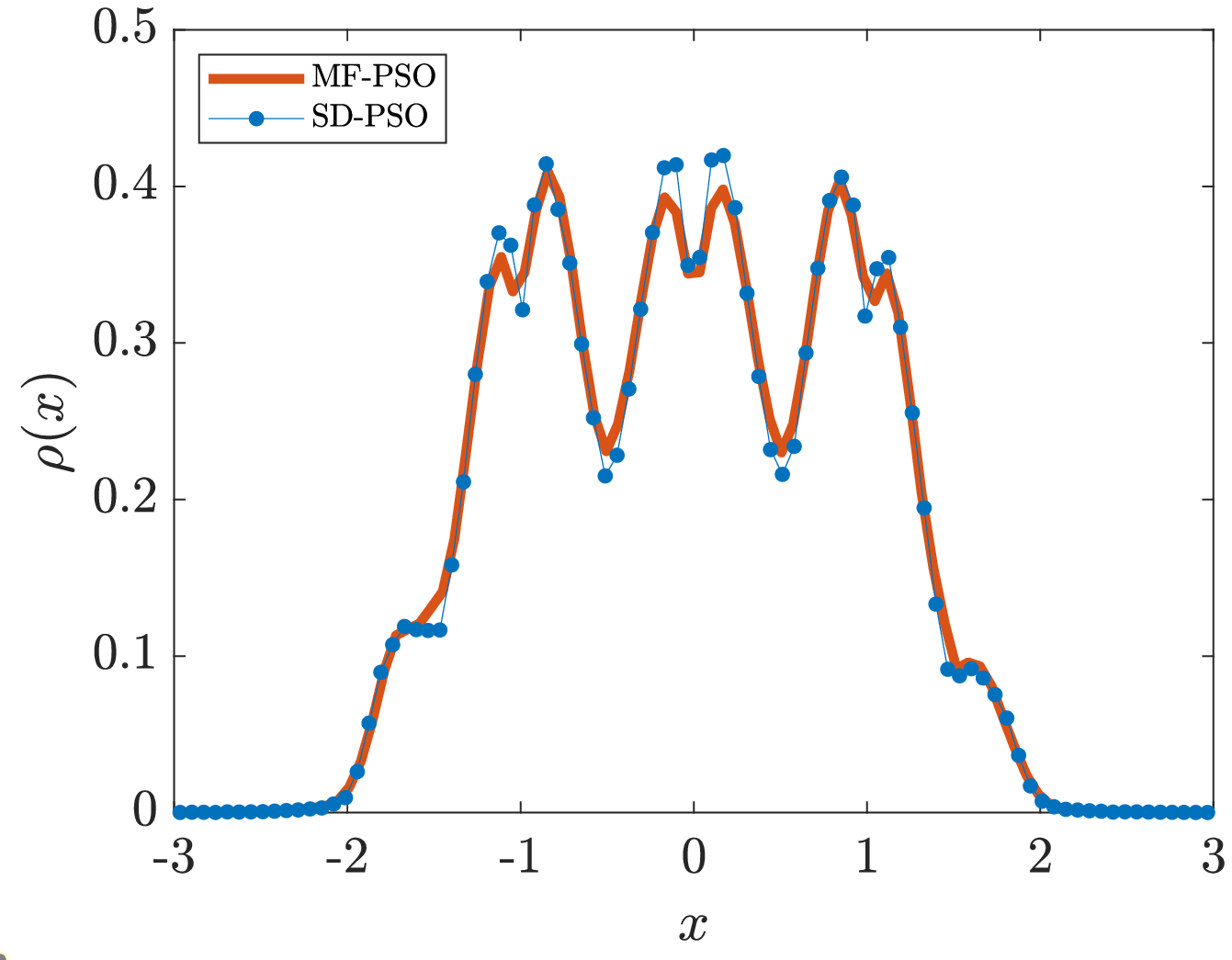}} \ 
\subcaptionbox{$\rho(x,t)$, $t = 6$}{\includegraphics[scale=0.35]{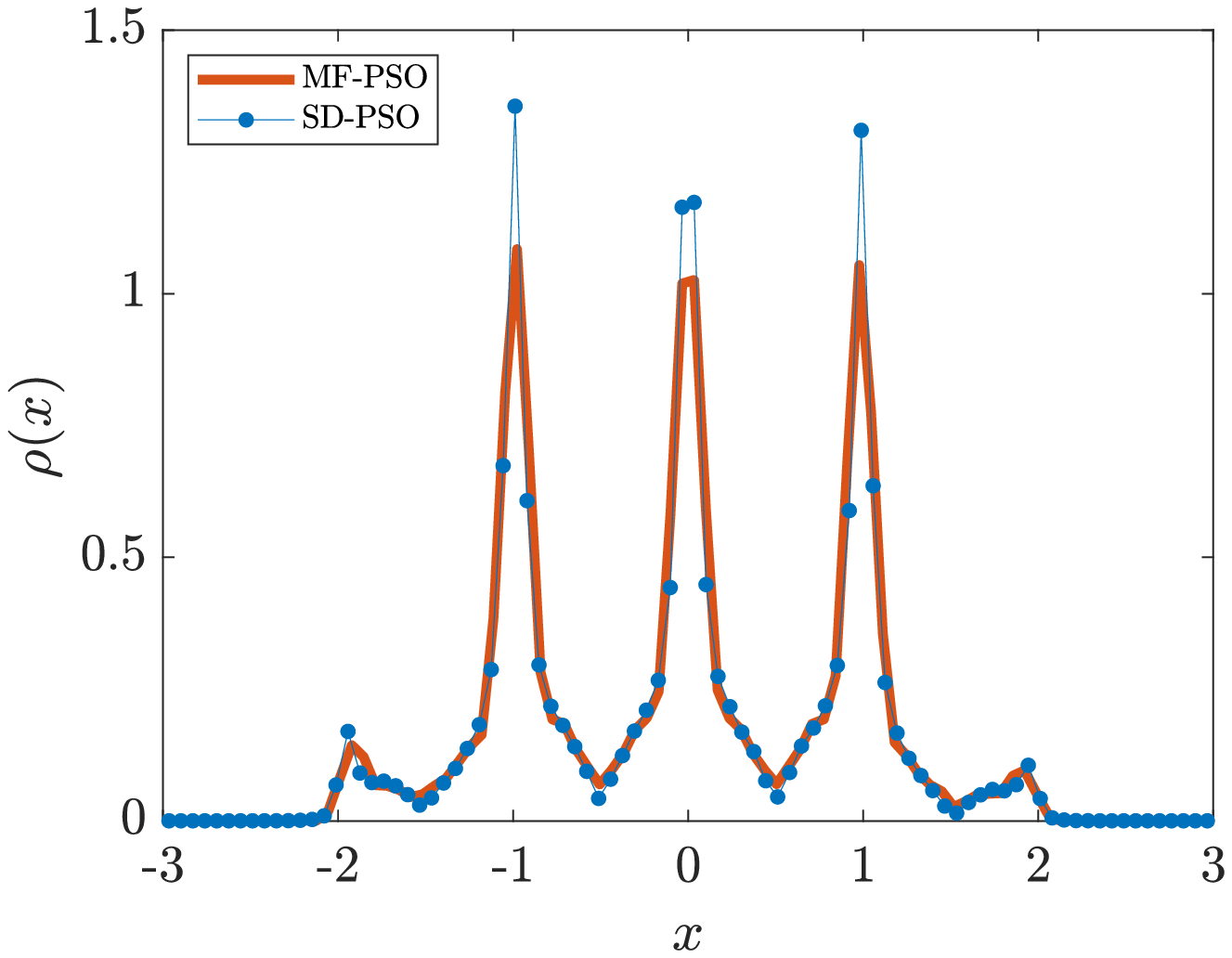}} \\
\caption{Case \#2 (local best only). Evolution of the density $\rho(x,t)$ of the SD-PSO system \eqref{eq:psocir} and the MF-PSO limit \eqref{PDE} for the one-dimensional Rastrigin function with minimum in $x=0$.} 
\label{Fig10}
\end{minipage}
\end{figure}

In Figures \ref{Fig7} and \ref{Fig9} we report the contour plots of the particle solution and the mean-field solution for the one-dimensional Ackley and Rastrigin functions with minimum in $x=0$ and using an uniform initial data. The final simulation time now is $t=6$.
We can note that in the presence of local best only, the particles tend to return to their local best position creating a "memory effect" that leads them to concentrate not only in the global minimum but also in the local minima. For large times we obtain a sequence of particle peaks with zero speed exactly in the positions of the local minima. Thus the dynamic allows us to identify each type of minimum present in the functions. 

This effect can be better appreciated if one compares the explicit form of the one-dimensional Ackley and Rastrigin function in Figure \ref{Fig1} and the corresponding density plots of the particle minimizer in Figures \ref{Fig8} and \ref{Fig10}. It is interesting to point out that from a computational point of view solving the mean field equation \eqref{PDE} in this setting (presence of local best only) has proved to be quite challenging due to the high dimensionality and the importance of avoiding dissipative effects in the discretization of memory terms to preserve the peaks structure in the asymptotic numerical solution. Second order schemes for the discretization of the mean-field equation are essential in this case to resolve correctly the structure of the solution. One can appreciate the good agreement between the particle and mean-field solutions in Figures \ref{Fig8} and \ref{Fig10}.

\subsubsection*{Case \#3: MF-PSO with memory, general case}
In the last test case we repeat the previous scenario by adding the action of the global best with the same weight as the local best. Therefore, we take $\lambda_1=\lambda_2=1$, $\sigma_1=\sigma_2={1}/{\sqrt{3}}$ and the same parameters in \eqref{eq:paramt} in our numerical experiments.

\begin{figure}[H] 
\begin{minipage}{\linewidth}
\centering
\subcaptionbox{Particle solution, $t = 0.5$}{\includegraphics[scale=0.35]{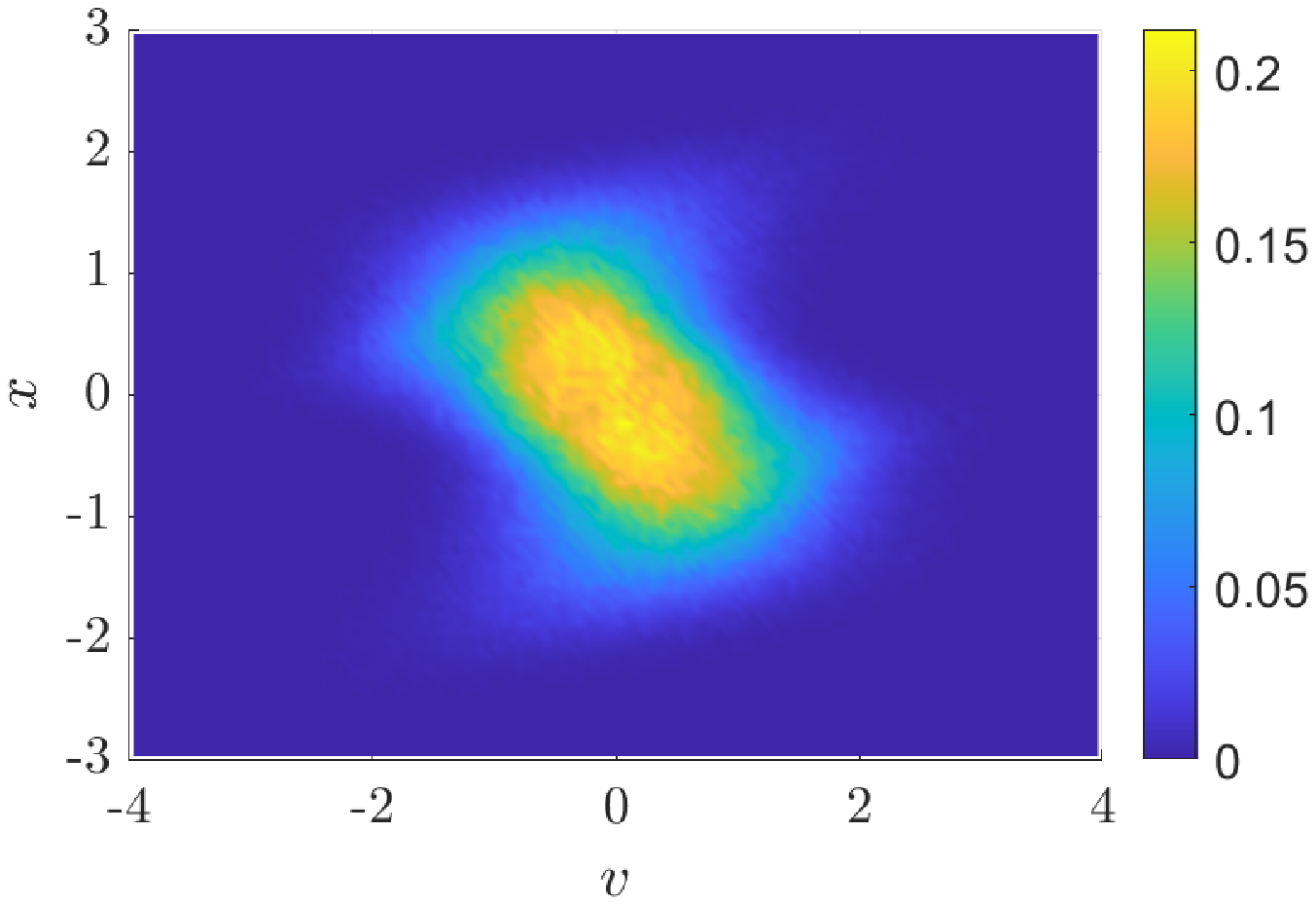}}
\subcaptionbox{Particle solution, $t = 1$}{\includegraphics[scale=0.35]{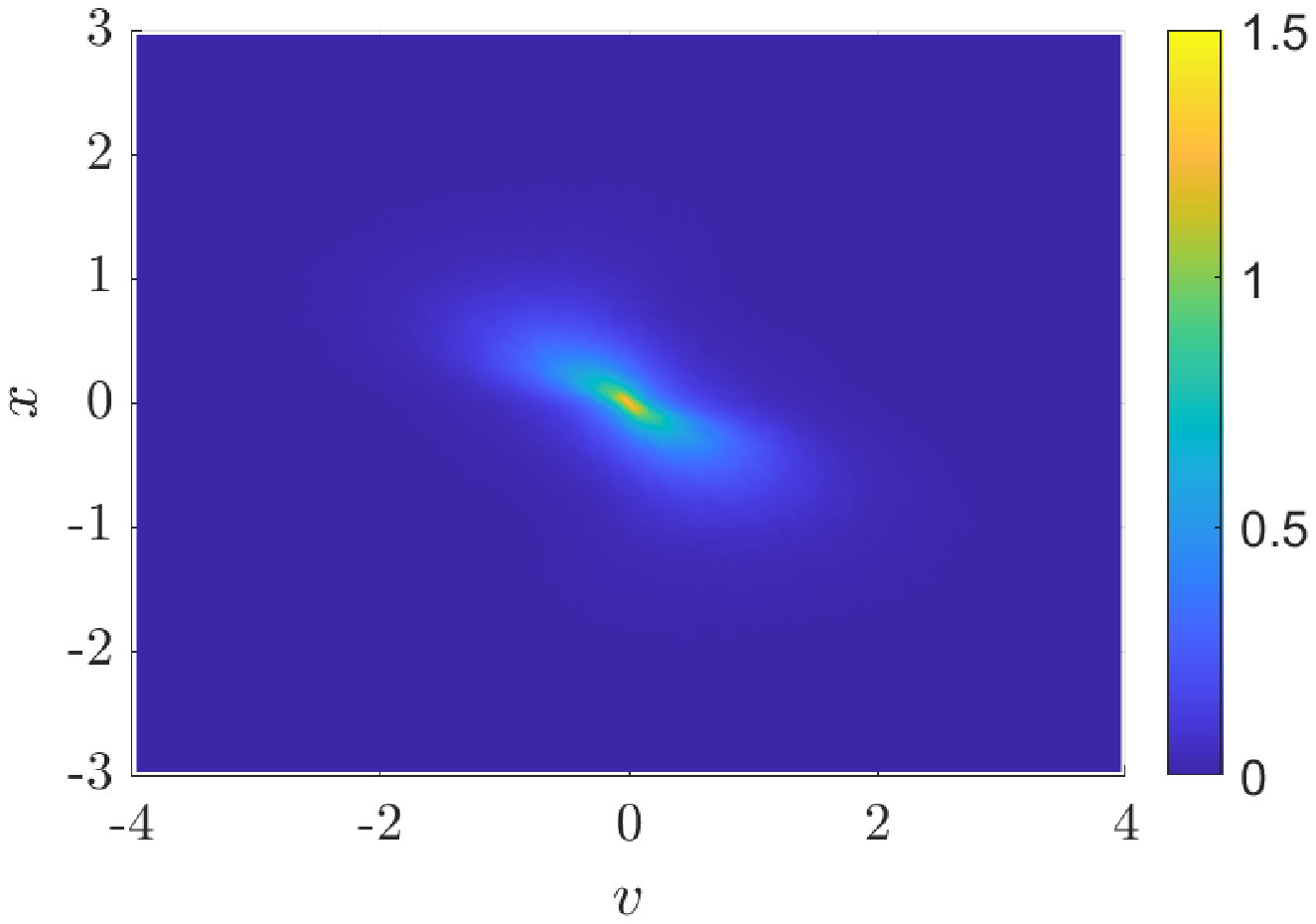}}
\subcaptionbox{Particle solution, $t = 3$}{\includegraphics[scale=0.35]{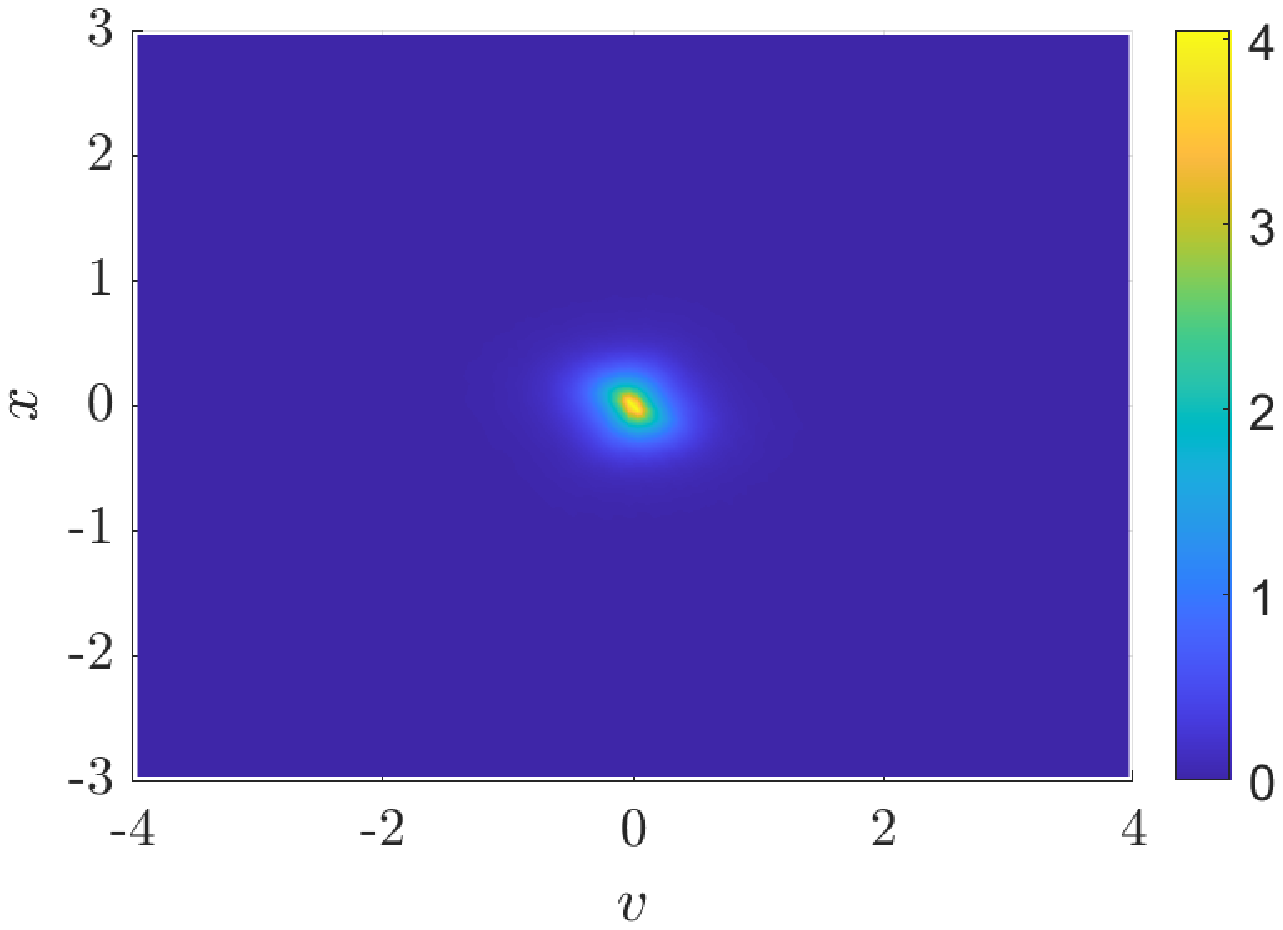}}\\
\subcaptionbox{Mean-field solution, $t = 0.5$}{\includegraphics[scale= 0.35]{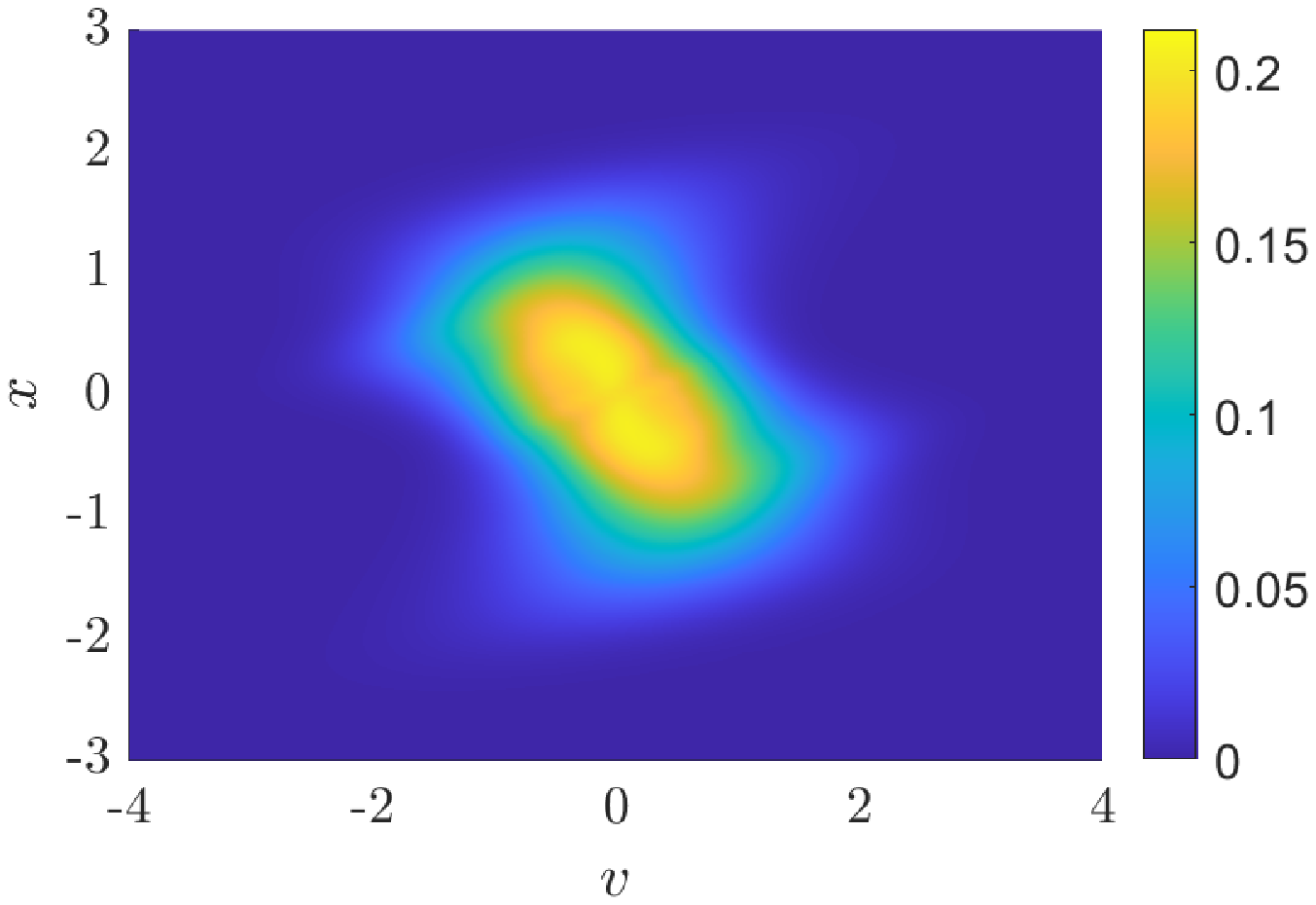}} 
\subcaptionbox{Mean-field solution, $t = 1$}{\includegraphics[scale= 0.35]{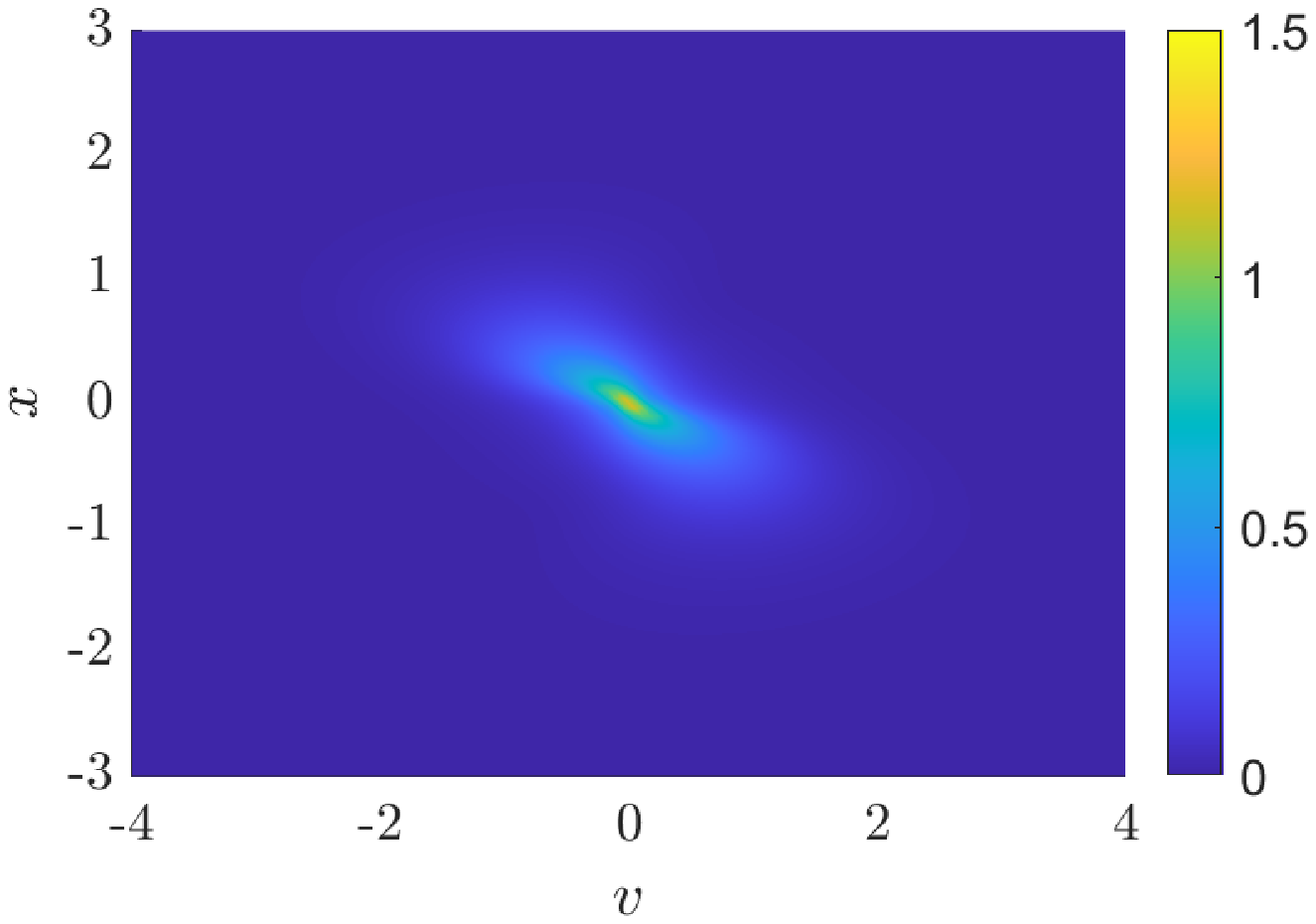}} 
\subcaptionbox{Mean-field solution, $t = 3$}{\includegraphics[scale= 0.35]{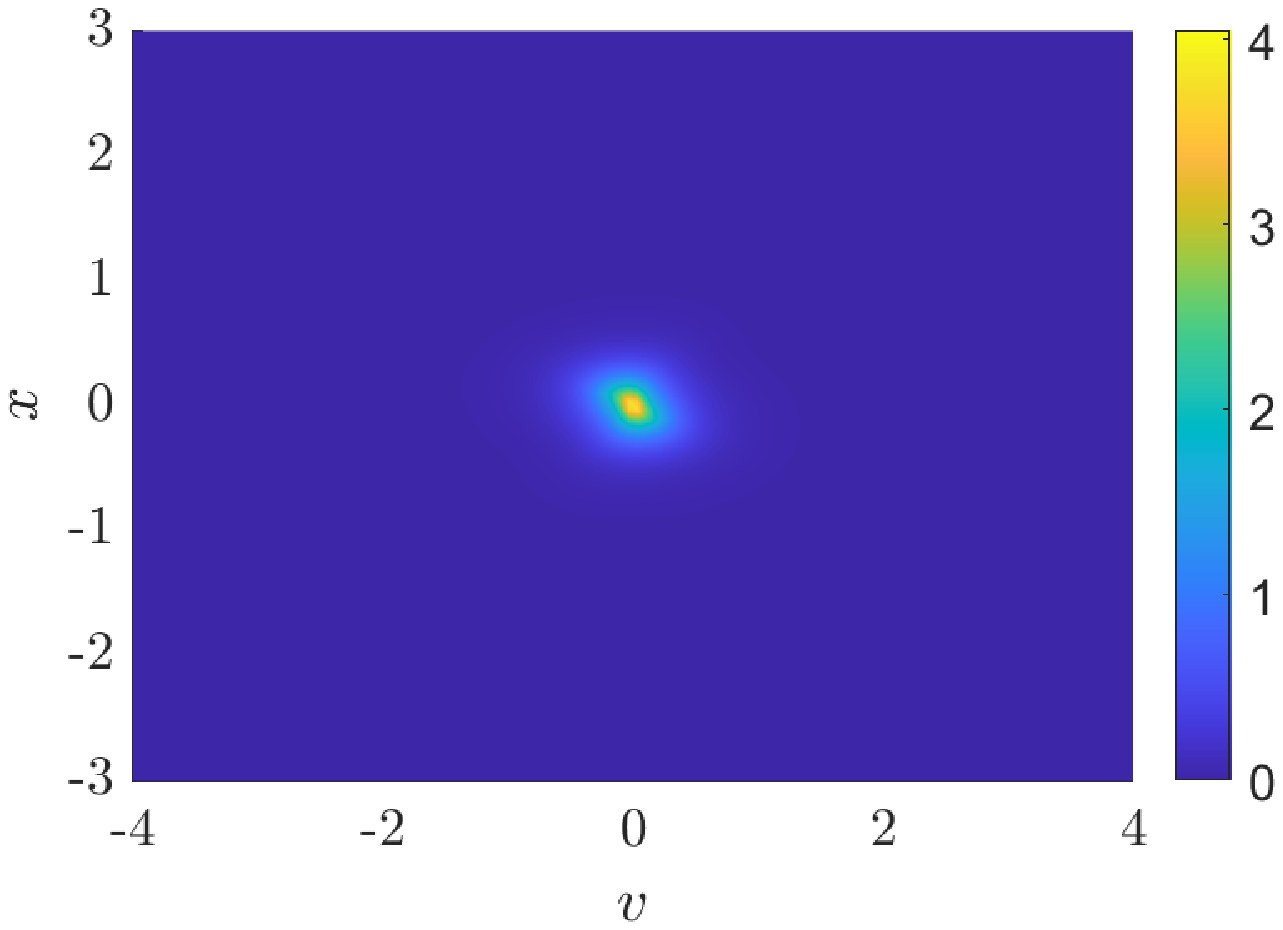}} 
\caption{Case \#3 (general case). Optimization of the one-dimensional Ackley function with minimum in $x=0$. First row: solution of the SD-PSO system \eqref{eq:psocir}. Second row: solution of the MF-PSO limit \eqref{PDE}.} 
\label{Fig11}
\end{minipage}
\vspace{5pt}
\\
\begin{minipage}{\linewidth}
\centering
\subcaptionbox{$\rho(x,t)$, $t = 0.5$}{\includegraphics[scale=0.35]{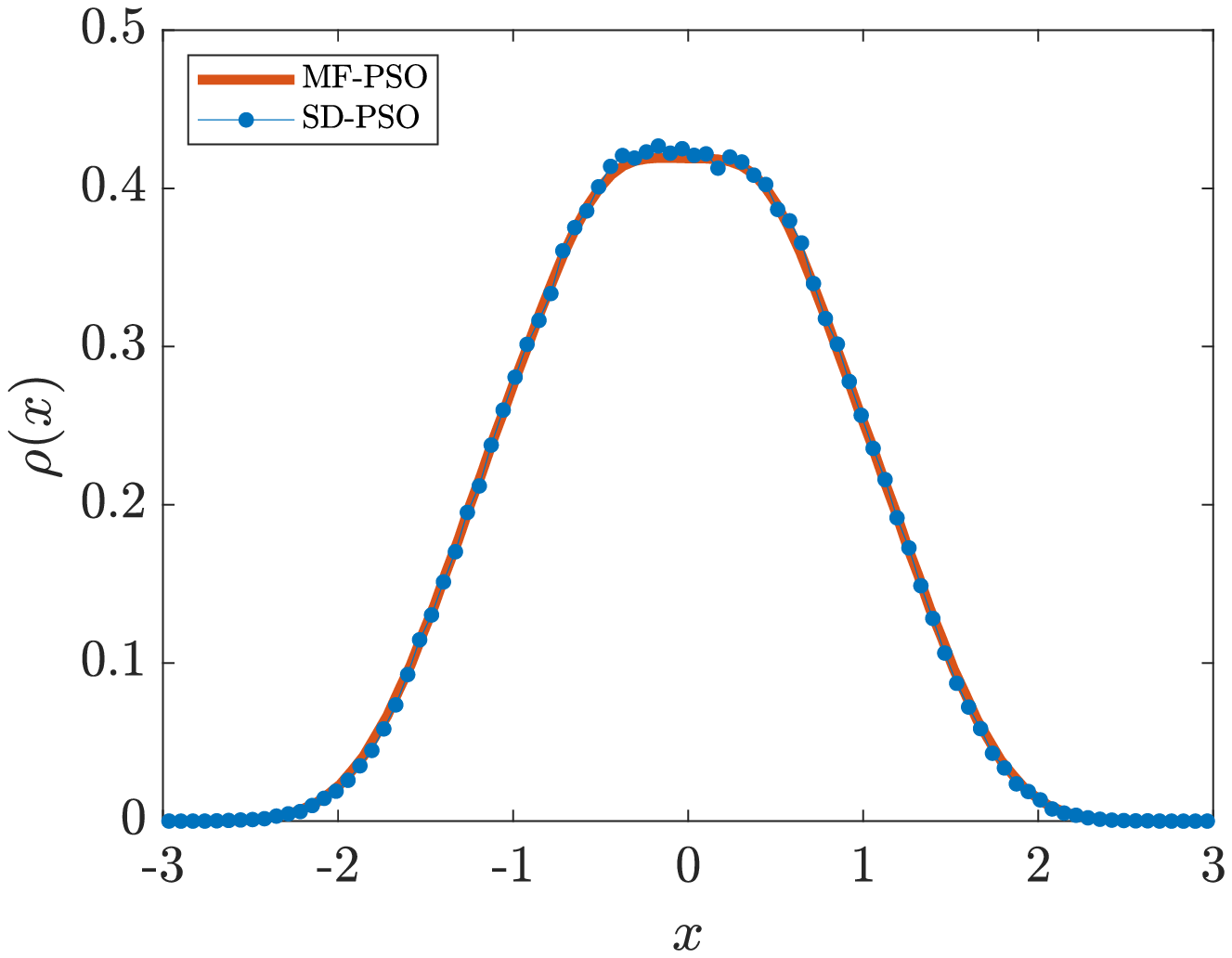}}\
\subcaptionbox{$\rho(x,t)$, $t = 1$}{\includegraphics[scale= 0.35]{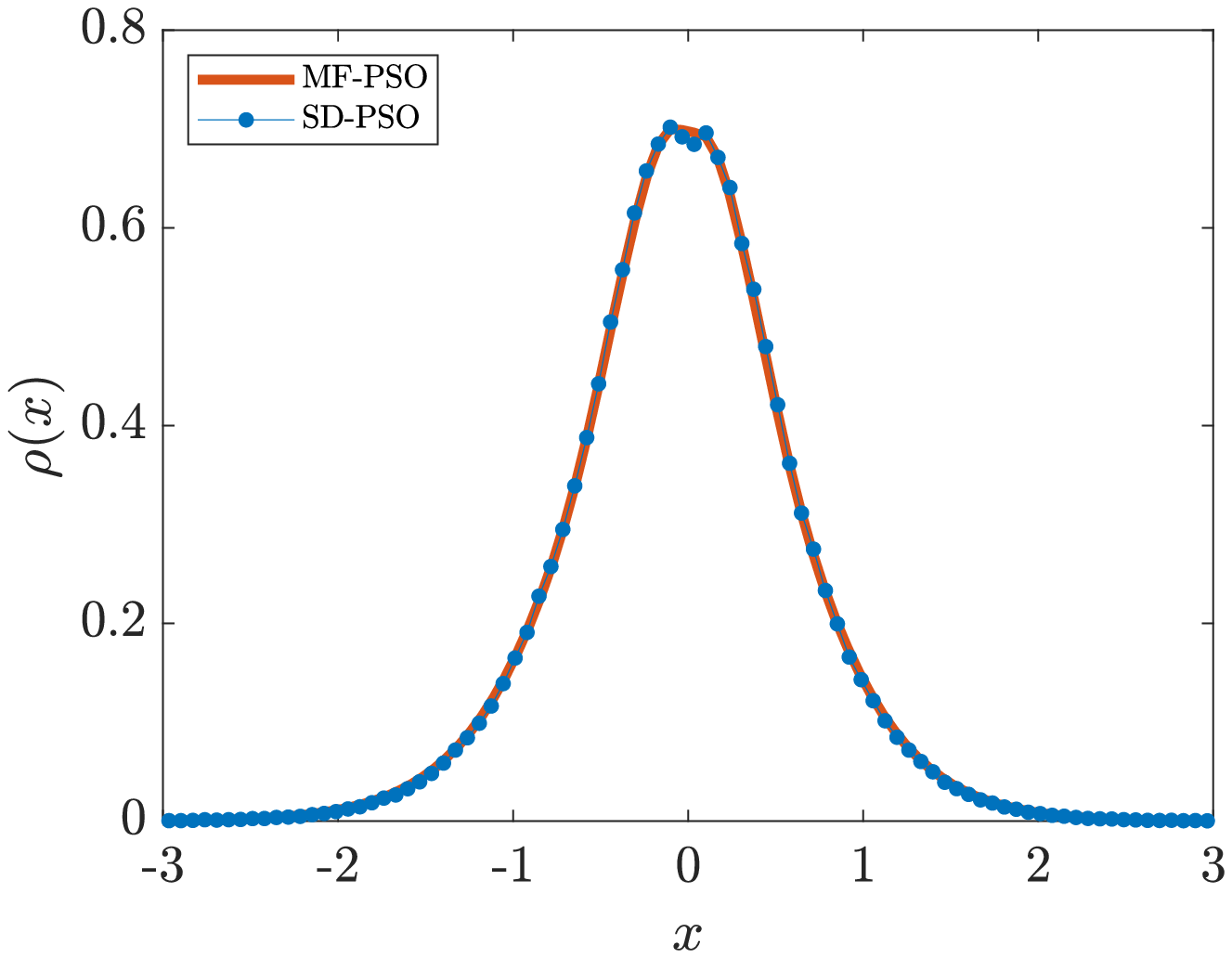}} \
\subcaptionbox{$\rho(x,t)$, $t = 3$}{\includegraphics[scale=0.35]{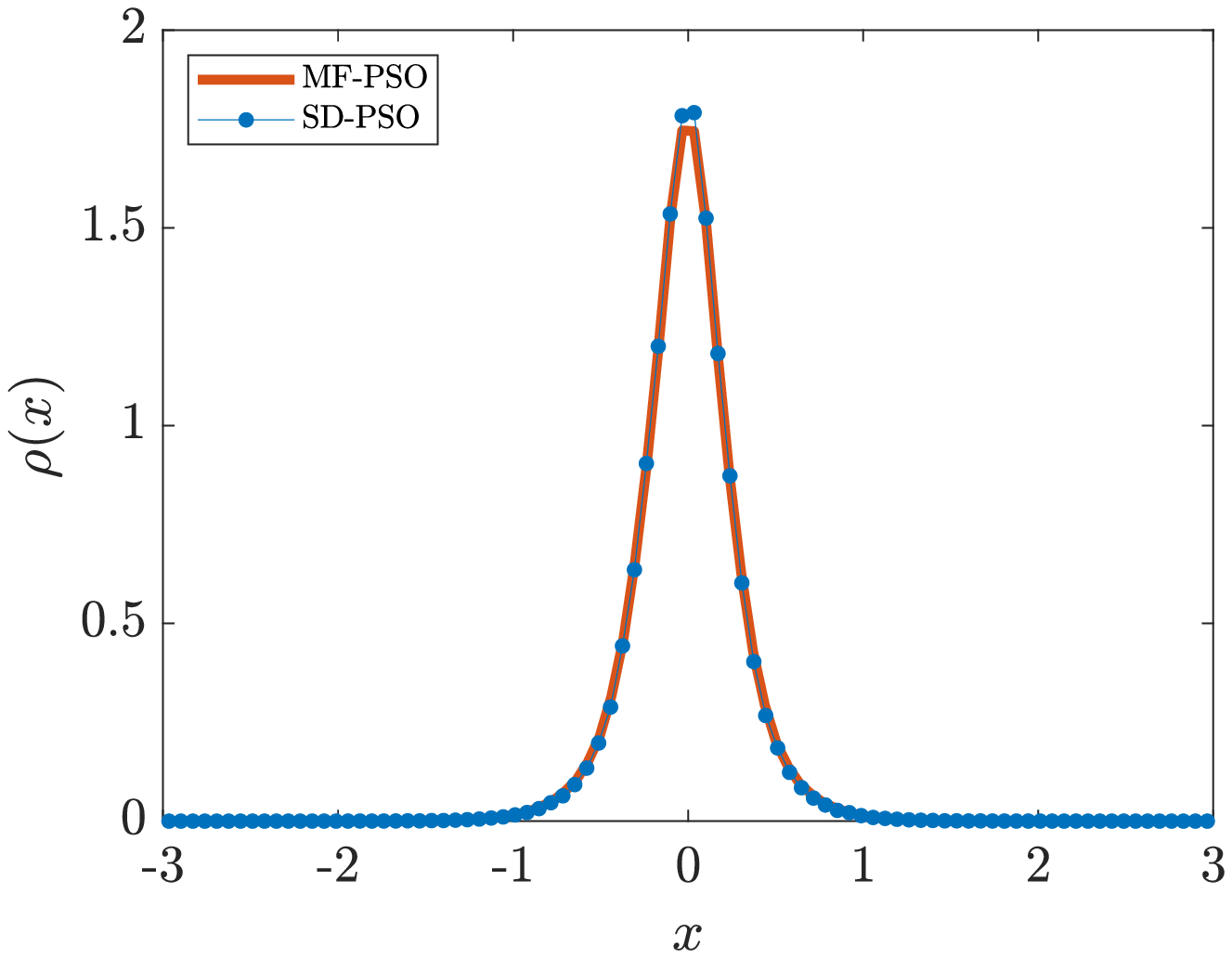}} 
\caption{Case \#3 (general case). Evolution of the density $\rho(x,t)$ of the SD-PSO system \eqref{eq:psocir} and the MF-PSO limit \eqref{PDE} for the one-dimensional Ackley function with minimum in $x=0$.}
\label{Fig12}
\end{minipage}
\end{figure}

\begin{figure}[H] 
\begin{minipage}{\linewidth}
\centering
\subcaptionbox{Particle solution, $t = 0.5$}{\includegraphics[scale=0.35]{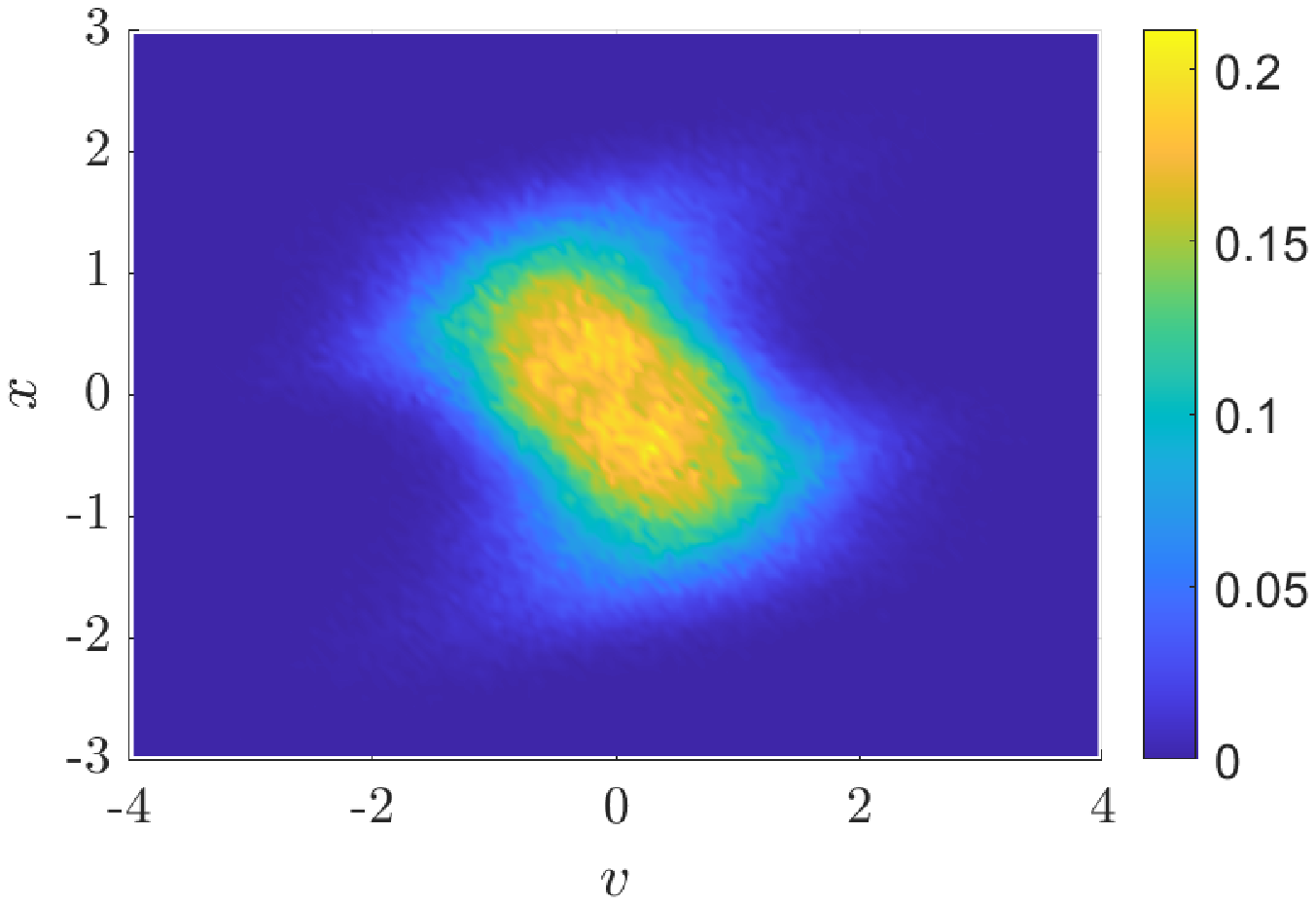}}
\subcaptionbox{Particle solution, $t = 1$}{\includegraphics[scale=0.35]{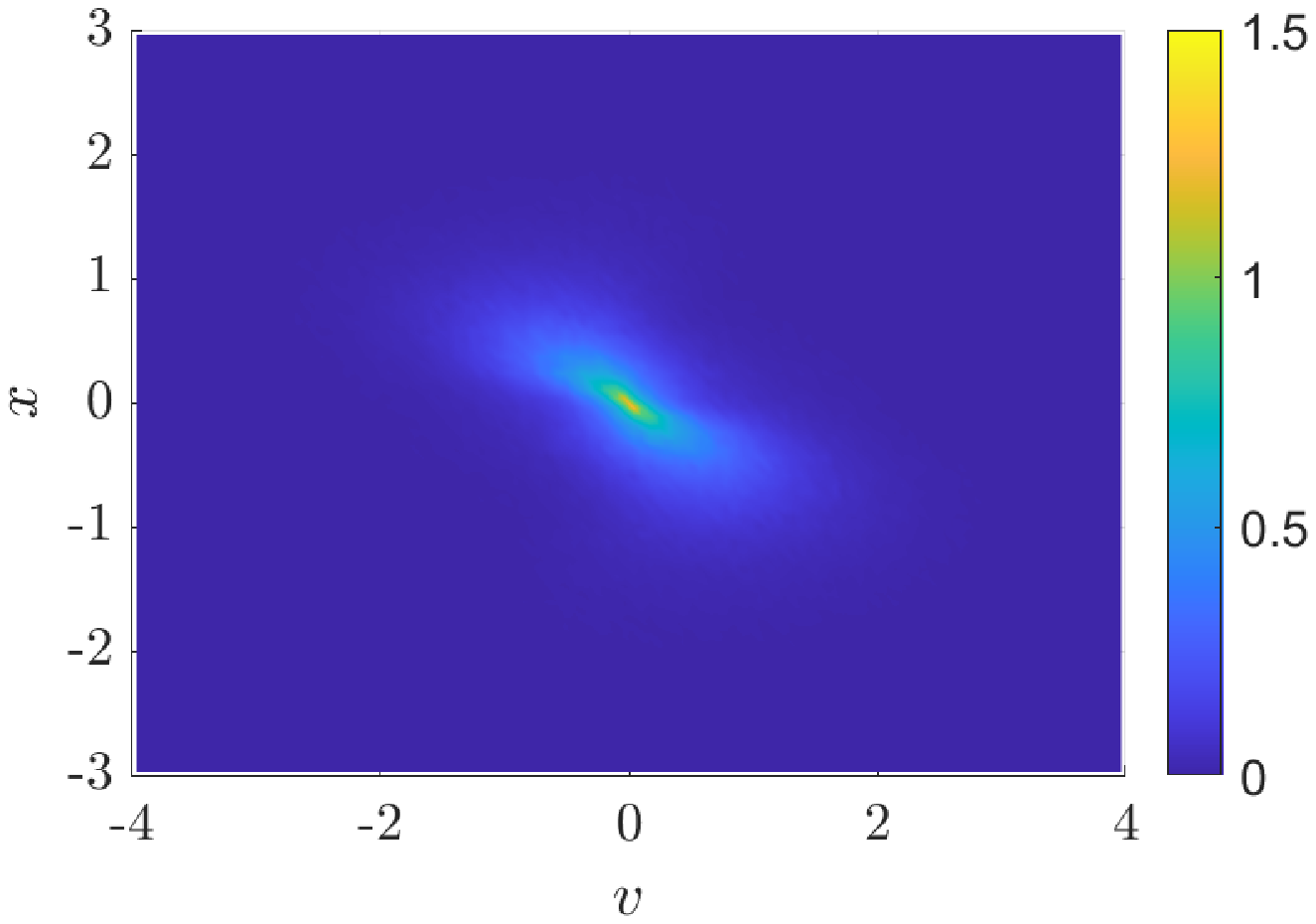}}
\subcaptionbox{Particle solution, $t = 3$}{\includegraphics[scale=0.35]{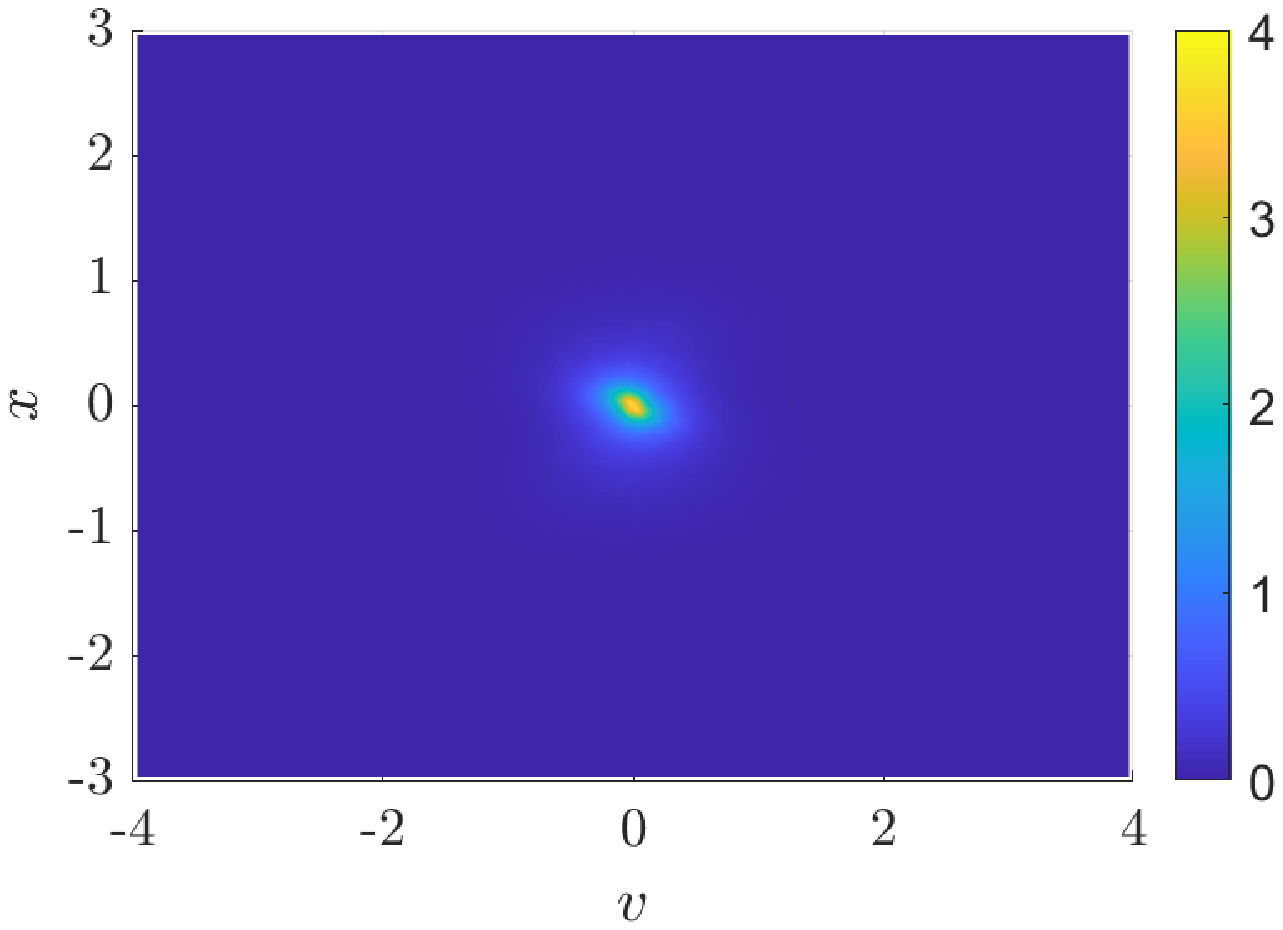}}\\
\subcaptionbox{Mean-field solution, $t = 0.5$}{\includegraphics[scale=0.35]{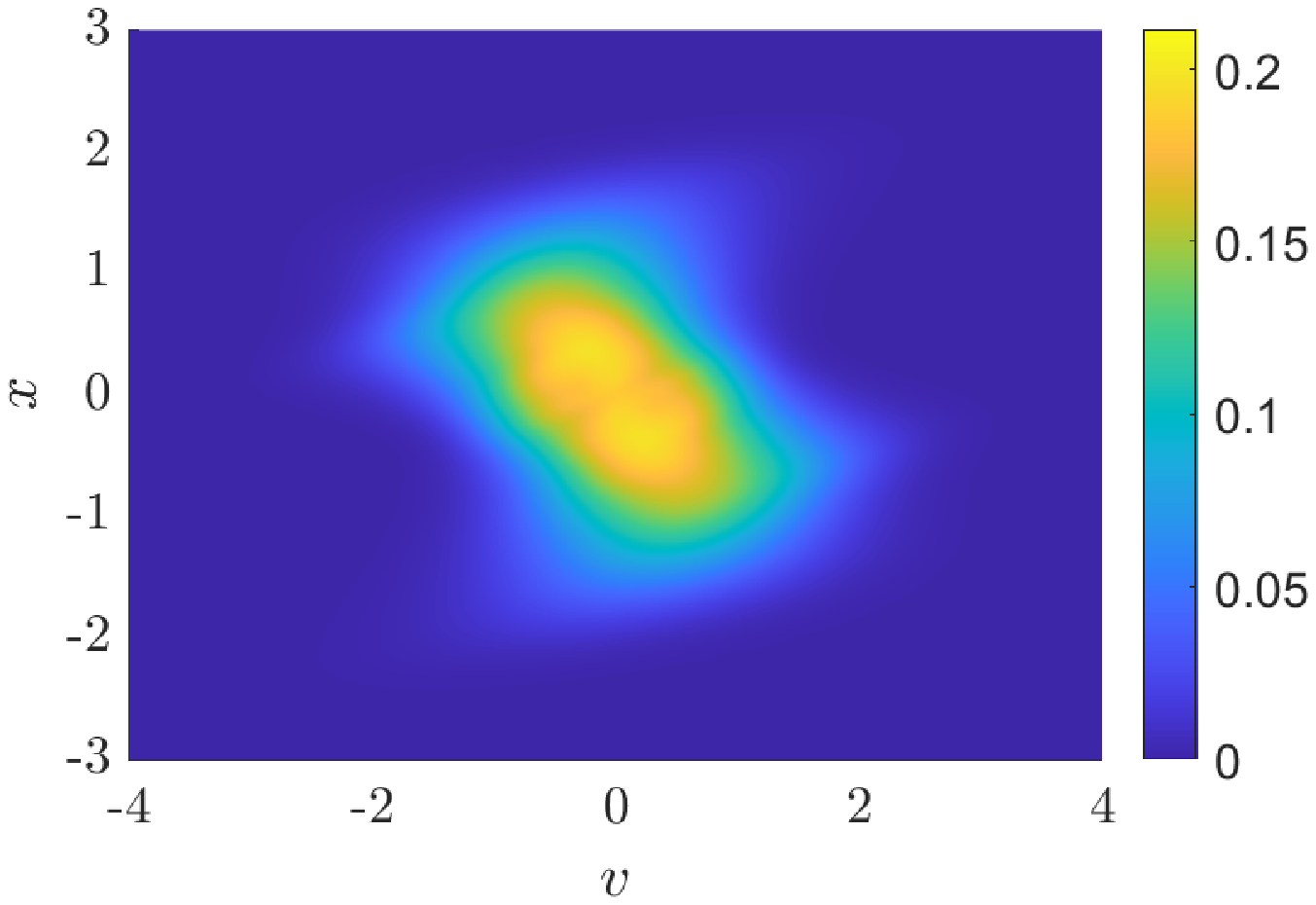}}
\subcaptionbox{Mean-field solution, $t = 1$}{\includegraphics[scale=0.35]{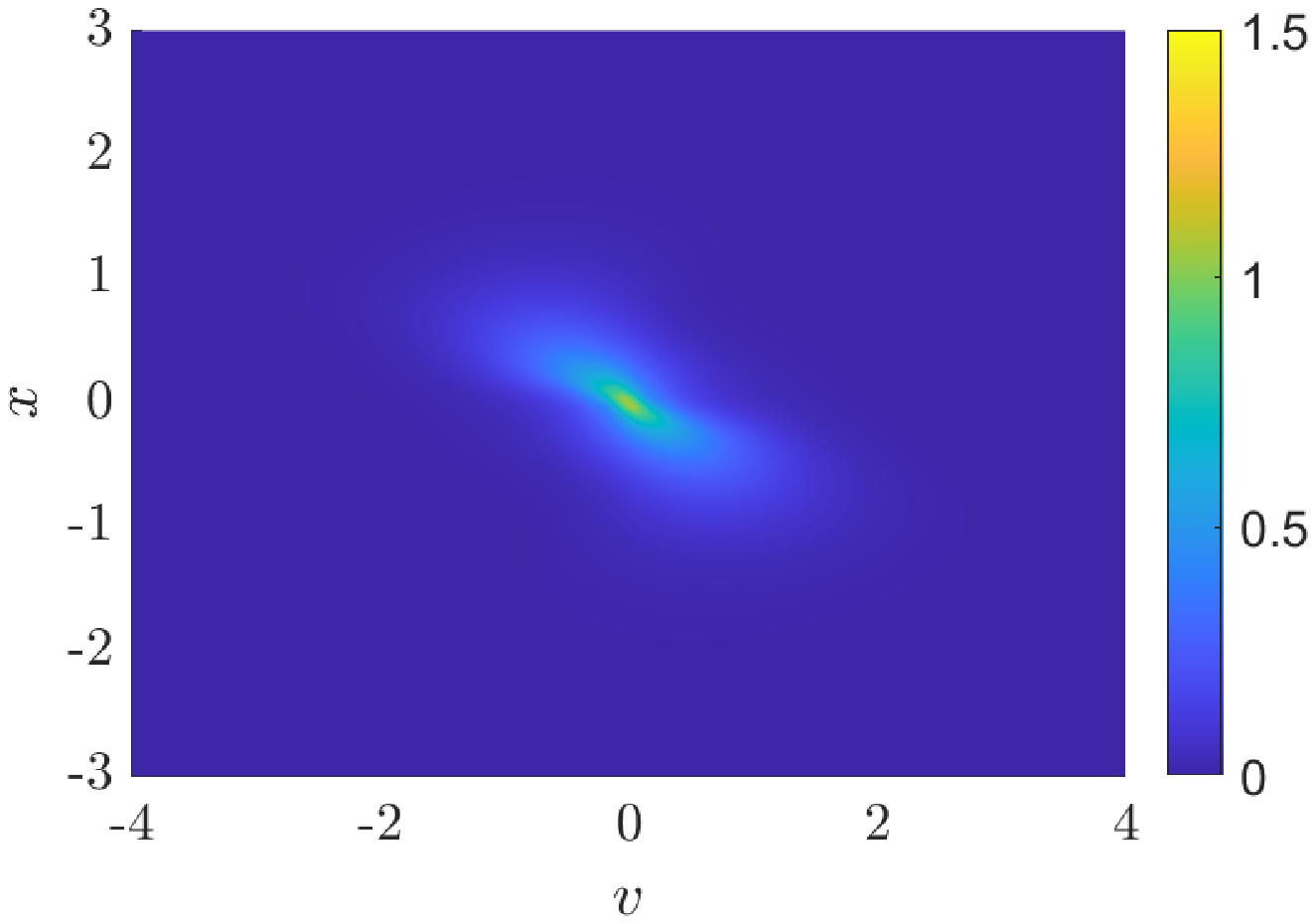}} 
\subcaptionbox{Mean-field solution, $t = 3$}{\includegraphics[scale=0.35]{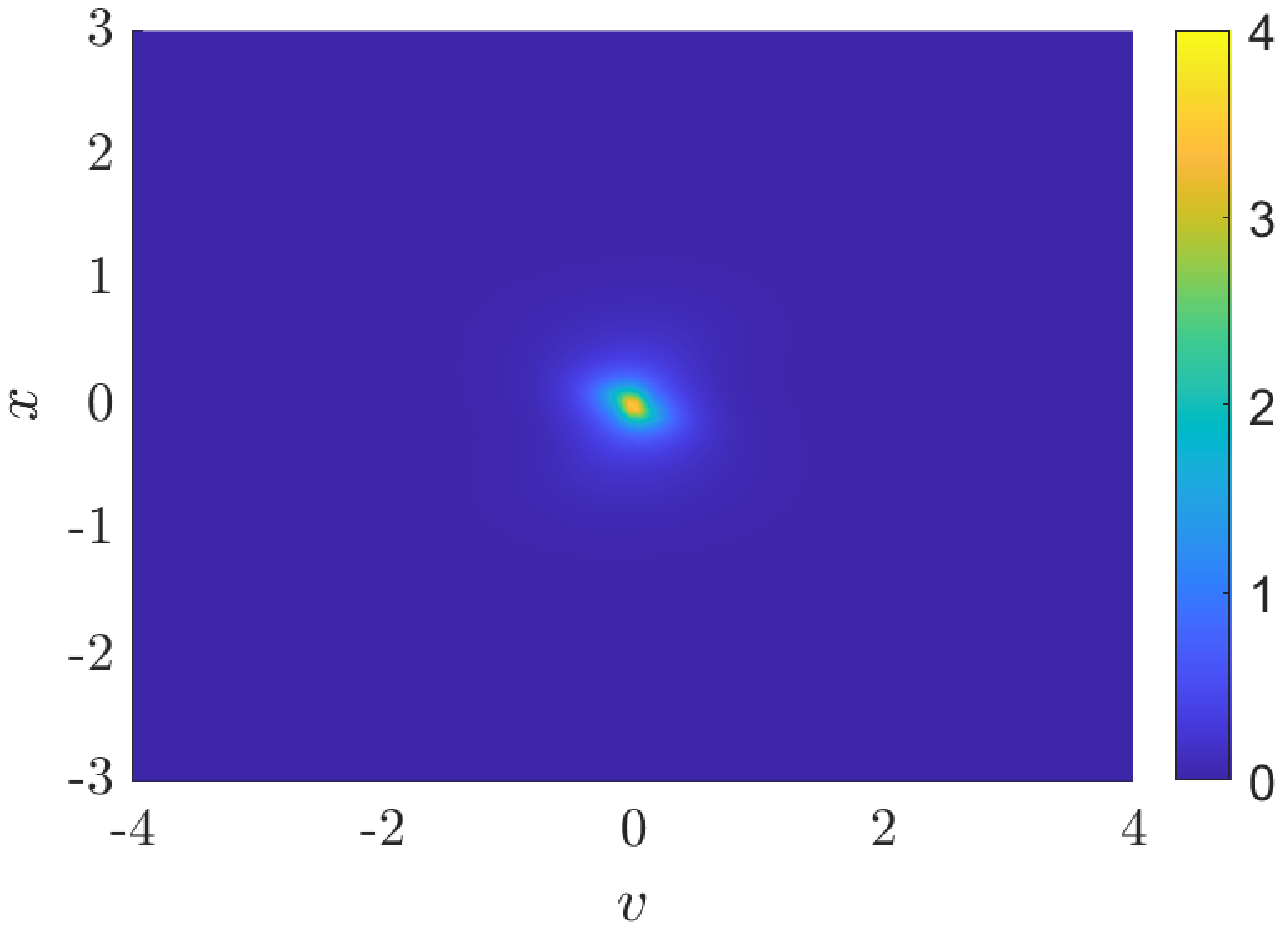}} 
\caption{Case \#3 (general case). Optimization of the one-dimensional Rastrigin function with minimum in $x=0$. First row: solution of the SD-PSO system \eqref{eq:psocir}. Second row: solution of the MF-PSO limit \eqref{PDE}.} 
\label{Fig13}
\end{minipage}
\vspace{5pt}
\\
\begin{minipage}{\linewidth}
\centering
\subcaptionbox{$\rho(x,t)$, $t = 0.5$}{\includegraphics[scale=0.35]{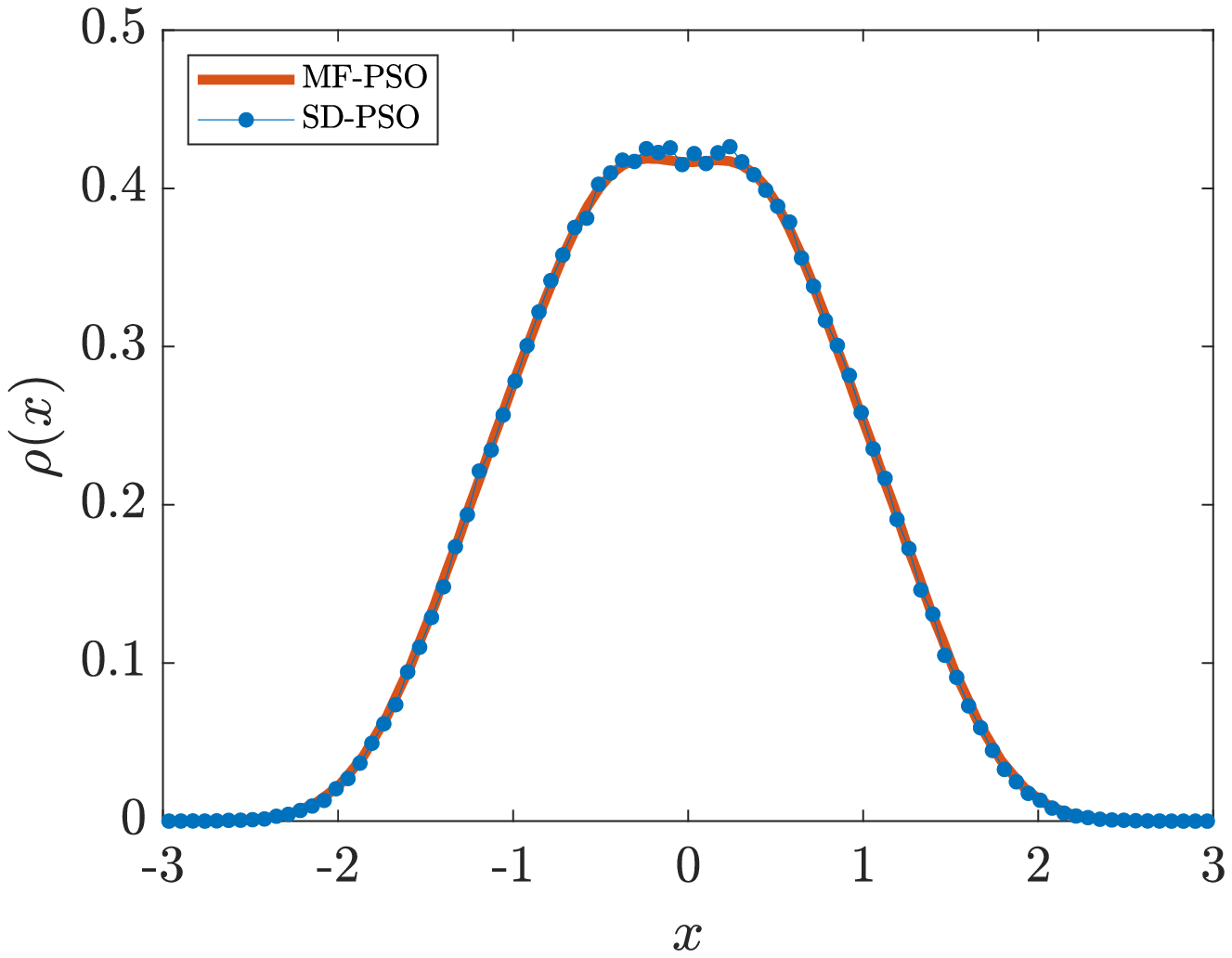}}\
\subcaptionbox{$\rho(x,t)$, $t = 1$}{\includegraphics[scale=0.35]{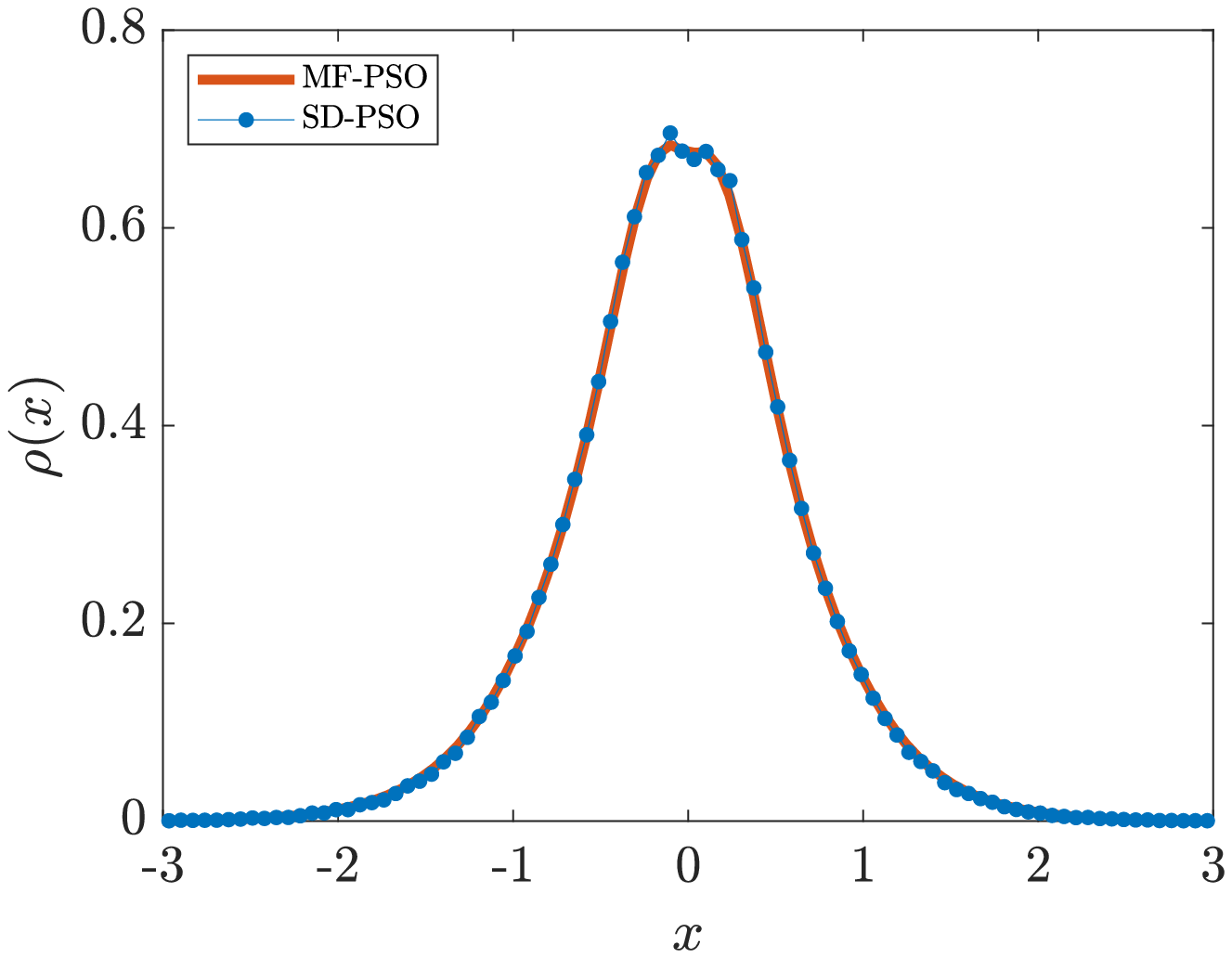}} \
\subcaptionbox{$\rho(x,t)$, $t = 3$}{\includegraphics[scale=0.35]{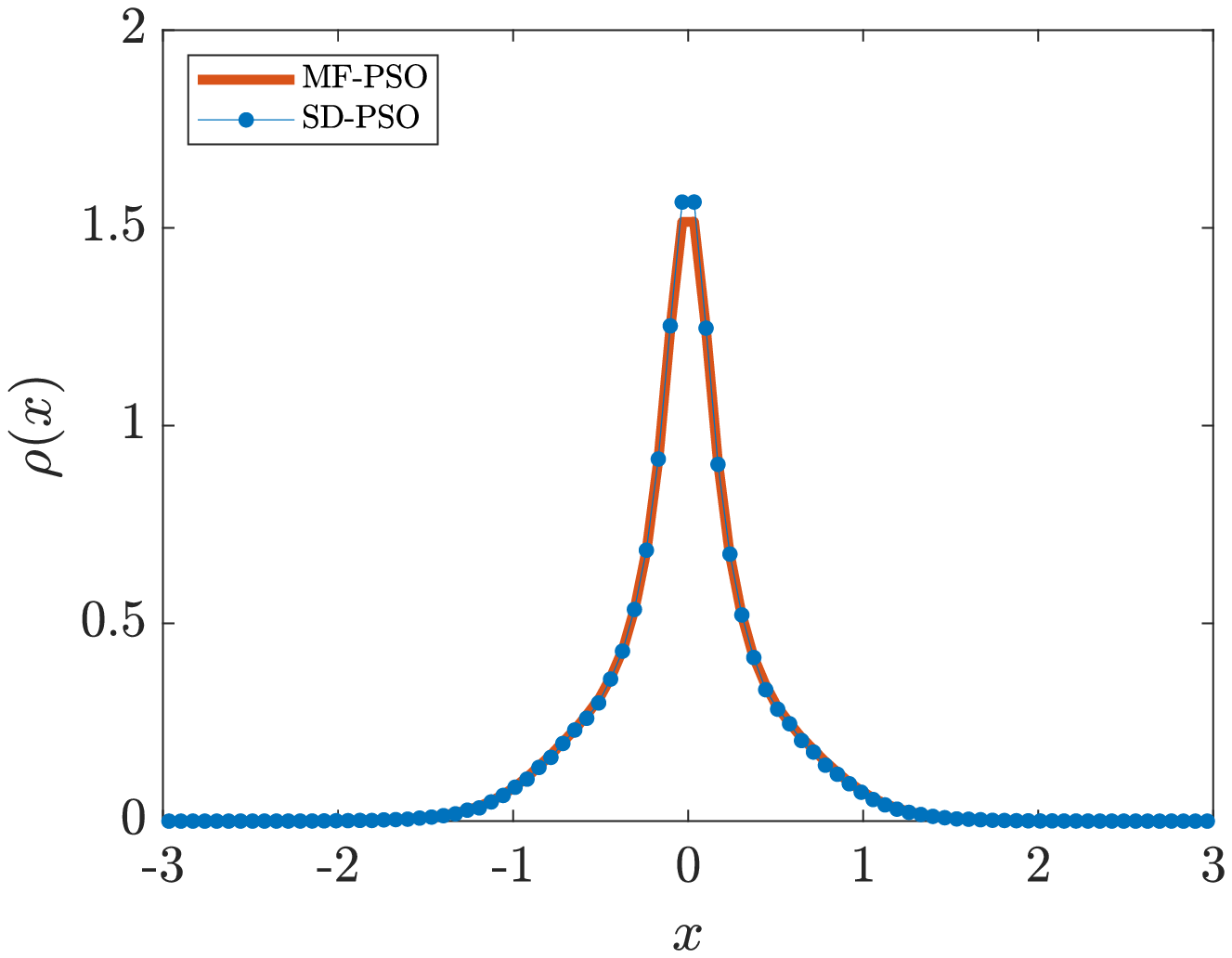}}
\caption{Case \#3 (general case). Evolution of the density $\rho(x,t)$ of the SD-PSO system \eqref{eq:psocir} and the MF-PSO limit \eqref{PDE} for the one-dimensional Rastrigin function with minimum in $x=0$.}
\label{Fig14}
\end{minipage}
\end{figure}

In Figures \ref{Fig11} and \ref{Fig13} we report the contour plots of the solutions obtained with the discretized stochastic particle system \eqref{eq:psoDiscr} and the deterministic solver of the mean field equation \eqref{PDE}.
One can immediately observe that the local minima effect disappears and the systems converge consistently towards the global minima for both the Ackley and the Rastrigin functions. The good agreement between the particle and the mean-field solutions, as before, is emphasized by the density plots in Figures \ref{Fig12} and \ref{Fig14}. Note that, by comparing the results in Figure \ref{Fig12} and those in Figure \ref{Fig4} obtained by solving the same problem in absence of memory terms and presence of global best only, at the same time instants, a faster convergence towards the global minimum is observed in Figure \ref{Fig12} thanks to the inclusion of the memory effects in the dynamic.

\subsection{Numerical small inertia limit}
From the analysis in Section 4, in the limit of small inertia the classical CBO model \eqref{eq:CBOp} is obtained as hydrodynamic approximation of the mean-field PSO system \eqref{PDEi}. Therefore, starting from the discretization of the stochastic particle model without memory effect \eqref{eq:psoiDiscr}, we decrease the inertial weight $m \to 0$ ($\gamma \to 1$) and compare the particle solution with a direct discretization of the limiting mean-field CBO system \eqref{eq:CBOp}.

\begin{figure}[H]
\begin{minipage}{\linewidth}
\centering
\subcaptionbox{$\rho(x,t)$, $t = 0.2$}{\includegraphics[scale=0.35]{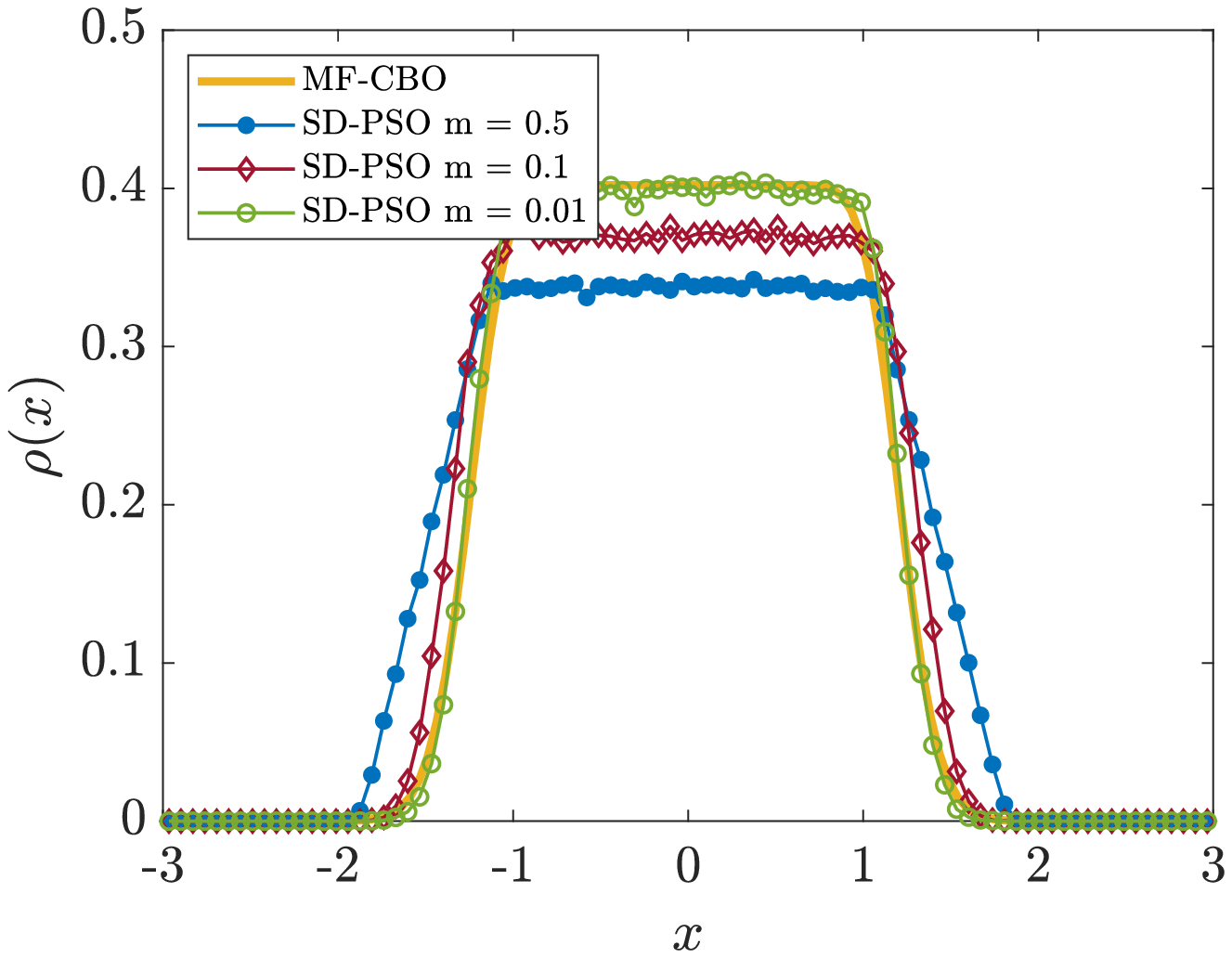}}\
\subcaptionbox{$\rho(x,t)$, $t = 1$}{\includegraphics[scale=0.35]{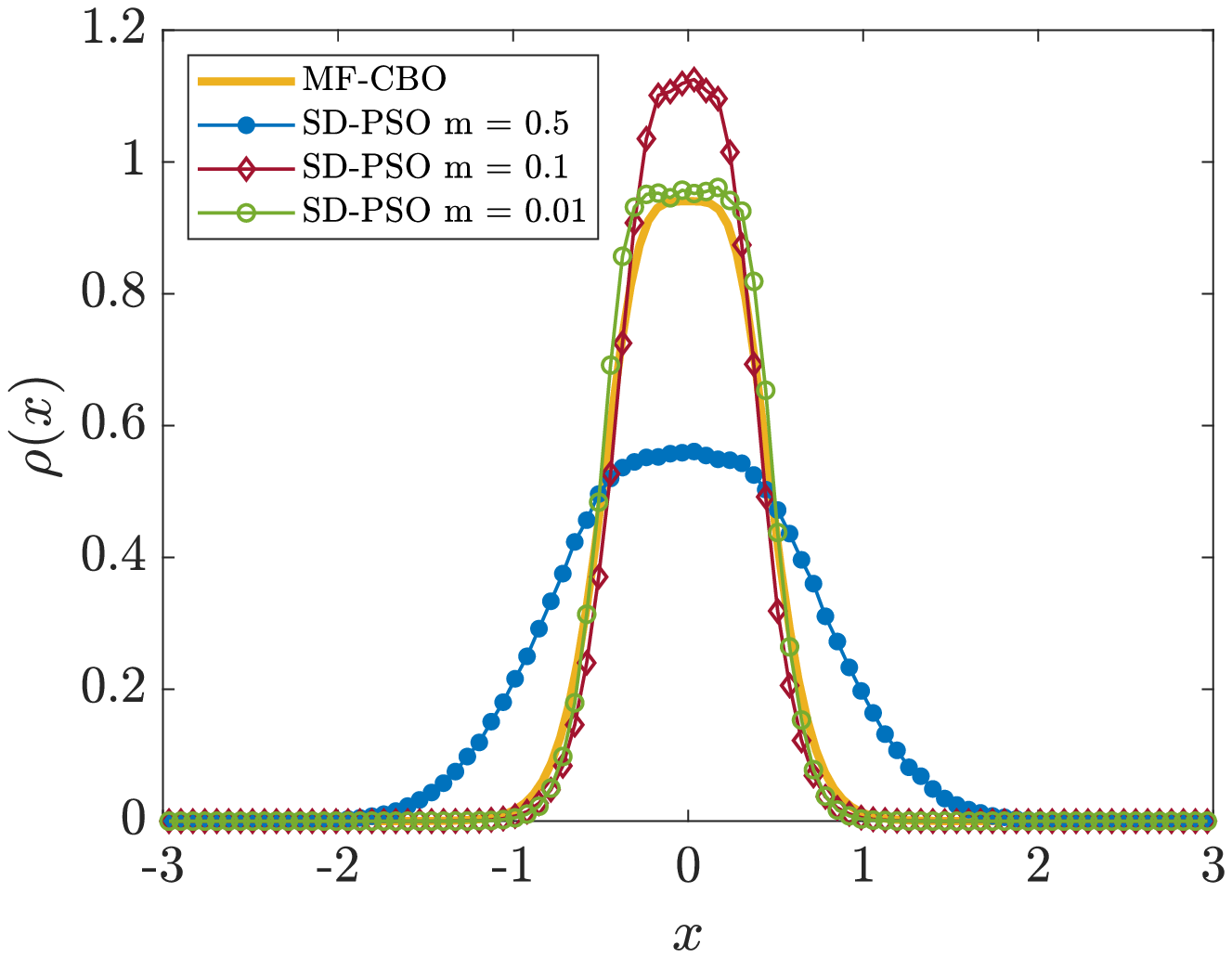}} \
\subcaptionbox{$\rho(x,t)$, $t = 2$}{\includegraphics[scale=0.35]{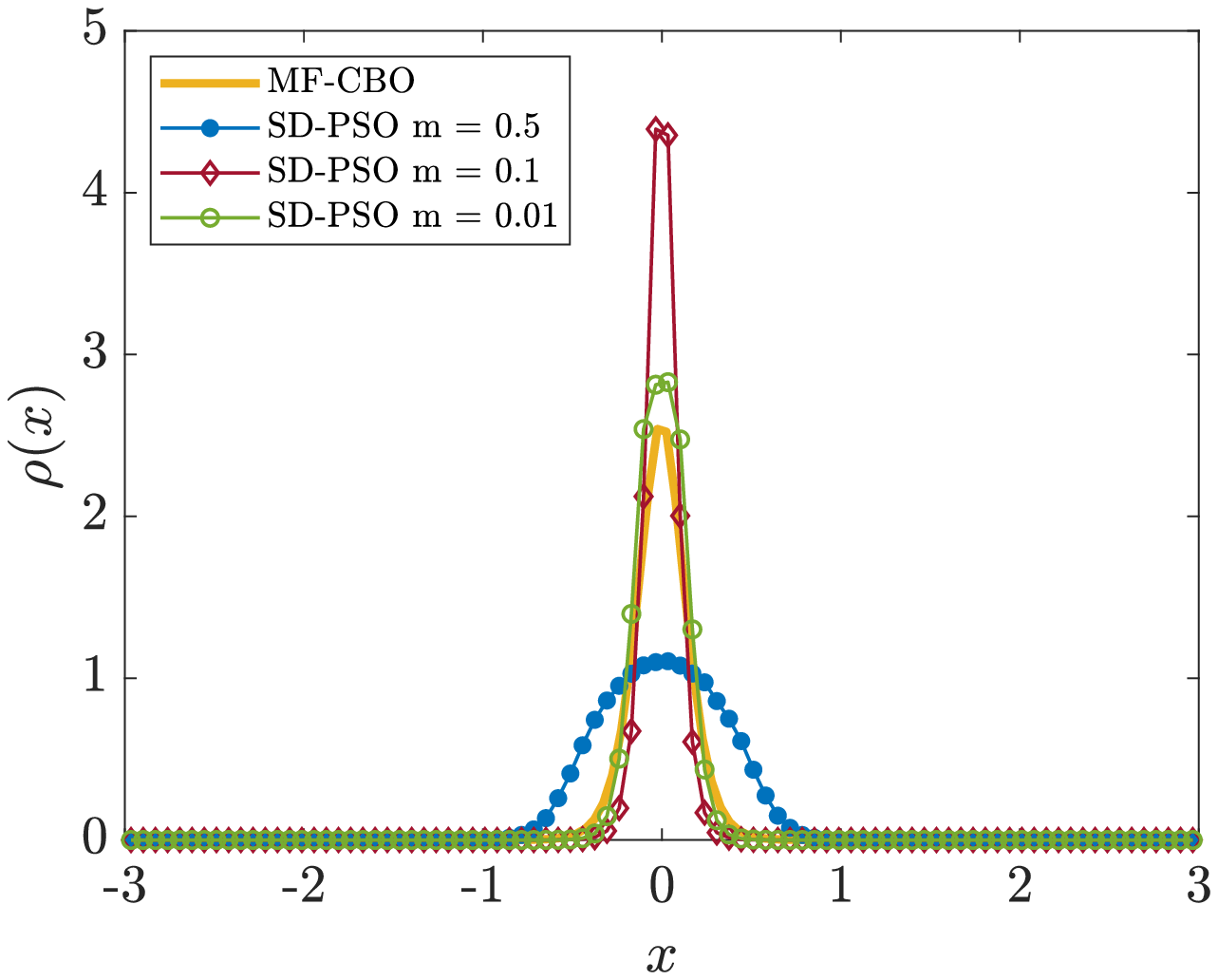}} \\
\subcaptionbox{$\rho(x,t)$, $t = 0.2$}{\includegraphics[scale=0.35]{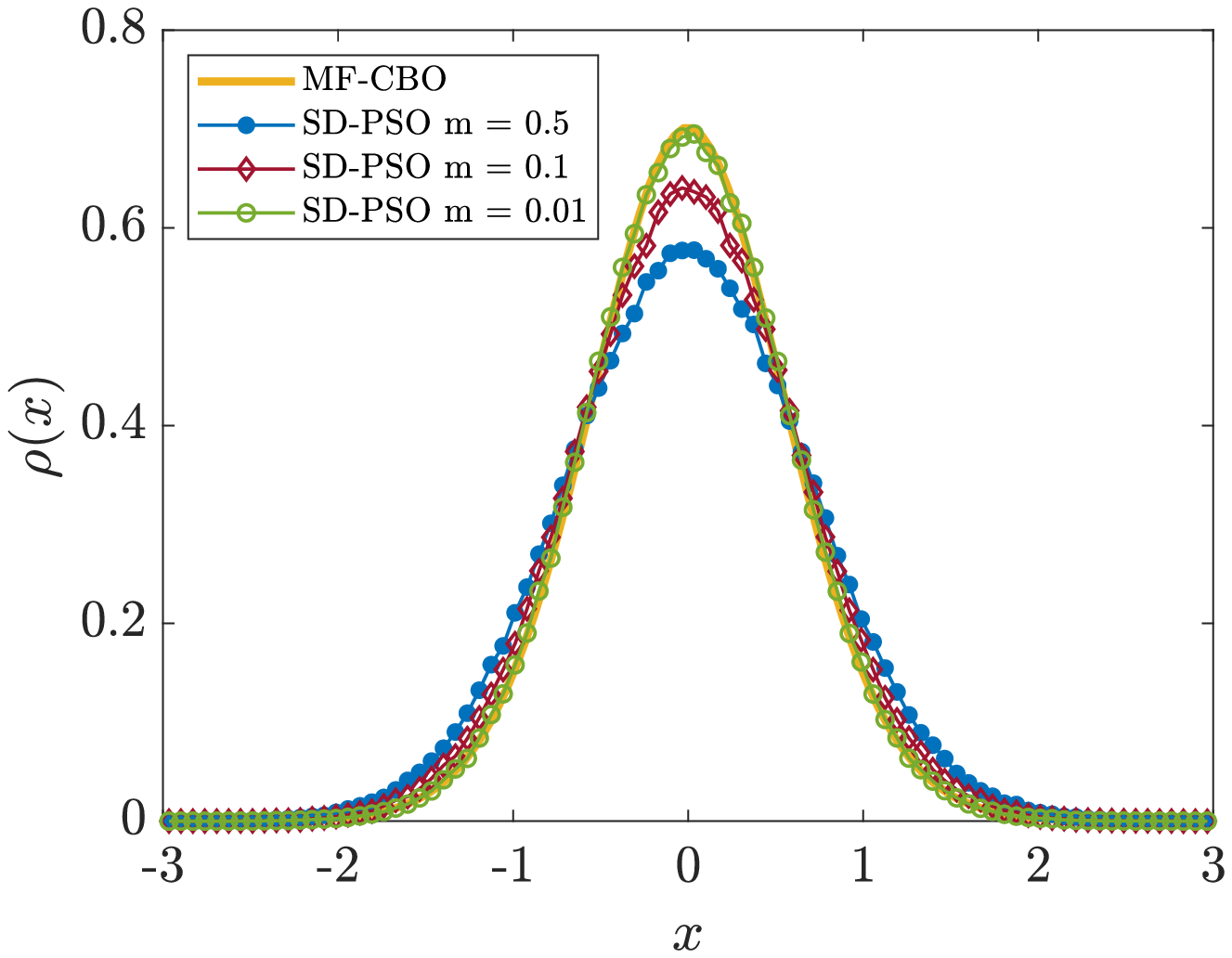}}\
\subcaptionbox{$\rho(x,t)$, $t = 1$}{\includegraphics[scale=0.35]{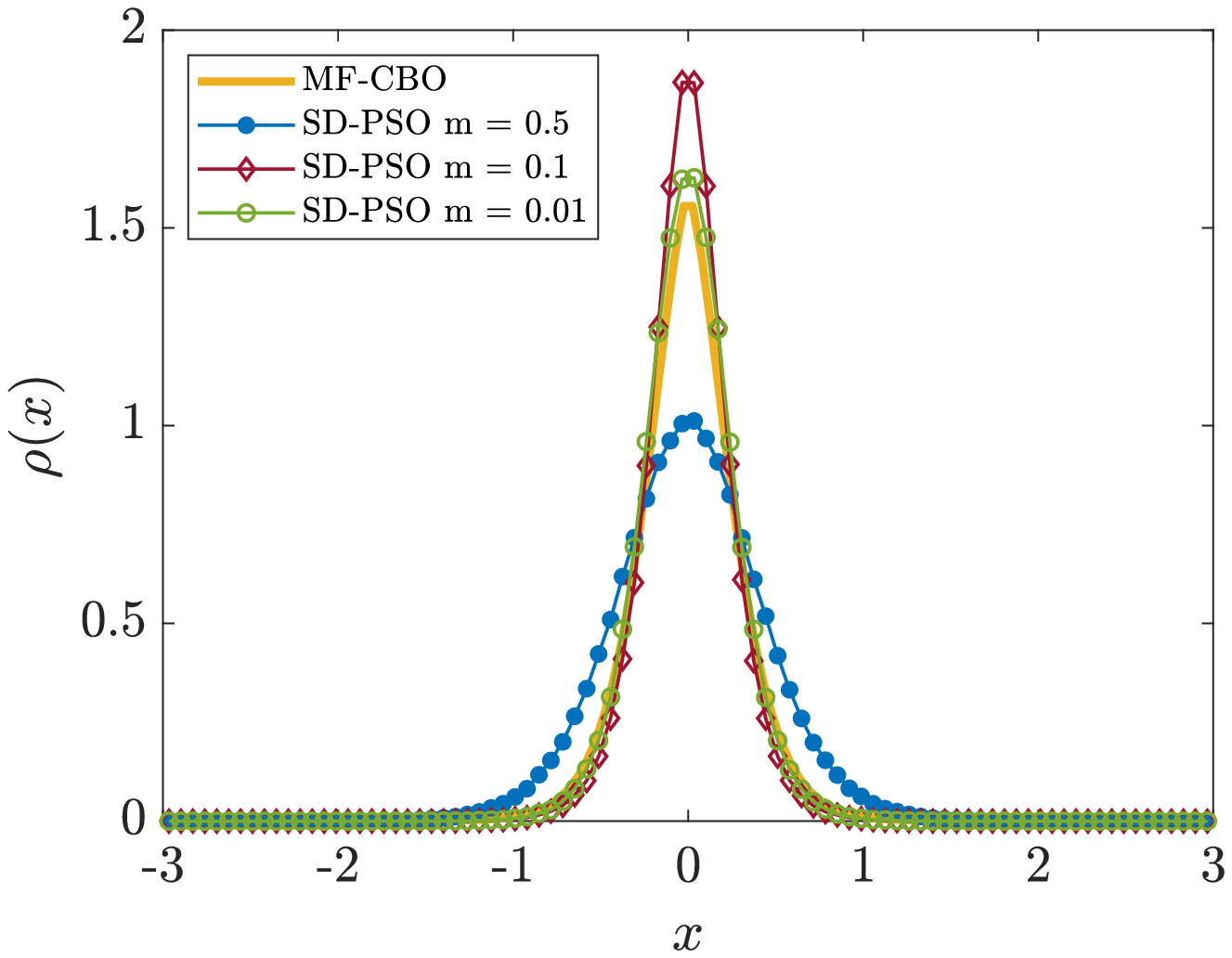}} \
\subcaptionbox{$\rho(x,t)$, $t = 2$}{\includegraphics[scale=0.35]{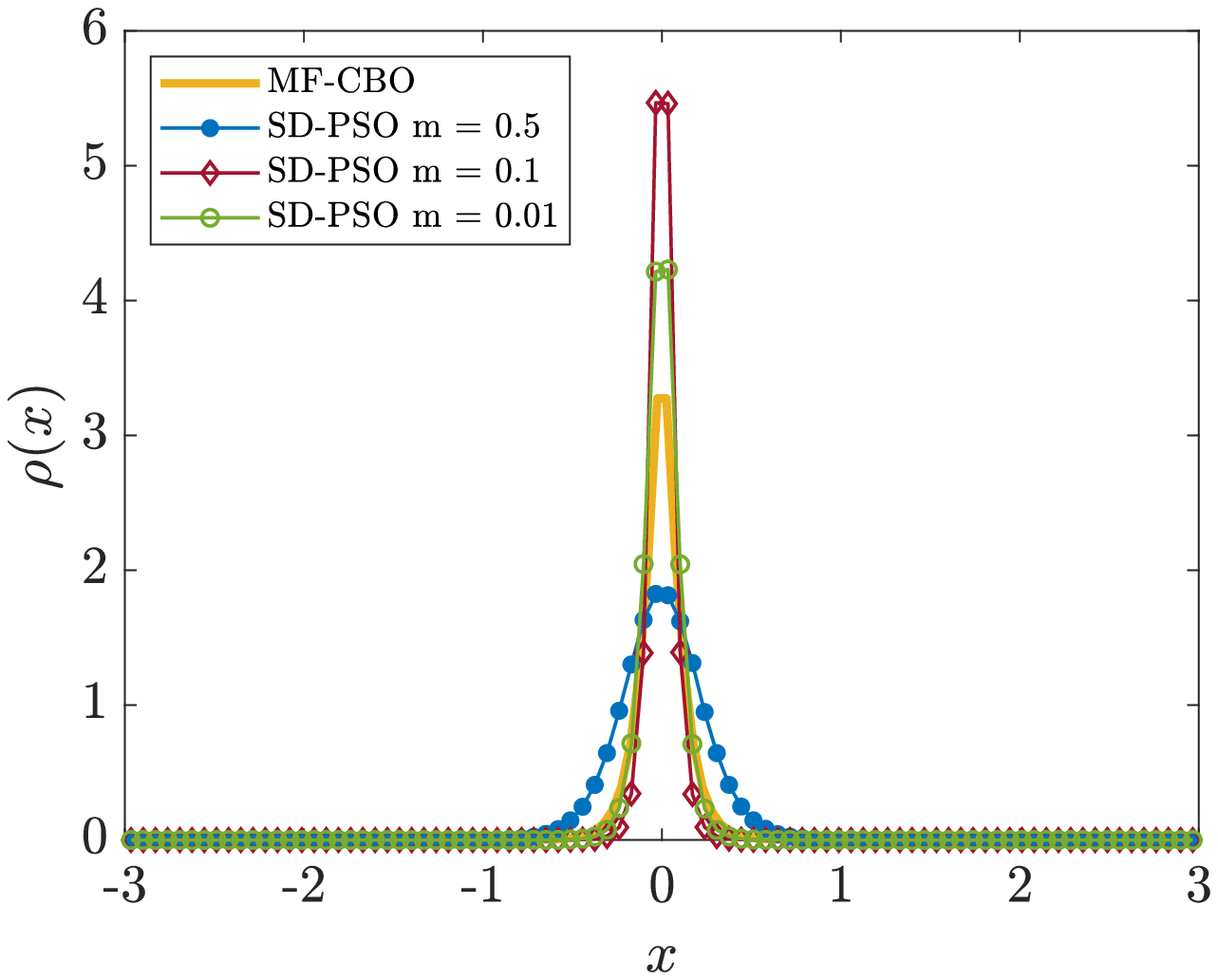}} 
\caption{Low inertia limit. Evolution of the density $\rho(x,t)$ of the SD-PSO discretization \eqref{eq:psoiDiscr2}, for decreasing inertial weight $m = 0.5, 0.1, 0.01$, and the mean-field CBO model \eqref{eq:CBOp} for the Ackley function with minimum in $x=0$. First row: uniform initial data. Second row: Gaussian initial data.}
\label{Fig15}
\end{minipage}
\end{figure}

\begin{figure}[H]
\begin{minipage}{\linewidth}
\centering
\subcaptionbox{$\rho(x,t)$, $t = 0.2$}{\includegraphics[scale=0.35]{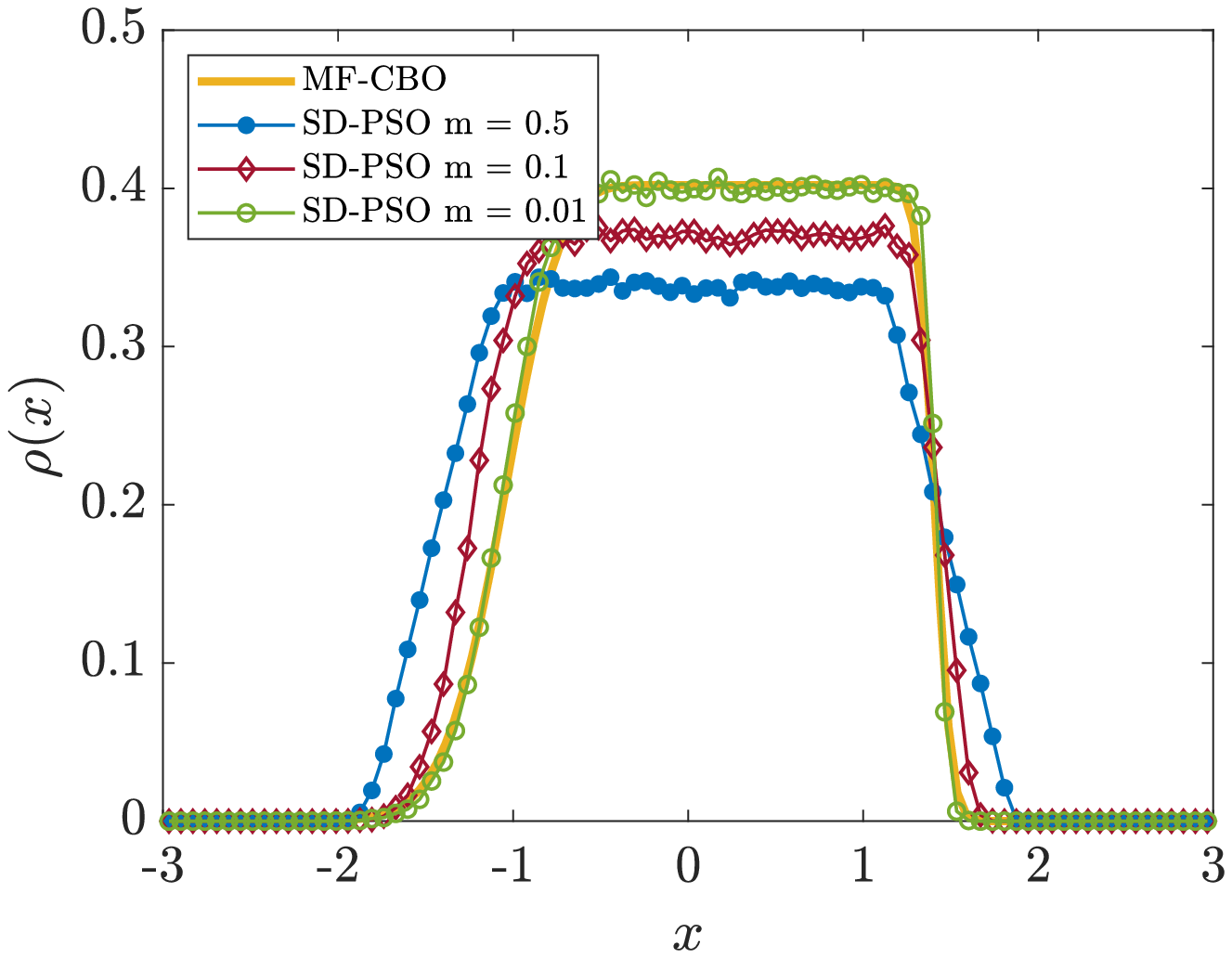}}\
\subcaptionbox{$\rho(x,t)$, $t = 1$}{\includegraphics[scale=0.35]{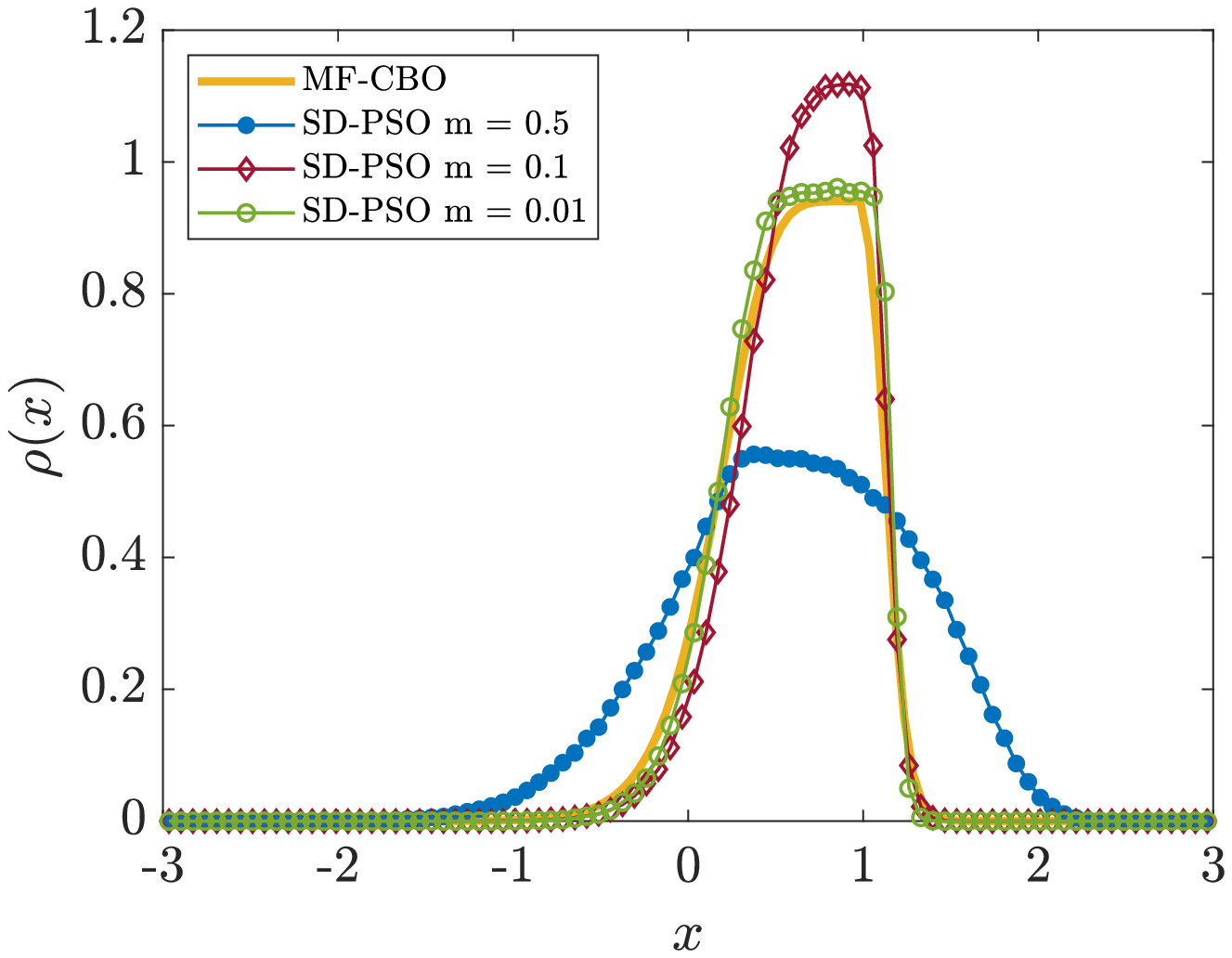}} \
\subcaptionbox{$\rho(x,t)$, $t = 2$}{\includegraphics[scale=0.35]{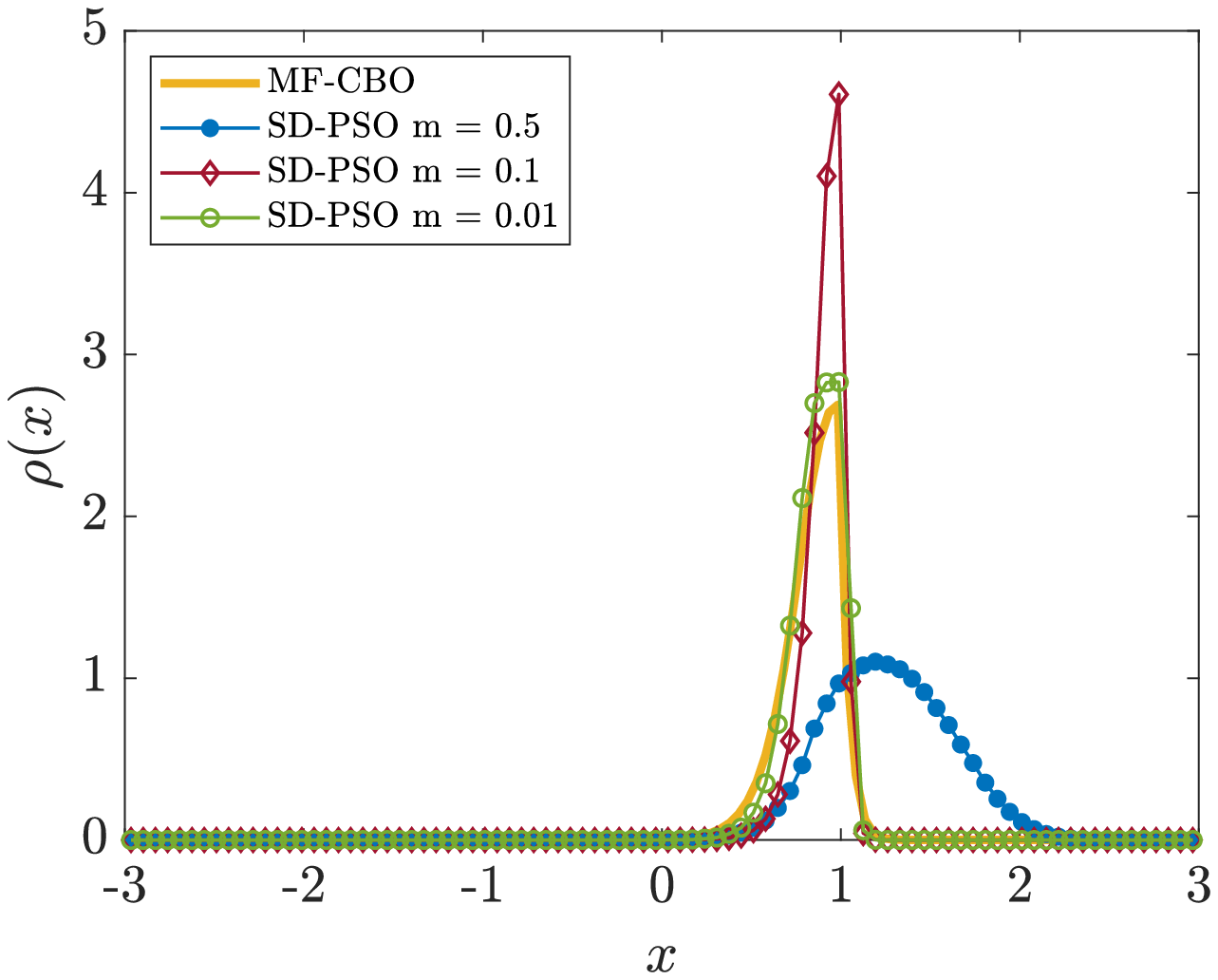}} \\
\subcaptionbox{$\rho(x,t)$, $t = 0.2$}{\includegraphics[scale=0.35]{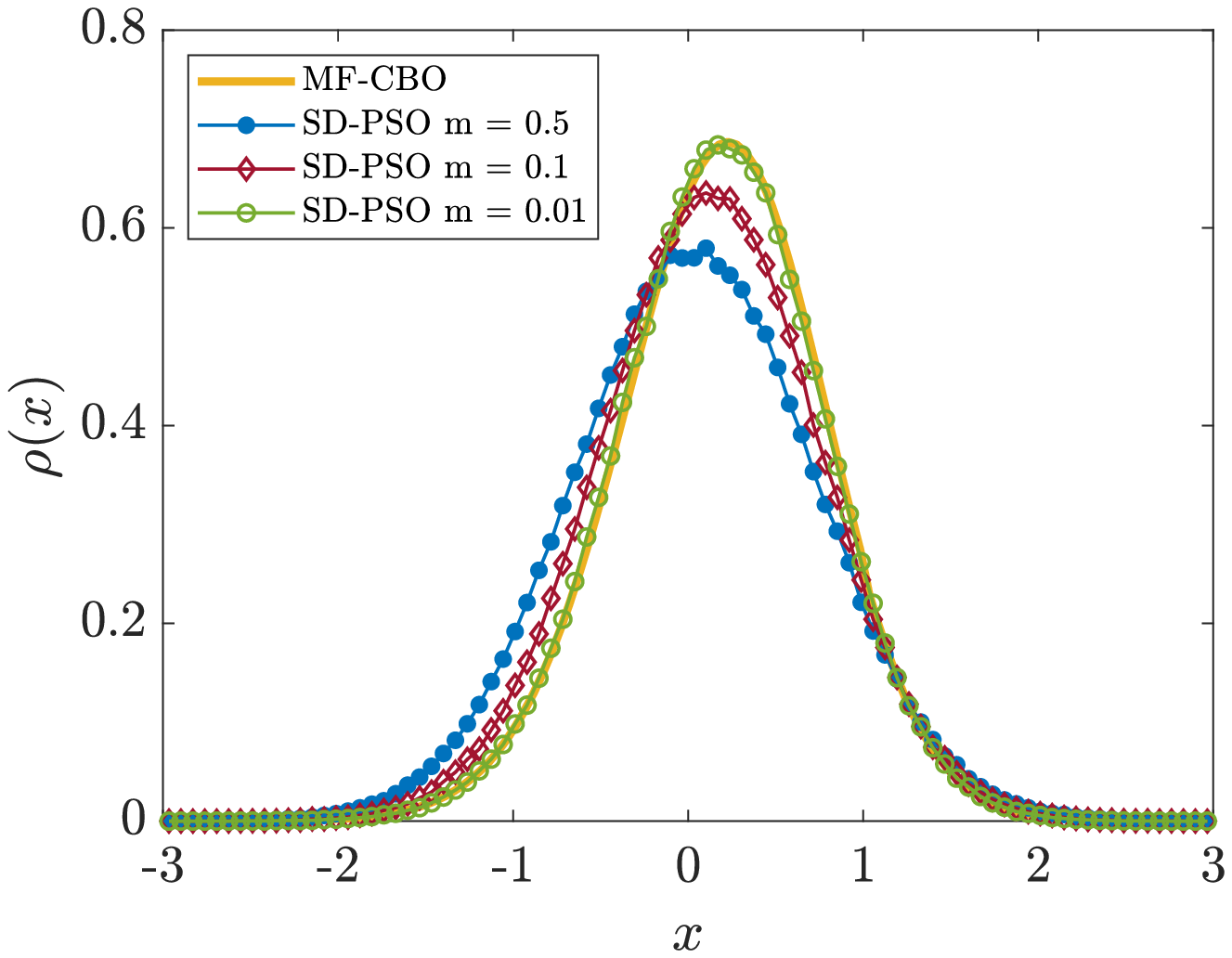}}\
\subcaptionbox{$\rho(x,t)$, $t = 1$}{\includegraphics[scale=0.35]{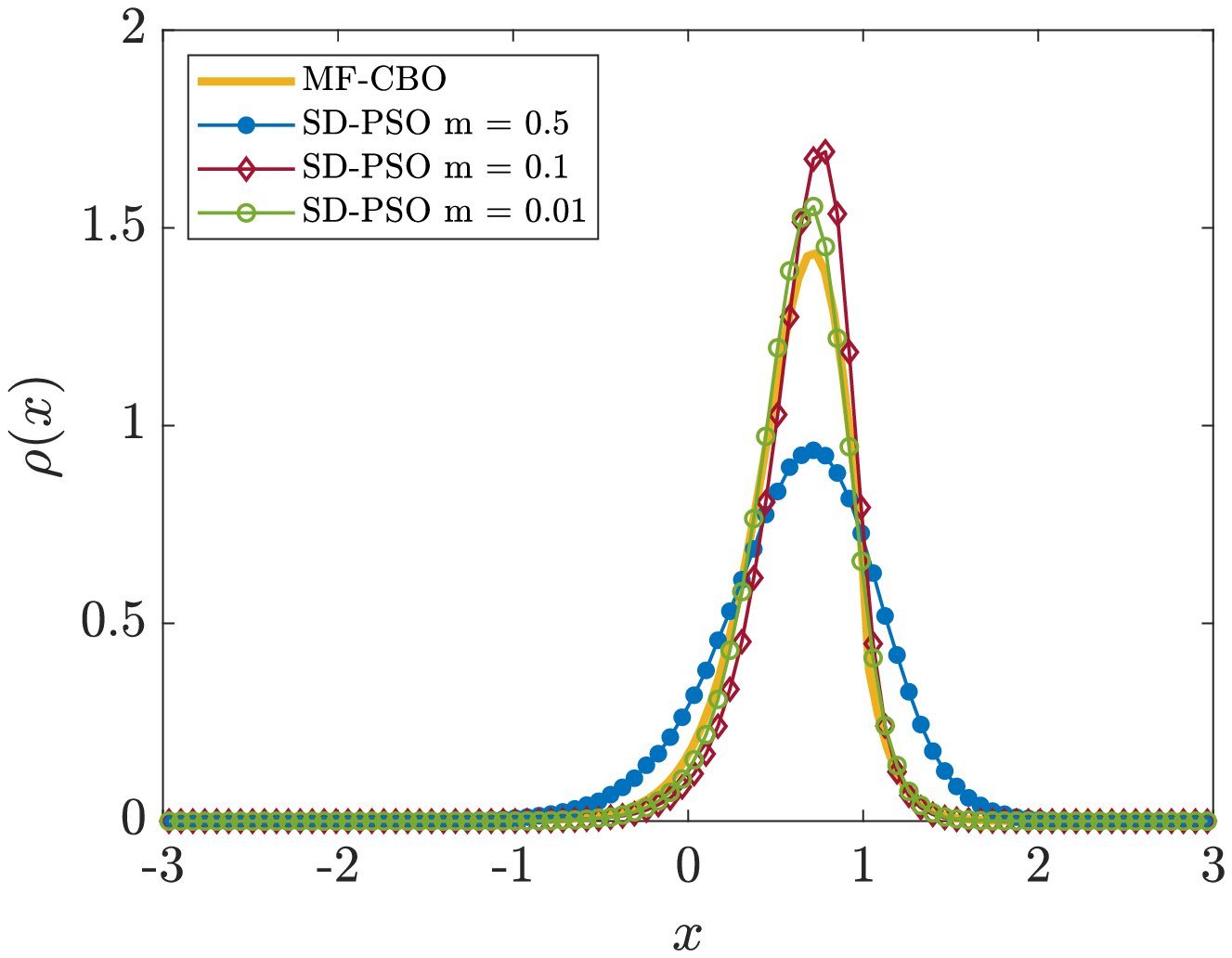}} \
\subcaptionbox{$\rho(x,t)$, $t = 2$}{\includegraphics[scale=0.35]{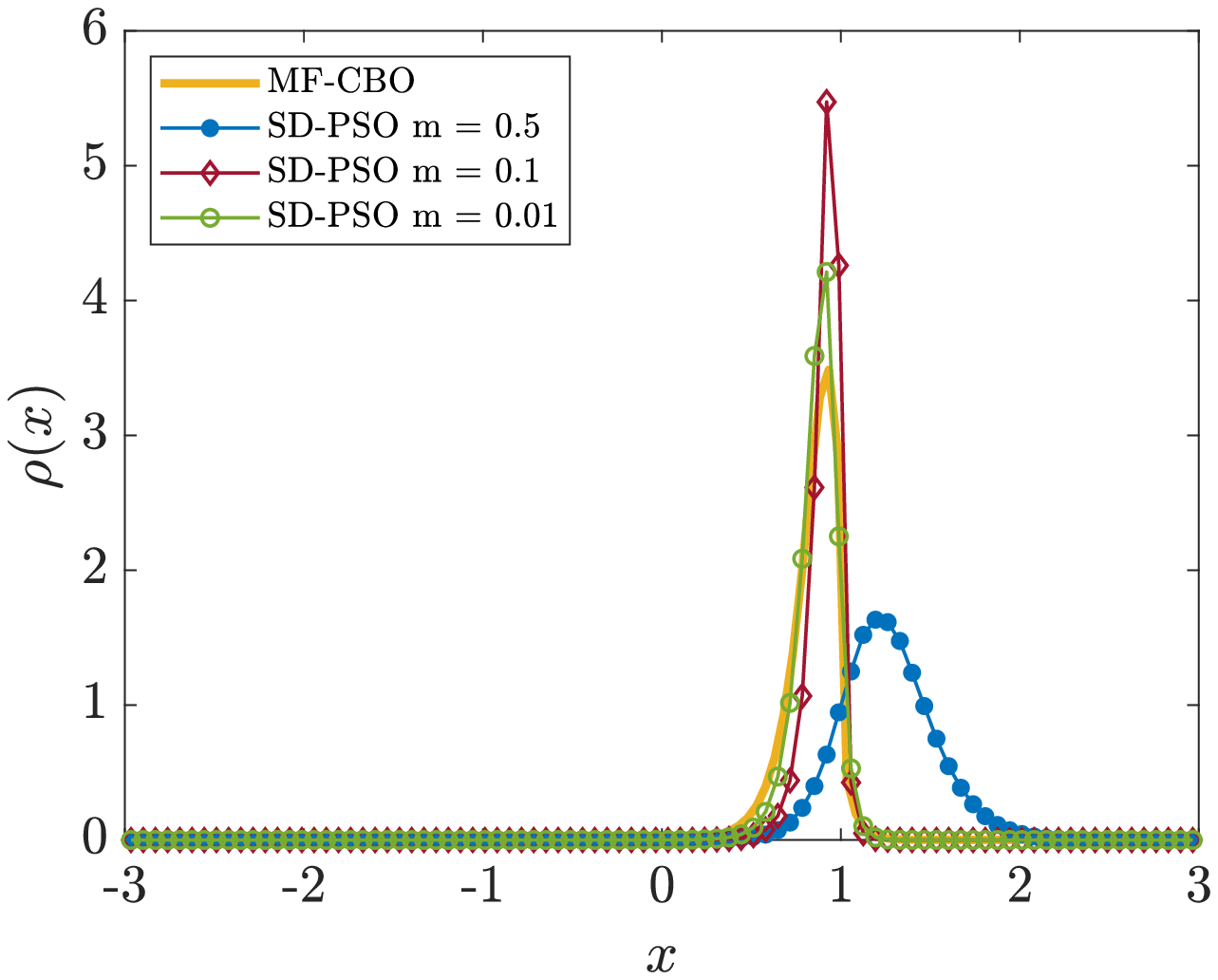}} 
\caption{Low inertia limit. Evolution of the density $\rho(x,t)$ of the SD-PSO discretization \eqref{eq:psoiDiscr2}, for decreasing inertial weight $m = 0.5, 0.1, 0.01$, and the mean-field CBO model \eqref{eq:CBOp} for the Ackley function with minimum in $x=1$. First row: uniform initial data. Second row: Gaussian initial data.}
\label{Fig16}
\end{minipage}
\end{figure}

First, let us observe that the semi-implicit discretization scheme \eqref{eq:psoiDiscr}  
\begin{eqnarray}
\nonumber
X^{n+1}_i &=& X^{n}_i + \Delta t \ V^{n+1}_i ,
\\[-.25cm]
\label{eq:psoiDiscr2}
\\[-.25cm]
\nonumber
V^{n+1}_i &=& \left( \frac{m}{m + \gamma \ \Delta t}\right) V^{n}_i + \frac{\lambda \ \Delta t }{m + \gamma \ \Delta t}\left(\Q_{\alpha}^n-X^n_i\right)+ \frac{\sigma \ \sqrt{\Delta t} }{m + \gamma \ \Delta t} D(\Q_{\alpha}^n-X^n_i) \ \theta^n_{i},
\end{eqnarray}
satisfies an asymptotic-preserving type property, allowing to pass to the limit $m\to 0$ without any restriction on $\Delta t$. In fact, passing to the limit, the second equation in \eqref{eq:psoiDiscr2} gives 
\[
V^{n+1}_i =  \lambda\left(\Q_{\alpha}^n-X^n_i\right)+ \frac{\sigma}{\sqrt{\Delta t}} D(\Q_{\alpha}^n-X^n_i) \ \theta^n_{i},
\]
which substituted into the first equation of \eqref{eq:psoiDiscr2} corresponds to the Euler-Maruyama scheme applied to the CBO system \eqref{eq:CBO}
\begin{eqnarray}
\label{eq:cboiDiscr2}
X^{n+1}_i &=& X^{n}_i + \Delta t  \lambda\left(\Q_{\alpha}^n-X^n_i\right)+ \sqrt{\Delta t}\sigma D(\Q_{\alpha}^n-X^n_i).
\end{eqnarray}
In Figures \ref{Fig15} and \ref{Fig16} we report the plots of the density that describes the solution of the mean-field CBO model and the stochastic PSO model for different inertial weights ($m = 0.5$, $m = 0.1$ and $m = 0.01$). We considered the minimization problem for the Ackley function with minimum in $x = 0$ and in $x = 1$ with $N=5 \times 10^{5}$ particles for the SD-PSO discretization, a grid of $120$ points in space for the mean field CBO solver and the same set of parameters \eqref{eq:paramt} with two different initial data: a uniform distribution and a Gaussian distribution. 

It is clear that in the case of $ m = 0.5 $ the two densities at the final time $ t = 2 $ are considerably different and a slower convergence is observed in the SD-PSO system, for $ m = 0.1 $ the agreement is higher and the particle solution seems to converge faster to the minimum, finally in the case $ m = 0.01 $ both densities simultaneously grow towards a Dirac delta centered in the minimum. As expected an initial Gaussian profile, being more concentrated, leads to a faster convergence. For smaller values of $m$ the two solutions becomes indistinguishable and we omitted the results.

\subsection{Comparison on high dimensional test cases}
In this section we report the results of several experiments concerning the behavior of the stochastic PSO models, discretized using \eqref{eq:psoiDiscr} in absence of memory or \eqref{eq:psoDiscr} in the general case, in high dimension ($d=20$) for various prototype test functions (see Appendix A). For the sake of simplicity we will focus our attention mostly to the case of the Ackley function and the Rastrigin function, used also in the previous examples, and report additional results for other global optimization test functions at the end of the Section. These two functions, in fact, although with several local minima presents very different levels of difficulty and have been used as test functions for CBO methods in various other papers \cite{pinnau2017consensus,carrillo2018analytical,carrillo2019consensus,fhps20-1,fhps20-2,TW}. 

In all tables reported in this section we will use the following terminology:
\begin{enumerate}
\item  the \textit{success rate}, computed averaging over $n_r=500$ runs and using as convergence criterion 
\[
\begin{split}
\Vert\Q_\alpha^{n_*}-x_{min}\Vert_{\infty}< \delta_{err},\qquad {\rm or}\qquad \Vert\bar\P_\alpha^{n_*}-x_{min}\Vert_{\infty}< \delta_{err} 
  \end{split}
  \label{ConvCrit}
\]
where $x_{min}$ is the position of the minimum, $n_*$ the final time, and $\delta_{err}=0.25$ as in \cite{pinnau2017consensus,carrillo2019consensus}.
 \item the \textit{error}, calculated as expected value in the $L_2$ norm over the successful runs
\[
\begin{split}
\mathbb{E}(\Vert \Q_\alpha^{n_*}-x_{min} \Vert_2 ),\qquad {\rm or}\qquad \mathbb{E}(\Vert \bar\P_\alpha^{n_*}-x_{min} \Vert_2 );
  \end{split}
\] 
\item the \textit{number of iterations}, where for a given tolerance $\delta_{stall}=10^{-4}$, we stop the iteration if
\[
\Vert \Q^{n}_{\alpha}-\Q^{n-1}_{\alpha}\Vert < \delta_{stall},\qquad {\rm or}\qquad \Vert \bar\P^{n}_{\alpha}-\Q^{n-1}_{\alpha}\Vert < \delta_{stall}
\]
for $n_{stall}=250$ consecutive iterations or a maximum number of $ 10^4 $ iterations has been reached. 
\end{enumerate}

It is important to remark that, while in previous simple one-dimensional validation examples we have chosen a low $ \alpha $ value, increasing the dimension, a larger value of $ \alpha \gg 1 $ provides better performance \cite{pinnau2017consensus}. On the other hand, a large value of $ \alpha $ may generate numerical instabilities given by the definition of $ \Q_{\alpha}^n $ in \eqref{ValphaE}. To avoid this, we used the algorithm presented in \cite{fhps20-2} which prevents the value of $N^{\alpha}$ from being close to $0$ by substituting, in absence of memory,
\begin{equation}
\frac{\omega_\alpha(X_i^n)}{N^\alpha} =\frac{\hbox{exp}(-\alpha(\EE (X_i^n)-\EE (X_{\ast}^n)))}{\sum_{i=1}^{N} \hbox{exp}(-\alpha(\EE (X_i^n)-\EE (X_{\ast}^n)))},
\label{alphalarge}
\end{equation}
with $X_{\ast}^n := \hbox{argmin}_{X \in \lbrace X_i^n \rbrace^N_{i = 1}} \EE(X)$. In a similar way, the method extends to the evaluation of $\bar\P_{\alpha}^n$ in the general case including memory.

Let us emphasize that, it is beyond the scopes of the present manuscript to perform an extensive testing of the various methods performances and to discuss the additional practical algorithmic enhancements that can be adopted to increase the success rate, like for example the use of random batch methods \cite{AlPa, carrillo2019consensus, JLJ}, particle reduction techniques \cite{fhps20-2} and parameters adaptivity \cite{poli2007particle}. In contrast, in the following test cases, we will address the role of the various parameters, of the presence of memory and of the local best when solving high dimensional global optimization problems. We refer also to \cite{Totzeck2018ANC} for additional comparisons. In all simulations we used Gaussian samples $\mathcal{N}(0,1)$ in the exploration term but in principle it is possible to use another distribution having mean 0 and variance 1. We analyzed also the use of a uniform noise $\mathcal{U}(-\sqrt{3},\sqrt{3})$ similarly to the original discrete PSO dynamic but without noticing significant changes in the results except that smaller values of $\sigma$ are needed to achieve convergence. In our experiments, additionally, the PSO constraints \eqref{eq:param} have shown strong limitations in terms of success rates and have not been considered. This has been verified both by direct simulations using traditional PSO parameters in our numerical scheme and by using standard PSO libraries such as the {\tt particleswarm} function of the Matlab Global Optimization Toolbox. 
 
\subsubsection*{Effect of the inertial parameter $m$ }
Initially we report in Table \ref{tab:1R} and Table \ref{tab:1A} the results obtained with the formulation  \eqref{eq:psoiDiscr} that does not exploit the memory of particles. The initial data is a uniform distribution on the whole domain of definition that here is fixed to $[-3, 3]^d$. Since, typically, optimizing the Rastrigin function is far more difficult than the Ackley function, we explore the space of parameters searching for optimal values of $\sigma$ and $\Delta t$ for the Rastrigin function,  then we used the same values for the Ackley function. This optimization has been done roughly through simple variations of a given step size for the parameters.
In the numerical examples we consider $\alpha = 5 \times 10$ and $\alpha = 5 \times 10^4$ to emphasize the role of such parameter in improving both the convergence rate as well as the efficiency of the solver. 

The structure of the tables is the same: the results are given for two different choices of $ \alpha $ and for different numbers of particles $N$. Since we are using a different solver we computed a different optimal value of $\sigma$ for the Rastrigin function which is then used also for the Ackley function. 
The other two important parameters to be set are $\beta$ and $\nu$, chosen respectively equal to $3\times 10^3$ and ${1}/{(2\Delta t)}=50$ so that the memory dynamics is close enough to the classical one of discrete PSO methods.
In this case, thanks to the memory effect there is no need to bound the computational domain since particle tend to converge to the global minimum without spreading in the whole space. Introducing the boundary conditions, the situation improves but not significantly, and we omit the corresponding results. Finally, from Table \ref{tb:2R} it is clear that, within the memory formulation, the choice of $ \alpha \gg 1 $ is essential to achieve good performances and therefore, the use of the algorithm described in \eqref{alphalarge} is fundamental in the practical implementation of the method. 

\subsubsection*{Introduction of local best dynamics}
Next we have considered the formulation \eqref{eq:psoDiscr} introducing the dynamics that lead the particles to move towards their local historical minimum. To reduce the number of free parameters we made the following assumption between the parameter defining the local best and the global best
\begin{align}
\lambda_1 = \xi \cdot \lambda_2,\qquad 
\sigma_1 = \xi \cdot \sigma_2
\label{eq:assum}
  \end{align}
with $\xi \in [0,1]$ so that the local best is always weighted less than the global best. In this test we keep the inertial value $m = 0$ and $\lambda_1=1$, so that we are solving the generalized stochastic differential CBO model with memory \eqref{eq:psocirlb}. For each value of $\xi$ reported, we have computed an optimal $\sigma_2$ achieving the maximum rate of success. We chose $\beta = 3\times 10^3$, $\Delta t=0.01$ and $\nu=0.5/\Delta t=50$ as in the previous case and consider $\alpha = 5 \times 10^4$ in evaluating the global best. 

\begin{table}[H]
\begin{center}
\begin{tabular}{l||l|l|ccc|l|ccc}
\multicolumn{3}{l}{\hspace{-5pt}\textbf{Rastrigin}} &  \multicolumn{4}{c}{\hspace{-45pt}Case $\alpha = 5 \times 10$} &\multicolumn{3}{c}{Case $\alpha = 5 \times 10^4$} \\
\hline
\hline
\hspace{5pt}$m$& & \hspace{4.5pt}$\sigma$ & $N=50$ & $N=100$ & $N=200$ & \hspace{4pt}$\sigma$ & $N=50$ & $N=100$ & $N=200$\\                  
\hline
{\rm $0.00$}& Rate & {\rm $7.0$}&  100.0\% & 100.0\% & 100.0\% &{\rm $9.0$}& 100.0\% & 100.0\% & 100.0\%\\
& Error &&  6.10e-04 & 3.91e-04 & 2.52e-04& & 1.19e-04 & 1.11e-04 & 9.68e-05\\
& $n_{iter}$& &10000.0 & 10000.0 & 10000.0 &&10000.0 & 10000.0 & 9912.4\\
\hline
{\rm $0.01$}& Rate &{\rm $6.5$} &  100.0\% & 100.0\% & 100.0\% &{\rm $7.0$}& 100.0\% & 100.0\% & 100.0\%\\
& Error &&   8.57e-04 & 4.94e-04 & 3.08e-04& & 9.74e-05 & 2.01e-05 & 1.62e-05\\
& $n_{iter}$ & &10000.0 & 10000.0 & 10000.0 &&10000.0 & 6899.2 & 2060.1\\
\hline
{\rm $0.05$}&Rate& {\rm $3.5$} &  42.5\% & 85.0\% & 92.0\% &{\rm $3.5$}& 37.0\% & 74.0\% & 94.0\%\%\\
& Error && 1.02e-03 & 7.98e-04 & 6.40e-04& &4.27e-04 & 1.26e-04 & 1.14e-04\\
& $n_{iter}$& &10000.0 & 10000.0 & 10000.0 &&8233.2 & 7814.0 & 7326.6\\
\hline
{\rm $0.10$}& Rate &{\rm $2.0$} &  0.7\% & 2.5\% & 12.5\% &{\rm $ 2.0$}& 1.0\% & 5.5\% & 29.5\%\\
& Error& & 3.52e-03 & 2.24e-03 & 2.05e-03& & 2.00e-04 & 1.28e-04 & 1.11e-04\\
& $n_{iter}$& &6818.3 & 7495.8 & 8680.9 &&6155.4 & 6221.9 & 6214.3 \\
\hline
\hline
\end{tabular}
\end{center}
\caption{SD-PSO without memory and with b.c. for $\lambda=1$ and $  \Delta t = 0.01$.}
\label{tab:1R}
\begin{center}
\begin{tabular}{l||l|l|ccc|l|ccc}
\multicolumn{3}{l}{\textbf{Ackley}} &  \multicolumn{4}{c}{\hspace{-30pt}Case $\alpha = 5 \times 10$} &\multicolumn{3}{c}{\hspace{5pt}Case $\alpha = 5 \times 10^4$} \\
\hline
\hline
\hspace{5pt}$m$& & \hspace{4.5pt}$\sigma$ & $N=50$ & $N=100$ & $N=200$ & \hspace{4pt}$\sigma$ & $N=50$ & $N=100$ & $N=200$\\          
\hline
{\rm $0.00$}& Rate & {\rm $7.0$} &  100.0\% & 100.0\% & 100.0\% &{\rm $ 9.0$}& 100.0\% & 100.0\% & 100.0\%\\
& Error & & 3.43e-03 & 1.90e-03 & 1.18e-03& & 8.46e-05 & 4.20e-05 & 1.27e-05\\
& $n_{iter}$& &10000.0 & 10000.0 & 10000.0&&1364.9 & 1032.4 & 869.2\\
\hline
{\rm $0.01$}& Rate& {\rm $ 6.5$} &  100.0\% & 100.0\% & 100.0\% &{\rm $ 7.0$}& 100.0\% & 100.0\% &  100\%\\
& Error && 5.03e-03 & 2.52e-03 & 1.36e-03& &9.49e-05 & 5.89e-05 & 2.81e-05\\
& $n_{iter}$& &10000.0 & 10000.0 & 10000.0 &&2192.9 & 1886.7 & 1723.6\\
\hline
{\rm $0.05$}& Rate&{\rm $ 3.5$} &  100.0\% & 100.0\% & 100.0\% &{\rm $ 3.5$}& 100.0\% & 100.0\% & 100.0\% \\
& Error && 3.76e-03 & 2.82e-03 & 7.74e-03& &  2.27e-04 & 1.48e-04 & 1.03e-04\\
& $n_{iter}$& & 10000.0 & 9857.2 & 6031.1 &&5367.3 & 4459.4 & 3928.4\\
\hline
{\rm $ 0.10$}&Rate& {\rm $2.0$} &  99.5\% & 100.0\% & 100.0\% &{\rm $ 2.0$}& 100.0\% & 100.0\% & 100.0\%\\
&Error& & 2.34e-03 & 2.28e-03 & 2.24e-03& &8.31e-04 & 2.76e-04 & 1.91e-04\\
&$n_{iter}$&  & 5914.0 & 3856.7 & 2909.1 &&5480.8 & 4514.1 & 3909.4\\
\hline
\hline
\end{tabular}
\end{center}
\caption{SD-PSO without memory and with b.c. for $\lambda=1$ and $  \Delta t = 0.01$.}
\label{tab:1A}
\end{table}
First, in Tables \ref{tab:3R} and \ref{tab:3A} we report the behavior of the particle optimizer on the Ackley and Rastrigin functions for different positions of the minimum $x_{min} = 0$, $x_{min} = 1$ and $x_{min} = 2$. In this test we need to use the boundary conditions in order to achieve a high success rate when the minimum is close to the boundary of the domain. Since for large values of $\xi$ we must decrease $\sigma_2$ to avoid a reduction of the convergence rate we expect that the total number of iterations may decrease. This is the case of the Ackley function in Table \ref{tab:3A} where a considerable speed-up is obtained thanks to the local best when the minimum is close to the boundary. 

\begin{table}[H]
\begin{center}
\begin{tabular}{l||l|l|ccc|l|ccc}
\multicolumn{3}{l}{\hspace{-5pt}\textbf{Rastrigin}} &  \multicolumn{4}{c}{\hspace{-45pt}Case $\alpha = 5 \times 10$} &\multicolumn{3}{c}{Case $\alpha = 5 \times 10^4$} \\
\hline
\hline
\hspace{5pt}$m$& & \hspace{5pt}$\sigma_2$ & $N=50$ & $N=100$ & $N=200$ & \hspace{5pt}$\sigma_2$ & $N=50$ & $N=100$ & $N=200$\\                  
\hline
{\rm $0.00$}& Rate & {\rm $11.0$}&  18.8\% & 16.8\% & 20.0\% &{\rm $11.0$}& 100.0\% & 100.0\% & 100.0\%\\
& Error &&  1.30e-03 & 5.09e-03 & 7.31e-03& & 6.83e-04 & 4.70e-04 & 4.69e-04\\
& $n_{iter}$& & 2331.9 & 1382.5 & 1289.6 && 10000.0 & 9878.2 & 3290.2\\
\hline
{\rm $0.01$}& Rate & {\rm \ $9.0$}&  25.4\% & 25.6\% & 39.5\% &{\rm \ $9.0$}& 100.0\% & 100.0\% & 100.0\%\\
& Error &&  4.53e-03 & 7.44e-03 & 9.06e-03& & 8.60e-04 & 8.56e-04 & 8.81e-04\\
& $n_{iter}$& &3536.0 & 3016.4 & 3128.9 &&9939.5 & 7012.2 & 5422.1\\
\hline
{\rm $0.05$}& Rate & {\rm \ $4.5$}&  30.4\% & 34.8\% & 44.4\% &{\rm \ $4.5$}& 100.0\% & 100.0\% & 100.0\%\\
& Error &&  4.51e-03 & 5.74e-03 & 9.87e-03& &  1.15e-03 & 6.67e-04 & 6.54e-04\\
& $n_{iter}$& &4646.0 & 4277.9 & 3598.5 &&9978.0 & 7657.6 & 5639.7\\
\hline
{\rm $0.10$}& Rate & {\rm \ $3.0$} &  8.6\% & 20.8\% & 35.2\% &{\rm \ $ 3.0$}& 80.8\% & 96.8\% & 100.0\%\\
& Error & & 2.72e-02 & 1.71e-02 & 1.31e-02& & 2.94e-03 & 8.96e-04 & 8.24e-04\\
& $n_{iter}$& &3686.7 & 4577.5 & 5361.4 &&9661.5 & 8676.5 & 7331.8\\
\hline
\hline
\end{tabular}
\end{center}
\caption{SD-PSO with memory for $\lambda_1=\sigma_1=0$, $\lambda_2=1$, $\Delta t = 0.01$, $\nu=50$, $\beta = 3 \times 10^3$.}
\label{tb:2R}
\begin{center}
\begin{tabular}{l||l|l|ccc|l|ccc}
\multicolumn{3}{l}{\textbf{Ackley}} &  \multicolumn{4}{c}{\hspace{-30pt}Case $\alpha = 5 \times 10$} &\multicolumn{3}{c}{\hspace{5pt}Case $\alpha = 5 \times 10^4$} \\
\hline
\hline
\hspace{5pt}$m$& & \hspace{5pt}$\sigma_2$ & $N=50$ & $N=100$ & $N=200$ & \hspace{5pt}$\sigma_2$ & $N=50$ & $N=100$ & $N=200$\\          
\hline
{\rm $0.00$}& Rate & {\rm $11.0$} &  100.0\% & 100.0\% & 100.0\% &{\rm $ 11.0$}& 100.0\% & 100.0\% & 100.0\%\\
& Error & & 2.84e-03 & 3.96e-03 & 5.47e-03& &1.02e-04 & 7.66e-05 & 5.44e-05\\
& $n_{iter}$& &2260.2 & 1762.0 & 1346.5 && 2457.0 & 1778.0 & 1513.1\\
\hline
{\rm $0.01$}& Rate & {\rm \ $9.0$} &  100.0\% & 100.0\% & 100.0\% &{\rm\  $ 9.0$}& 100.0\% & 100.0\% & 100.0\%\\
& Error & & 3.34e-03 & 4.70e-03 & 6.46e-03& &2.34e-03 & 1.91e-04 & 1.61e-04\\
& $n_{iter}$& &3722.2 & 2809.9 & 2104.2 && 6430.4 & 5447.8 & 4598.3\\
\hline
{\rm $0.05$}& Rate & {\rm \ $4.5$} &  100.0\% & 100.0\% & 100.0\% &{\rm \ $ 4.5$}& 100.0\% & 100.0\% & 100.0\%\\
& Error & & 4.17e-03 & 6.06e-03 & 8.24e-03& &2.41e-04 & 1.84e-04 & 1.48e-04\\
& $n_{iter}$& &5300.6 & 4059.0 & 3113.1 &&7186.1 & 5996.0 & 5074.6\\
\hline
{\rm $0.10$}& Rate & {\rm \ $3.0$} &  100.0\% & 100.0\% & 100.0\% &{\rm \ $ 3.0$}& 100.0\% & 100.0\% & 100.0\%\\
& Error & &6.72e-03 & 9.58e-03 & 1.25e-02& &3.90e-03 & 2.64e-03 & 2.06e-03\\
& $n_{iter}$& &6411.8 & 4856.5 & 3783.5 &&8590.6 & 7326.4 & 6350.2\\
\hline
\hline
\end{tabular}
\end{center}
\caption{SD-PSO with memory for $\lambda_1=\sigma_1=0$, $\lambda_2=1$, $\Delta t = 0.01$, $\nu=50$, $\beta = 3 \times 10^3$.}
\label{tb:2A}
\end{table}
Finally, in Table \ref{tab:TestFunctions} we report the results obtained by solving simultaneously a set of different optimization functions considered in their standard domains (see Appendix A). Here, instead of trying to find an optimal set of parameters for each function we use the same parameters for all functions. Even if further investigations are necessary in terms of identifying optimal set of parameters, through the previous simplifications assumptions we reduced our analysis to a minimal choice of parameters which seems the more relevant for the success of the algorithm. 

\begin{table}[H]
\begin{center}
{\renewcommand{\arraystretch}{1.1}
\begin{tabular}{l|lccc|ccc}
{\vspace{-2pt}\textbf{Rastrigin}} &  \multicolumn{4}{c}{Case $\xi = 0$, $\sigma_2 = 11.0$} &\multicolumn{3}{c}{Case $\xi = 0.25$, $\sigma_2 = 8.5$} \\
\hline
\hline
& & $N=50$ & $N=100$ & $N=200$ & $N=50$ & $N=100$ & $N=200$\\  
\hline
{\rm $x_{min} = 0$} & Rate  &  100.0\% & 100.0\% & 100.0\% & 100.0\% & 100.0\% & 100.0\%\\
 & Error & 7.04e-04 & 4.58e-04 & 3.29e-04&9.28e-04 & 6.11e-04 & 4.31e-04\\
& $n_{iter}$ & 10000.0 & 9963.9 & 4635.1&9978.0 & 8311.5 & 5754.1\\
\hline
{\rm $x_{min} = 1$}& Rate  &  98.8\% & 100.0\% & 100.0\% & 99.2\% & 100.0\% & 100.0\%\\
 & Error & 7.08e-04 & 4.60e-04 & 3.27e-04&9.31e-04 & 6.74e-04 & 4.59e-04\\
& $n_{iter}$ &10000.0 & 10000.0 & 4670.0 & 9987.0 & 9746.7 & 7460.1\\
\hline
{\rm $x_{min} = 2$}& Rate  &96.0\% & 99.1\% & 100.0\% & 93.5\% & 100.0\% & 100.0\%\\
 & Error &  6.91e-04 & 4.52e-04 & 3.28e-04&8.78e-04 & 6.74e-04 & 5.66e-04\\
& $n_{iter}$ &10000.0 & 10000.0 & 5035.5 & 9980.3 & 9854.1 & 8971.9\\
\hline  
\hline       
\end{tabular}}
\end{center}
\caption{SD-PSO with memory ($m=0$) for $\lambda_1$ and $\sigma_1$ given by \eqref{eq:assum}, $\lambda_2=1$, $\Delta t = 0.01$, $\nu=50$, $\beta = 3 \times 10^3$, $\alpha = 5 \times 10^4$.}
\label{tab:3R}
\begin{center}
{\renewcommand{\arraystretch}{1.1}
\begin{tabular}{l|lccc|ccc}
{\textbf{Ackley}} &  \multicolumn{4}{c}{Case $\xi = 0$, $\sigma_2 = 11.0$} &\multicolumn{3}{c}{Case $\xi = 0.25$, $\sigma_2 = 8.5$} \\
\hline
\hline
& & $N=50$ & $N=100$ & $N=200$ & $N=50$ & $N=100$ & $N=200$\\  
\hline
{\rm $x_{min} = 0$ \ } & Rate  &  100.0\% & 100.0\% & 100.0\% & 100.0\% & 100.0\% & 100.0\%\\
 & Error & 7.36e-05 & 5.13e-05 & 3.26e-05&2.54e-05 & 1.13e-05 & 1.07e-05\\
& $n_{iter}$ &2778.6 & 2030.0 & 1623.0 & 1942.9 & 1663.8 & 1442.5\\
\hline
{\rm $x_{min} = 1$ \ }& Rate  &  100.0\% & 100.0\% & 100.0\% & 100.0\% & 100.0\% & 100.0\%\\
 & Error & 7.31e-05 & 5.14e-05 & 3.26e-05&2.58e-05 & 1.12e-05 & 1.02e-05\\
& $n_{iter}$ & 5298.5 & 3640.6 & 2575.9 &  2465.3 & 1948.5 & 1632.5\\
\hline
{\rm $x_{min} = 2$ \ }& Rate  & 100.0\% & 100.0\% & 100.0\% &  100.0\% & 100.0\% & 100.0\%\\
 & Error & 7.30e-05 & 5.07e-05 & 3.22e-05&2.64e-05 & 1.09e-05 & 1.01e-05\\
& $n_{iter}$ &7819.8 & 5771.3 & 4235.9 &3126.8 & 2286.0 & 1803.8\\
\hline  
\hline       
\end{tabular}}
\end{center}
\caption{SD-PSO with memory ($m=0$) for $\lambda_1$ and $\sigma_1$ given by \eqref{eq:assum}, $\lambda_2=1$, $\Delta t = 0.01$, $\nu=50$, $\beta = 3 \times 10^3$, $\alpha = 5 \times 10^4$.}
\label{tab:3A}
\end{table}
Thus, we let most parameters fixed as in previous test case, namely $\alpha = 5 \times 10^4$, $\beta = 3 \times 10^3$, $\nu=0.5/\Delta t$ since these essentially define the modeling process of the local best and global best. Additionally we keep $m=0$, $\Delta t = 0.01$, and for a given value of $\xi=0$ (absence of local best) and $\xi=0.25$ (local best weighted 1/4 of global best) estimate the value for $\sigma_2$ in order to maximize the average convergence rate among all functions. Again this has been done roughly with simple variations of step $0.5$ for $\sigma_2$ in the simulations.

The results confirm the potential of the method in identifying correctly the global minima for different test functions. Overall, with the exception of the Rastrigin function for which  the local best produces a reduction in the convergence rate using this set of parameters, the importance of the local best is evident. In particular, the presence of the local best yields a significant reduction in the number of iterations for the Salomon function and an increase in the convergence rate for the XSY random function.
 
\begin{table}[H]
\begin{center}
{\renewcommand{\arraystretch}{1.2}
\begin{tabular}{l|lccc|ccc}
 & \multicolumn{4}{c}{Case $\xi = 0$, $\sigma_2 = 8.0$} &\multicolumn{3}{c}{ \hspace{5pt} Case $\xi = 0.25$, $\sigma_2 = 6.5$ \hspace{5pt} } \\
\hline
\hline
& & $N=50$ & $N=100$ & $N=200$ & $N=50$ & $N=100$ & $N=200$ \\  
\hline
{\rm \textbf{Ackley}}& Rate  &100.0\% & 100.0\% & 100.0\%  &100.0\% & 100.0\% & 100.0\% \\
& Error &  7.74e-05 & 6.12e-05 & 5.22e-05 & 2.08e-04 & 1.98e-04 & 1.65e-04\\
& $n_{iter}$ &1325.0 & 1114.8 & 924.9 &1263.2 & 992.3 & 902.4 \\
\hline
{\rm \textbf{Rastrigin}}& Rate  &31.4\% & 65.7\% & 95.6\% &5.3\% & 10.5\% & 27.9\% \\
& Error &  5.59e-04 & 7.70e-04 & 9.68e-04 &  1.46e-03 & 9.73e-04 & 1.06e-03\\
& $n_{iter}$ &1404.7 & 1107.8 & 954.1 &4390.0 & 4756.2 & 4643.6\\
\hline
{\rm \textbf{Griewalk}}& Rate  &100.0\% & 100.0\% & 100.0\% &100.0\% & 100.0\% & 100.0\% \\
& Error &  9.12e-02 & 7.55e-02 & 5.78e-02 &  8.22e-02 & 5.34e-02 & 4.12e-02\\
& $n_{iter}$ &10000.0 & 10000.0 & 10000.0 &10000.0 & 10000.0 & 9978.4 \\
\hline
{\rm \textbf{Schwefel}}& Rate  &99.6\% & 100.0\% & 100.0\% &100.0\% & 100.0\% & 100.0\% \\
& Error & 1.16e-05 & 1.22e-05 & 1.39e-05 &  2.69e-05 & 2.80e-05 & 2.85e-05\\
& $n_{iter}$ &1211.2 & 1044.1 & 987.3 &1127.9 & 964.8 & 849.7 \\
\hline
{\rm \textbf{Salomon}}& Rate  &96.7\% & 98.3\% & 100.0\% &100.0\% & 100.0\% & 100.0\% \\
& Error & 9.26e-02 & 8.73e-02 & 8.02e-02 & 8.22e-02 & 6.12e-02 & 5.23e-02\\
& $n_{iter}$ &9443.0 & 8176.2 & 7443.3 &6476.0 & 3009.7 & 1923.2\\
\hline
{\rm \textbf{XSY random \ }}& Rate  &35.5\% & 59.5\% & 94.2\% &75.1\% & 94.2\% & 100.0\% \\
& Error &  1.12e-01 & 9.81e-02 & 8.82e-02 &  1.06e-01 & 9.79e-02 & 8.57e-02\\
& $n_{iter}$ &10000.0 & 10000.0 & 10000.0 &10000.0 & 10000.0 & 10000.0\\
\hline  
\hline       
\end{tabular}}
\end{center}
\caption{SD-PSO with memory ($m=0$) for $\lambda_1$ and $\sigma_1$ given by \eqref{eq:assum}, $\lambda_2=1$, $\Delta t = 0.01$, $\nu=50$, $\beta = 3 \times 10^3$, $\alpha = 5 \times 10^4$.}
\label{tab:TestFunctions}
\end{table}

\section{Conclusions}
In this work we attempted to make a contribution to the construction of a general mathematical theory that will allow the rigorous analysis of optimization methods based on particle swarms (PSO). To this aim, starting from the original discrete formulation, we derived, by approximating in an appropriate way the memory dynamics, the corresponding systems of SDEs. In the large particle limit, using a regularized version of these systems we obtain a mean-field PDE of Vlasov-Fokker-Planck type describing the PSO dynamic. The new mean-field formalism, for small values of the inertia parameter, permits to compute as hydrodynamic approximation a generalization of consensus-based optimization models (CBO) with local best, thus highlighting the relationships between these two classes of metaheuristic optimization methods. These results are numerically validated through several examples that compare the mean-field and the particle dynamics. The methods are then tested against some prototype high dimensional global optimization functions with the goal of understanding the effects of the various parameters and the main differences between the novel stochastic differential models and standard CBO systems, namely the presence of memory effects together with the local best. The numerical results confirmed the ability of the local best to improve the performance of the methods in terms of speed of convergence and rate of success. It is worth noting that the SD-PSO models here introduced, thanks to the increased independence of the search parameters, allow better performances of the classic PSO methods (which represent a particular case included in the choices of possible optimization parameters). 

In perspective, the introduction of a memory variable opens interesting possibilities towards the construction of novel SD-PSO minimizers and CBO minimizers for multi-objective functions. In addition to this, the main research directions we intend to deal with in the near future concern the analysis of the convergence properties of the mean-field PSO system to global minimizers, the rigorous derivation of the small inertia limit, and the study of the convergence ranges for the parameters that characterize the system.

\section*{Acknowledgment}
The authors are grateful to Y-P.~Choi for helpful discussions concerning the overdamped limit of Vlasov-Fokker-Planck equations. This work has been written within the
activities of GNCS group of INdAM (National Institute of
High Mathematics). The support of MIUR-PRIN Project 2017, No. 2017KKJP4X "Innovative numerical methods for evolutionary partial differential equations and applications" and of the ESF PhD grant "Mathematical and statistical methods for machine learning in biomedical and socio-sanitary applications" is acknowledged.

\appendix
\section{Test Functions}
In this Appendix we report the global optimization test functions that were used in the article to validate the performance of the SD-PSO algorithms \cite{JYZ}. Each function is reported within the typical search domain and, if $B,C = 0$, has the global minima located at $x^{\ast} = (0,\dots,0)$ with $\EE(x^{\ast}) = 0$.
\begin{enumerate}
\item \textbf{Ackley function} (Continuous, Differentiable, Non-convex, Non-Separable, Multimodal)
\begin{align}
\EE(x) = -20\ \mbox{exp}\left( -0.2\sqrt{\frac{1}{d}\sum_{i=1}^{d}{(x_i-B)^2}}\right)-\mbox{exp}\left(\frac{1}{d}\sum_{i=1}^d{ \cos \left( 2\pi (x_i-B)\right) }\right)+20 +\mbox{exp}(1)+C,
\label{eq:Ackley}
\end{align}
subject to $-32 \leq x \leq 32$. 
\item  \textbf{Griewalk function} (Continuous, Differentiable, Non-convex, Separable, Unimodal)
\begin{align}
\EE(x) = 1+\sum_{i=1}^{d} \frac{(x_i-B)^2}{4000}-\prod_{i=1}^{d} \mbox{cos}\left(\frac{x_i-B}{i}\right) +C,
\label{eq:Griewalk}
\end{align}
subject to $-100 \leq x \leq 100$. 
\item  \textbf{Rastrigin function} (Continuous, Differentiable, Convex, Separable, Multimodal)
\begin{align}
\EE(x) = 10d + \sum_{i=1}^d \left[(x_i-B)^2 - 10 \cos\left(2\pi (x_i-B) \right)\right]+C,
\label{eq:Rastrigin}
\end{align}
subject to $-5.12 \leq x \leq 5.12$. .
\item  \textbf{Salomon function} (Continuous, Differentiable, Non-Convex, Non-Separable, Multimodal)
\begin{align}
\EE(x) =  1-\cos{\left( 2\pi \sqrt{\sum_{i=1}^{d}{(x_i-B)^2}}\right)} + 0.1 \sqrt{\sum_{i=1}^{d}{(x_i-B)^2}}  +C,
\label{eq:Salomon}
\end{align}
subject to $-100 \leq x \leq 100$. 
\item  \textbf{Schwefel function} (Continuous, Non-Differentiable, Convex, Separable, Unimodal)
\begin{align}
\EE(x) =  \sum_{i = 1}^d \vert x_i-B \vert +C ,
\label{eq:Schwefel}
\end{align}
subject to $-100 \leq x \leq 100$. 
\item  \textbf{Xin She-Yang random function} (Random, Non-Differentiable, Non-convex, Separable, Multimodal)
\begin{align}
\EE(x) = \sum_{i=1}^{d} \eta_i \vert x_i - B \vert^{i} + C
  \label{eq:XSY random}
\end{align}
with $\eta_i, i = 1,\dots,d$ random variable uniformly distributed in $[0,1]$. The standard domain size is $-5 \leq x \leq 5$. 
\end{enumerate}

\end{document}